\documentclass[x11names]{article}

\usepackage{graphicx}
\usepackage{amssymb}
\usepackage{bbm}
\usepackage[pdftex]{hyperref}
\usepackage{cite}
\usepackage[hmargin={2.5cm,2.5cm}, vmargin={3cm,3cm}, dvips]{geometry}

\usepackage{fancyhdr}
\usepackage[english]{babel} % definisce la lingua del documento
\usepackage{amsmath} %Formules
\usepackage{amssymb} %Formules
\usepackage{amsthm} %Environnements
\usepackage{bm}
\usepackage{array}
\usepackage{placeins}
\usepackage{enumerate}
\usepackage{esint}
\usepackage{graphics}
\usepackage{graphicx}
\usepackage{graphbox}
\usepackage{tikz}
\usetikzlibrary{patterns}
\usetikzlibrary{perspective}
\usetikzlibrary{3d}
\usepackage{pgfplots}
\usepackage{subcaption} % serve per aggiungere la didascalia alle figure composte (carica anche CAPTION)
\usepackage{listings}
\usepackage[shortlabels]{enumitem}
\usepackage{multirow}
\usepackage{multicol}
\usepackage{wrapfig}
\usepackage{algorithm} % Per scrivere gli algoritmi utilizzati
\usepackage[noend]{algpseudocode}
\usepackage{eqparbox}

\usepackage{afterpage}
\usepackage{mathrsfs}
\usepackage{mathtools}

\newtheoremstyle{droit}
{}%Space above
{}%Space below
{\upshape}%Body font
{}%Indent amount (empty = no indent, \parindent = para indent)
{\bfseries}%theoreme head font
{}%Punctuation after theoreme head
{ }%{\newline}%Space after theoreme head: " " = normal interword space; \newline = linebreak
{}%theoreme head spec (can be left empty, meaning `normal')
\newtheoremstyle{italique}
{}%Space above
{}%Space below
{\itshape}%Body font
{}%Indent amount (empty = no indent, \parindent = para indent)
{\bfseries}%theoreme head font
{}%Punctuation after theoreme head
{ }%{\newline}%Space after theoreme head: " " = normal interword space; \newline = linebreak
{}%theoreme head spec (can be left empty, meaning `normal')
\theoremstyle{italique}
\newtheorem{theorem}{Theorem}[section]

\theoremstyle{droit}

\newtheorem{definition}[theorem]{Definition}
\newtheorem{assumption}[theorem]{Assumption}

\usepackage{geometry}
\usetikzlibrary{shapes,backgrounds}

\newcommand{\R}{\ensuremath{\mathbb R}}%Real numbers
\DeclarePairedDelimiter\abs{\lvert}{\rvert}%Absolute value
\newcommand{\spn}[1]{\ensuremath{\mathrm{span}\left\{#1\right\}}}
\newcommand{\supp}[1]{\ensuremath{\mathrm{supp}\left(#1\right)}}

\hypersetup{
	colorlinks=true,       % false: boxed links; true: colored links
	linkcolor=blue,        % color of internal links
	citecolor=SpringGreen4,        % color of links to bibliography
	filecolor=magenta,     % color of file links
	urlcolor=blue         
}

\begin{document}

%\begin{frontmatter}

%% Title, authors and addresses

\title{\textbf{A shape derivative approach to domain simplification}}
%% use the tnoteref command within \title for footnotes;
%% use the tnotetext command for the associated footnote;
%% use the fnref command within \author or \address for footnotes;
%% use the fntext command for the associated footnote;
%% use the corref command within \author for corresponding author footnotes;
%% use the cortext command for the associated footnote;
%% use the ead command for the email address,
%% and the form \ead[url] for the home page:
%%
%% \title{Title\tnoteref{label1}}
%% \tnotetext[label1]{}
%% \author{Name\corref{cor1}\fnref{label2}}
%% \ead{email address}
%% \ead[url]{home page}
%% \fntext[label2]{}
%% \cortext[cor1]{}
%% \address{Address\fnref{label3}}
%% \fntext[label3]{}

%% use optional labels to link authors explicitly to addresses:
%% \author[label1,label2]{<author name>}
%% \address[label1]{<address>}
%% \address[label2]{<address>}

\author{{J. Hinz$^{1}$, O. Chanon$^2$, A. Arrigoni$^{1}$, A. Buffa$^{1,3}$}\\ \\
	\footnotesize{$^1$ MNS, Institute of Mathematics, \'Ecole Polytechnique F\'ed\'erale de Lausanne, Switzerland}\\
    \footnotesize{$^2$ Institute of Analysis and Scientific Computing, TU Wien, Austria}\\
	\footnotesize{$^3$ Istituto di Matematica Applicata e Tecnologie Informatiche `E. Magenes' (CNR), Pavia, Italy}
}

\maketitle
\vspace{-0.8cm}
\noindent\rule{\linewidth}{0.4pt}
\thispagestyle{fancy}
\begin{abstract}
The objective of this study is to address the difficulty of simplifying the geometric model in which a differential problem is formulated, also called defeaturing, while simultaneously ensuring that the accuracy of the solution is maintained under control. This enables faster and more efficient simulations, without sacrificing accuracy. More precisely, we consider an isogeometric discretisation of an elliptic model problem defined on a two-dimensional hierarchical B-spline computational domain with a complex boundary. Starting with an oversimplification of the geometry, we build a goal-oriented adaptive strategy that adaptively reintroduces continuous geometrical features in regions where the analysis suggests a large impact on the quantity of interest. This strategy is driven by an \textit{a posteriori} estimator of the defeaturing error based on first-order shape sensitivity analysis, and it profits from the local refinement properties of hierarchical B-splines. The adaptive algorithm is described together with a procedure to generate (partially) simplified hierarchical B-spline geometrical domains. Numerical experiments are presented to illustrate the proposed strategy and its limitations. 
\end{abstract}

\noindent\rule{\linewidth}{0.4pt}

%%%%%%%%%%%%%%%%%%%%%%%%%%%%%%%%%%%%%%%%%%%%%%%%%%%%%%%%%%%%%%%%%%%%%

\section{Introduction} \label{s:intro}

% \textcolor{red}{ we would like to enrich the space to create more and more complex geometries by selecting a set of functions to add, so that the error on the functional decreases. We define an estimator of the modelling error. We consider only changes in the boundary, not topological changes. \\
% Say what we mean by exact boundary (spline model or mathematical representation).\\
%DEFINE modelling error (change name?) and defeaturing problem/process (change name wrt Ondine?)
% Compared to Ondine, her defeaturing is local, while ours is global.} \\
%We tested only the 2D, but the algorithm can be extended to the 3D provided we have an efficient mesher\\}\\

% Domanda: nel problem setting metto solo neumann BC oppure anche dirichlet? In teoria solo Neumann, pensando al problema vero, ma poi nei test abbiamo anche Dirichlet.

The finite element method (FEM) has become an indispensable numerical technique for engineering research based on computational simulation. With the advancement of modern computer hardware as well as algorithmic improvements, FEM is employed in applications of increasing geometrical complexity. Mathematically, geometries are represented as two- or three-dimensional computer-aided design (CAD) models. The transition from a CAD model to a (typically simplicial) FEM-suitable computational mesh is associated with a number of accuracy and robustness bottlenecks \cite{cottrell2009isogeometric}. As such, isogeometric analysis (IGA) \cite{hughes2005isogeometric} has emerged as a promising variant of classical FEM since it employs the CAD-typical spline functions both to represent the geometry and as a basis for the numerical analysis of the underlying partial differential equations (PDEs). In IGA, the CAD model is immediately used for computational simulation, thereby avoiding the translation into a classical FEM-suitable mesh. \\

\noindent Unfortunately, increased geometrical complexity (in both standard FEM and IGA) is associated with novel challenges. In particular, a complex CAD model often causes the mesh generation to fail or to introduce an unfeasible number of elements. In FEM, for instance, White et al. \cite{white2003meshing} demonstrate that adding a single feature (such as a small cavity in the geometry's interior) can increase the required number of mesh elements by up to a factor of ten. Similarly in IGA, a complex CAD model of the geometry's boundary $\partial \Omega$ often leads to prohibitive computational costs. The reason for this is two-fold: firstly, one needs to solve the so-called surface-to-volume problem $\partial \Omega \rightarrow \Omega$ to find a valid description of the geometry's interior starting from the CAD input which is usually a boundary representation. Due to the (B-)splines' tensor-product structure, the increased complexity of the boundary contour's spline representation greatly increases the computational costs of the surface-to-volume problem. Secondly, requiring a denser spline basis to represent the geometry's interior $\Omega$ inevitably increases the computational costs of the numerical analysis step. The latter is particularly critical in settings where the solution's accuracy is only of importance on a subset of the whole domain. \\

\noindent To address this problem, one can use more advanced and potentially unstructured spline technologies to represent the geometry: for instance, the original IGA method has been extended to geometries obtained by trimming \cite{trimmingbletzinger,trimming, antolinvreps,weimarrusigantolin, antolin2021quadrature}, unions \cite{unionkargaran,unionzuo,antolinwei} or with multipatches \cite{multipatch,bracco2020isogeometric}, or one can use for instance (truncated) hierarchical B-splines (THB-splines) that allow for local refinement, see \cite{overview_thb} and the references therein. However, practical applications often require a more fundamental intervention. The concept considered in this work is called \textit{defeaturing} and concerns itself with the removal of geometric features without critically altering the PDE's solution in the regions of interest. In this context, features can be divided into two broad categories: 1) discrete features that can be regarded as performing a Boolean operation on a reference geometry (such as creating a hole in the domain's interior or adding / removing a discrete boundary feature) and 2) continuous features that manifest themselves in a higher local density of spline functions in the geometry's CAD model. The second case enables the examination of geometries where the boundary exhibits complex or multiscale characteristics everywhere.\\

\noindent Traditionally, defeaturing has relied on heuristical \textit{a priori} criteria often drawing on the engineer's practical experience \cite{thakur2009survey}. More deterministic \textit{a priori} indicators are based on constitutive or conservation laws \cite{fine2000automated, foucault2004mechanical, rahimi2018cad}. However, to integrate defeaturing into a fully automated workflow, \textit{a posteriori} estimators are of larger interest. In \cite{ferrandes2009posteriori}, Ferrandes et al. propose an \textit{a posteriori} error estimator based on approximating the energy norm between the exact and the defeatured geometries' solutions. A further technique is based on feature sensitivity analysis (FSA). FSA uses topological sensitivity analysis \cite{sokolowski1999topological, choi2004structural} to estimate the first order change in the quantity of interest due to the removal of small interior or boundary features. However, due to its first-order nature, this technique is limited to small (Boolean) features. In \cite{li2011estimating}, Li et al. propose an error estimator for internal holes by formulating two different PDEs, one corresponding to the original and one to the defeatured geometry, on a unique reference geometry. The methodology then employs the dual weighted residual \cite{becker2001optimal, oden2002estimation} \textit{a posteriori} technique to estimate the modelling error. The same technique is adopted in \cite{li2013engineering, li2013goal, zhang2016estimation} and generalised to various linear and quasi-linear model problems along with support for a wider range of features. However, all aforementioned techniques rely on some degree of heuristic and a precise investigation of the estimator's accuracy is lacking. To mitigate this issue, Buffa et al. introduce the concept of \textit{analysis-aware} defeaturing \cite{buffa2022analysis, antolin2022multifeature, buffa2022adaptive}. In this series of articles, the authors propose an \textit{a posteriori} error estimator for Boolean operations on isogeometric discretisations of elliptic model problems, and they build an adaptive strategy based on this estimator. The reliability of the estimator is proved for general geometric configurations, thus largely avoiding heuristics. \\

\noindent Compared to Boolean features, analysis-aware continuous defeaturing techniques are nevertheless under-represented in the literature, and to the best of the author's knowledge, this work is amongst the first ones to approach this challenge. As such, this manuscript presents a goal-oriented adaptive algorithm for the continuous defeaturing of planar geometries stemming from complex CAD models. We focus primarily on isogeometric discretisations of elliptic model problems. Starting with a gross oversimplification of the CAD input represented by THB-splines, the proposed methodology employs the concept of \textit{shape gradients} to adaptively reintroduce continuous features in regions where first-order sensitivity analysis suggests a large impact on a quantity of interest. This strategy has been inspired by the work of Buffa et al. \cite{buffa2021adaptiveapproxshapes}. For the isogeometric discretisation and since standard B-splines, which have a tensor-product structure, makes it difficult to achieve localized refinement, we employ THB-splines. This enables us to accurately capture important boundary features without introducing too many degrees of freedom in the domain's interior. The methodology does not require a parameterisation of the complex geometry's interior. Instead, only a parameterisation of the defeatured geometry is computed in each iteration using the concept of harmonic maps. Heydarov and al. \cite{heydarov2022unrefinement} have also used (T)HB-splines to simplify geometrical models; however in that article and in contrast to the present work, the defeaturing is not driven by the analysis of the PDE at hand. \\

This manuscript is structured as follows. In the remaining part of the introduction, we precisely describe the tackled problem and the main objective of this work, namely the construction of an adaptive algorithm for continuous defeaturing. Then in Section~\ref{sec:math_tools}, we introduce the necessary mathematical tools on which relies the proposed adaptive strategy, namely the notion of shape derivatives, THB-splines, and a parametrisation and mesh generation method. The core of the manuscript can then be found in Section~\ref{sec:algo}, in which the proposed adaptive gradient-based defeaturing algorithm is described in details, together with implementation considerations. Subsequently in Section~\ref{sec:numerical}, we present numerical experiments that illustrate the proposed method, its capabilities and its limitations. We finally draw some conclusions in Section~\ref{sec:conclusions}.

\subsection{Problem setting}
\label{subsec:pb_setting}

\noindent 
Let $\Omega \subset \R^2$ be an open bounded Lipschitz domain with a piecewise smooth, potentially complicated boundary. In this manuscript, we are interested in approximating a quantity of interest $J(\, \cdot \,, \Omega)$ computed over the solution $u \in H^1(\Omega)$ of an elliptic differential problem posed over $\Omega$. More precisely, let $\Gamma^N, \Gamma^D \subset \partial \Omega$ such that $\overline{\Gamma^N} \cup \overline{\Gamma^D} = \partial\Omega$ and $\Gamma^N \cap \Gamma^D = \emptyset$, and let us consider the problem of the form:
\begin{equation}
\label{eq:strong_state_pb}
\mathcal{P}(\Omega): \begin{cases} 
\mathcal{L} (u) &= f \qquad \text{in} \ \Omega, \\
u   &= u_D \quad \text{ on} \ \Gamma^D, \\
\partial_n u &= g \qquad \text{on} \ \Gamma^N, \\
\end{cases}
\end{equation}
where $\mathcal{L}$ is an elliptic differential operator and $\partial_n$ denotes the directional derivative normally outward to $\Omega$. For the sake of simplicity, we assume that $u_D \equiv 0$ and that $f$ and $g$ are sufficiently smooth and defined everywhere in $\mathbb{R}^2$, that is, $f \in L^2(\R^2)$ and $g\in H^1(\R^2)$. \\

This paper's aim is to introduce a strategy that reliably produces domain simplifications $\Omega_S \subset \mathbb{R}^2$ such that if $u_S$ is the function that solves the defeatured analogue of~\eqref{eq:strong_state_pb} over $\Omega_S$:
\begin{equation}
\label{eq:strong_state_pb_simplified}
\mathcal{P}(\Omega_S): \begin{cases} 
\mathcal{L} (u_S) &= f \qquad \text{in} \ \Omega_S, \\
u_S   &= 0 \qquad \text{on} \ \Gamma^D_S, \\
\partial_n u_S &= g \qquad \text{on} \ \Gamma^N_S, \\
\end{cases}
\end{equation}
then the discrepancy $\vert J(u, \Omega) - J(u_S, \Omega_S) \vert$ is small and furthermore controllable. In this paper, we assume that $\Omega$ is complex and that its (internal) parametrisation would hence require manual intervention and / or be prohibitively expensive. Therefore, a mandatory requirement is to completely avoid the parametrisation of $\Omega$ while solely operating on its boundary representation which is a standard input from a CAD system. The key steps of the methodology are summarised as follows:
\begin{enumerate}
    \item Creating an initial oversimplified domain representation $\Omega_0$ from the CAD input;
    \item Employing a gradient-based strategy that iteratively includes additional geometric features in regions of the computational model where a large impact on the quantity of interest can be expected.
\end{enumerate}
The result is a sequence of increasingly complex domains $\Omega_0, \Omega_1, \ldots, \Omega_n = \Omega_S$ with associated meshes $\mathcal{Q}_i$ and state problem approximations $u_i$ such that $J(u_i, \Omega_i), i \in \{0, \ldots, n\}$ is a sequence that approximates the quantity of interest with growing accuracy. \\

In the following, we choose $\mathcal{L}(u) = - \Delta u$ for simplicity. However, the methodology generalises to any generic elliptic operator in a straightforward way. Let us denote by $\mathcal U$ the \emph{universe set} of the problem at hand, that is, the collection of all admissible bounded domains $B \subset \mathbb{R}^2$ we may consider in the simplification process. Clearly $\Omega_0, \Omega \in \mathcal U$. We respectively denote by $\Gamma_B^D$ and $\Gamma_B^N$ the Dirichlet and the Neumann boundaries of a generic $B\in \mathcal U$ in which problem $\mathcal{P}(B)$ is solved. If we define $$H^1_{0, \Gamma^{D}_B}(B) := \{v \in H^1(B) \, \vert \, \text{the trace of } v \text{ vanishes on } \Gamma_B^D \},$$ the weak form of $\mathcal{P}(B)$ reads: \\
\begin{equation}
\label{eq:weak_state_pb}
    \text{Find } u \in H_{0, \Gamma_B^D}^1(B) \quad \text{such that} \quad a(u,v \,; B) = b(v \,; B) \quad \forall \ v \in H_{0, \Gamma_B^D}^1(B),
\end{equation}
with
\begin{equation}
    a(u,v \, ; B) \vcentcolon = \int_{B} \nabla u \cdot \nabla v \operatorname{d}\!\mathbf{x} \quad \text{and} \quad b(v\, ; B) = \int_{B} fv \operatorname{d}\!\mathbf{x} + \int_{\Gamma_B^N} g v \operatorname{d}\!x.
\end{equation}
The weak forms of problems~\eqref{eq:strong_state_pb} and~\eqref{eq:strong_state_pb_simplified} then follow from taking $B = \Omega$ and $B = \Omega_S$, respectively. \\

\noindent As previously mentioned, we are interested in computing a functional $J(\, \cdot \, ; \, \Omega) \vcentcolon H_{0, \Gamma^D}^1(\Omega) \rightarrow \R$ of the solution to problem $\mathcal{P}(\Omega)$. In practice, we are only able to compute the corresponding quantity on the simplified geometry $\Omega_S$, which, therefore, must be wisely chosen to minimize the error in the quantity of interest. For any element $w \in H^1(B)$, $B\in\mathcal U$, we assume that the functional $J$ takes the general form 
\begin{equation}
\label{eq:J_integral}
    J(w\,; B) = \int_{B^\prime} j\big(w\big)\operatorname{d}\!\mathbf{x} + \int_{\partial B^\prime}q\big(w\big) \operatorname{d}\!x,
\end{equation}
where $B^\prime \subseteq B$ and $\partial B^\prime \subseteq \partial B$ while $ j \vcentcolon H^1(B) \rightarrow L^1(B^\prime)$ and $q \vcentcolon H^1(B) \rightarrow L^1(\partial B^\prime)$ have locally Lipschitz continuous derivatives. Finally, we define the \emph{shape functional} $\mathcal{J} \vcentcolon \mathcal U \rightarrow \R$ as
\begin{equation}
\label{eq:def_shapefuncJ} 
\mathcal{J}(B) \vcentcolon = J\big(u(B)\,; B\big), 
\end{equation}
where $u(B)$ is the unique weak solution of $\mathcal{P}(B)$ as defined in~\eqref{eq:weak_state_pb}.\\

%\noindent In the same way, given a simplified domain $\Omega_S \subset U$, we consider the shape functional $\mathcal{J}(\Omega_S) = J\big(u(\Omega_S), \Omega_S \big)$, where $u(\Omega_S)$ is now the exact solution of $\mathcal{P}(\Omega_S)$, that is, of the model problem~\eqref{eq:strong_state_pb} defined on the simplified geometry. \\
\noindent Hence, what we call the modelling error due to the simplification of the geometry is the quantity: 
\begin{equation}
\label{eq:modelling_err}
\abs{\mathcal{J}(\Omega) - \mathcal{J}(\Omega_S)}.
\end{equation}

\begin{assumption}
In most cases, if $B = \Omega_i$, then $\Omega_i^\prime$ and $\partial \Omega_i^\prime$ in~\eqref{eq:J_integral} will constitute relatively small subsets of $\Omega$ and $\partial \Omega$, respectively. Therefore, we assume that $\Omega_i^\prime \subset \Omega$ and $\partial \Omega_i^\prime \subset \partial \Omega$ for any simplified geometry $\Omega_i$ we may construct. In this way, we can always compute the error \eqref{eq:modelling_err} and the simplification problem is well defined.
\end{assumption}

%%%%%%%%%%%%%%%%%%%%%%%%%%%%%%%%%%%%%%%%%%%%%%%%%%%%%%%%%%%%%%%%%%%%%

\section{Mathematical tools}\label{sec:math_tools}

% ----------------------------------------------------------------------------------------------------------------------------------------

\subsection{Shape gradients}\label{subsec:shapeder}
Before detailing our simplification algorithm, let us introduce some basic definitions and tools of first order shape calculus. For a complete discussion on the topic we refer the reader to \cite{ introshapeopti_zol, shape&geom_zol_delf}.\\

\noindent Let $B\subset \R^2$, $B \in \mathcal U$ be an open bounded domain with piecewise smooth boundary $\partial B$. To evaluate the sensitivity of a shape functional $\mathcal{F}(B)$ to a small change in the boundary $\partial B$, we define a class of admissible deformation maps following the \emph{perturbation of identity} approach  as in \cite{article_paganini}. That is,
let $\bm{\theta} \in C^1(\R^2; \R^2)$ be a vector field, and define an admissible perturbation map $\mathbf{T}_t \vcentcolon B \rightarrow \R^2$ as follows:
$$\mathbf{T}_t(\mathbf{x}) \vcentcolon = \mathbf{x} + t \bm{\theta}(\mathbf{x}), \qquad \text{with} \  \abs{t} \leq 1, \ \mathbf{x} \in \partial B$$
and suitably extended in the interior of $B$.

\begin{definition} 
\label{def:direct_shapeder} The shape functional $\mathcal{F} \vcentcolon \mathcal U \rightarrow \R$ admits a \emph{directional shape derivative} at $B$ along the displacement field $\bm{\theta}$ if the following limit exists and is finite:
$$\frac{\operatorname{d}\! \mathcal{F}}{\operatorname{d} \!\Omega}\Big(B\, \Big|\, \bm{\theta}\Big) \vcentcolon = \lim_{t \searrow 0}\frac{\mathcal{F}\big(\mathbf{T}_t(B)\big) - \mathcal{F}(B)}{t} .$$
\end{definition}
\noindent As it is quite common in the literature \cite{article_paganini, shape&geom_zol_delf}, we will also refer to the quantity above as the \emph{sensitivity} of $\mathcal{F}(B)$ with respect to the perturbation direction $\bm{\theta}$.
\begin{definition}
\label{def:shape_differentiable} 
The shape functional $\mathcal{F} \vcentcolon \mathcal U \rightarrow \R$ is \emph{shape differentiable} at $B$ if the mapping
$$ \frac{\operatorname{d}\! \mathcal{F}}{\operatorname{d} \!\Omega}\Big(B\, \Big|\, \cdot\Big) \vcentcolon C^1(\R^2; \R^2) \rightarrow \R, \qquad \bm{\theta} \mapsto \frac{\operatorname{d}\! \mathcal{F}}{\operatorname{d} \!\Omega}\Big(B\,\Big|\, \bm{\theta}\Big)$$ is a linear and continuous functional on $C^1(\R^2; \R^2)$. 
If $\mathcal{F}$ is shape differentiable, then \smash{$\displaystyle \tfrac{\operatorname{d}\! \mathcal{F}}{\operatorname{d} \!\Omega}(B\,|\, \cdot)$} is an element of the dual space of $C^1(\R^2; \R^2)$ and we refer to it as the \emph{shape gradient} of $\mathcal{F}$ at $B$.
\end{definition}
\noindent Explicit formulas for the shape gradient of $\mathcal{F}$ and its directional shape derivatives are available for several types of shape functionals \cite{shape&geom_zol_delf}. In particular, as we focus on PDE-constrained shape functionals $\mathcal{J}$ of the form \eqref{eq:def_shapefuncJ}, such formulas can be derived by introducing the \emph{Lagrangian} functional $L$ \cite{Verani_adaptiveFEM,Verani_discreteGradFlows} as follows.\\ %\vcentcolon H^1(\R^2) \times H^1(\R^2) \times U \rightarrow \R$ for any $w,v \in H^1(\R^2)$ and any $B \in U$ as
%\begin{equation}
%\label{eq:def_lagrangian}
% L(w,v \, ; B) \vcentcolon = J(w \, ; B) + a(w, v \, ; B) - b(v \, ; B).
% \end{equation}
% \noindent By considering a suitable extension to $H^1(\R^2)$ of the weak solution $u(D) \in H_{0, \Gamma_D}^1(D)$ to the \emph{state problem}~\eqref{eq:weak_state_pb} defined on $D$, and by computing the Lagrangian $L$ at $B=D$ and $w = u(D)$, we obtain:
% \begin{equation}
% \label{eq:Lagr_evaluated}
% L\big(u(D), v \, ; D\big) = J\big(u(D) \,; D \big) = \mathcal{J}(D), \qquad \forall \ v \in H^1(\R^2),
% \end{equation}
% that is, the shape functional whose derivative we are interested in computing.\\
% \noindent We also introduce the following auxiliary notation for the different types of derivatives we will consider:
% $ \displaystyle \frac{\partial L}{\partial \Omega}\Big (u(D), v\, ; D\,\Big |\, \bm{\theta}\Big)$ for the shape derivative at $D$ in the direction $\bm{\theta}$ of the shape functional obtained by evaluating $L(w, v ; B)$ at $w = u(D)$ and at any $v \in H^1(\R^2)$; moreover, let
% $\displaystyle \frac{\partial J}{\partial w}\Big(w \, ; D\, \Big|\, \psi\Big) \vcentcolon = \frac{\partial \big [\mathbb{J}(D)\big]}{\partial w}(w \,|\, \psi)$ and$ \displaystyle \frac{\partial L}{\partial w}\Big(w, v\, ; D \,\Big| \,\psi\Big) \vcentcolon$ be the G\^ateaux derivatives at $w$ in the direction $\psi \in H^1(D)$ of the functionals $\mathbb{J}(D)$ and $\mathbb{L}(D)$.\\

\noindent We denote by $[H^1(B) \times H^1(B)]^*$ the set of bilinear bounded functionals on $H^1(B) \times H^1(B)$ for $B \in \mathcal U$. Let $\mathbb{L} \vcentcolon \mathcal U \rightarrow \big\{[H^1(B) \times H^1(B)]^* \ | \ B \in \mathcal U\big\}$ be such that
\begin{equation}
    \begin{aligned}
    \mathbb{L}(B) \vcentcolon H^1(B) \times H^1(B) & \rightarrow \R \\
    (w, v) & \mapsto [\,\mathbb{L}(B)\,](w , v) \vcentcolon = L(w, v \,; B),
    \end{aligned}
\end{equation}
where the Lagrangian functional $L \vcentcolon H^1(B) \times H^1(B) \rightarrow \R $ is defined for any $w,v \in H^1(B)$ as
\begin{equation}
\label{eq:def_lagrangian}
L(w,v \, ; B) \vcentcolon = J(w \, ; B) + a(w, v \, ; B) - b(v \, ; B).
\end{equation}
\noindent By considering the weak solution $u(B) \in H_{0, \Gamma_B^D}^1(B)$ of the \emph{state problem}~\eqref{eq:weak_state_pb} defined on $B$, and by computing the Lagrangian functional $L$ at $w = u(B)$, we obtain:
\begin{equation}
\label{eq:Lagr_evaluated}
L\big(u(B), v \, ; B\big) = J\big(u(B) \,; B \big) = \mathcal{J}(B), \qquad \forall \ v \in H_{0,\Gamma_B^D}^1(B),
\end{equation}
that is, the shape functional whose derivative is the one we are interested in computing.\\

\noindent We furthermore introduce the following auxiliary notation for the different types of derivatives we will consider:
$\displaystyle{\frac{\partial L}{\partial \Omega}\left (u(B), v\, ; B\, |\, \bm{\theta} \right)}$ for the shape derivative at $B$ in the direction $\bm{\theta}$ of the shape functional obtained by evaluating $\mathbb{L}(B)$ at $u(B)$ and at any $v \in H_{0,\Gamma_B^D}^1(B)$; moreover, let
$$\frac{\partial L}{\partial w}\Big(w, v\, ; B \,\Big| \,\psi\Big) \vcentcolon = \frac{\partial \big [ \mathbb{L}(B)\big]}{\partial w}(w, v \,| \, \psi)% \qquad \text{and} \qquad \frac{\partial J}{\partial w}\Big(w \, ; B\, \Big|\, \psi\Big) % \vcentcolon = \frac{\partial \big [\mathbb{J}(B)\big]}{\partial w}(w \,|\, \psi)
$$
be the G\^ateaux derivative at $w$ in the direction $\psi \in H^1(B)$ of the functional $\mathbb{L}(B)$, and let $\displaystyle\frac{\partial J}{\partial w}\Big(w \, ; B\, \Big|\, \psi\Big)$ and $\displaystyle\frac{\partial a}{\partial w}\Big(w, v\, ; B \, \Big |\, \psi \Big)$ be defined in a similar fashion.\\

\noindent We can use property~\eqref{eq:Lagr_evaluated} of the Lagrangian and the chain rule to write for every $v \in H_{0,\Gamma_B^D}^1(B)$,
 \begin{equation}
 \label{eq:Lagr_evaluated_b}
 \frac{\operatorname{d}\! \mathcal{J}}{\operatorname{d} \!\Omega}\Big(B\,\Big|\, \bm{\theta}\Big) = \frac{\partial L}{\partial \Omega} \Big(u(B), v\, ; B \,\Big|\, \bm{\theta}\Big) + \frac{\partial L}{\partial w}\Big(u(B), v\, ; B \, \Big| \, \frac{\operatorname{d}\! u}{\operatorname{d}\!\Omega}(B \, |\, \bm{\theta})\Big).
 \end{equation}
The term $\displaystyle \frac{\operatorname{d}\! u}{\operatorname{d}\!\Omega}(B \,|\, \bm{\theta})\in H^1_{0,\Gamma_B^D}(B)$ denotes the shape derivative of the function $u(B)$. The shape derivative of a function can be defined similarly to the shape derivative of a functional from Definition~\ref{def:direct_shapeder}, under suitable regularity assumptions. The interested reader is referred to \cite{Henrot_bookshapeder_clear} for more details.
%The last term in the right hand side above writes, for any $\psi \in H_{0,\Gamma_D}^1(D)$,
%$$ \frac{\partial L}{\partial v}\Big(u(D), v\, ; D \, \Big| \, \psi \Big) = a(u(D),\, \psi \, ; D) - b(\psi \, ; D).$$
%Hence, since for our model problem $\displaystyle \frac{\operatorname{d}\! v}{\operatorname{d}\!\Omega}(D \,|\, \bm{\theta}) \in H_{0,\Gamma_D}^1(D)$ (see \cite{Verani_discreteGradFlows} and references therein), this G\^ateaux derivative vanishes since $u(D)$ solves the state problem~\eqref{eq:weak_state_pb}.\\
Then, the Lagrange multiplier $v$ is chosen so that the second term in the right hand side of~\eqref{eq:Lagr_evaluated_b} vanishes. That is, if we recall definition~\eqref{eq:def_lagrangian} of the Lagrangian $L$, it is chosen such that
$$ 0= \frac{\partial L}{\partial w}\Big(u(B), v\, ; B \, \Big| \, \psi \Big) = \frac{\partial J}{\partial w}\Big( u(B) \, ; B \,\Big|\, \psi \Big) + \frac{\partial a}{\partial w}\Big(u(B), v\, ; B \, \Big |\, \psi \Big), \quad \forall \ \psi \in H_{0,\Gamma_B^D}^1(B).$$
This equality holds in particular for the choice $\displaystyle \psi = \frac{\operatorname{d}\! u}{\operatorname{d}\!\Omega}(B \,|\, \bm{\theta})\in H^1_{0,\Gamma_B^D}(B)$, see also \cite{Verani_discreteGradFlows}. \\

\noindent Let us finally introduce the following \emph{adjoint problem}: find $p \equiv p(B) \in  H_{0,\Gamma_B^D}^1(B)$ such that
\begin{equation}
    \label{eq:weak_adj_pb}
    \frac{\partial a}{\partial w}\Big(u(B), p \, ; B \, \Big |\, \psi \Big) = - \frac{\partial J}{\partial w}\Big( u(B) \, ; B \,\Big|\, \psi \Big), \quad \forall \ \psi \in H_{0,\Gamma_B^D}^1(B),
\end{equation}
where $\displaystyle \frac{\partial a}{\partial w}\Big(u(B),\, p(B) \, ; B \, \big |\, \psi \Big) = a(\psi,\, p(B) \, ; B)$ for the bilinear form associated with the model problem considered in this paper, namely Poisson's equation. \\

\noindent Consequently, the formula for the shape derivative of $\mathcal{J}$ at $B$ in the direction $\bm{\theta}$ involves the solutions $u(B)$ and $p(B)$ of the state~\eqref{eq:weak_state_pb} and adjoint~\eqref{eq:weak_adj_pb} problems, and it reads
\begin{equation}
\label{eq:shape_gradient}
\frac{\operatorname{d}\! \mathcal{J}}{\operatorname{d} \!\Omega}\Big(B\,\Big|\, \bm{\theta}\Big) = \frac{\partial L}{\partial \Omega} \Big(u(B),\, p(B)\, ; B \,\Big|\, \bm{\theta}\Big),
\end{equation}
where the term on the right hand side can be computed by standard transformation techniques~\cite{article_paganini}.

\subsection{HB and THB splines}\label{subsec:hbsplines}
\noindent This section introduces the notation and reviews the main concepts related to hierarchical B-splines (or HB-splines) and their truncated counterparts (called THB-splines) in the bidimensional setting, closely following \cite{geopde_hier}. These functions are used to build hierarchical geometries in the same way as classical B-spline curves, surfaces and volumes, with the advantage of allowing for local refinement, overcoming the limitations intrinsic to the tensor-product structure of multidimensional B-splines. \\

\noindent Let us consider a nested sequence of Cartesian grids $\{\mathcal{Q}_{\ell}^c\}_{\ell \in \mathbb{N}_0}$ of the parametric domain $\widehat{\Omega} = (0,1)^2$. For a given $\ell \in \mathbb{N}_0$, we call \textit{cell of level $\ell$} any $Q \in \mathcal{Q}_{\ell}^c$.
Then, let $\{\mathcal{B}_\ell\}_{\ell \in \mathbb{N}_0}$ be a sequence of (bivariate) tensor product B-spline bases determined by their degree and knot vectors and associated to $\{\mathcal{Q}_{\ell}^c\}_{\ell \in \mathbb{N}_0}$, such that the spaces they generate are nested as follows:
$$ \spn{\mathcal{B}_0} \subset \spn{\mathcal{B}_1} \subset \spn{\mathcal{B}_2} \subset \dots .$$

\noindent For each $\ell \in \mathbb{N}_0$, we denote by $N_\ell$ the dimension of the space of level $\ell$, so that $\mathcal{B}_\ell \vcentcolon = \{ b_{i,\ell} \ | \ i = 1, \dots , N_\ell\}$.
A characteristic property of B-splines is the so-called \emph{two-scale relation}, which enables expressing the $\beta_{i, \ell} \in \mathcal{B}_{\ell}$ as linear combinations of the basis functions $\beta_{j, \ell +1} \in \mathcal{B}_{\ell +1}$ with non-negative coefficients $c_{j, \ell + 1}( \beta_{i, \ell}) \geq 0$, that is:
\begin{equation}
\label{eq:two_scale}
    \beta_{i,\ell} = \sum_{j = 1}^{N_{\ell + 1}} c_{j, \ell + 1}( \beta_{i, \ell}) \, \beta_{j, \ell + 1}, \qquad \forall \ \beta_{i, \ell} \in \mathcal{B}_\ell.
\end{equation}
Thanks to the nestedness of the $\left \{ \mathcal{B}_\ell \right \}_{\ell \in \mathbb{N}_0}$, this relation can be recursively applied to relate basis functions of non consecutive levels.  \\

\noindent Then, consider a hierarchy of subdomains of depth $M$, defined as the set $ \bm{\omega}_M \vcentcolon = \{ \omega_0, \omega_1 , \dots, \omega_M\} $ where
$$ \widehat{\Omega} = \omega_0 \supset \omega_1 \supset \dots \supset \omega_M = \emptyset.$$ 
Each subdomain $\omega_\ell$ is the union of closed cells of level $\ell -1$ while $\omega_0$ coincides with the parametric domain $\widehat{\Omega}$.
\noindent The hierarchical basis $\mathcal{H} \equiv \mathcal{H}(\bm{\omega}_M ) $ associated to the set $\bm{\omega}_M$ is then defined recursively with the algorithm:
\begin{equation*}
\label{eq:def_hbsplines}
    \begin{cases}
    \mathcal{H}_0 &\vcentcolon = \mathcal{B}_0, \\
    \mathcal{H}_{\ell + 1} &\vcentcolon = \{\beta \in \mathcal{H}_\ell \ | \ \supp{\beta} \not\subset \omega_{\ell+1}\} \cup \{\beta \in \mathcal{B}_{\ell+1} \ | \ \supp{\beta} \subset \omega_{\ell+1}\}, \quad \ell = 0, \dots, M-2, \\
    \mathcal{H} &\vcentcolon = \mathcal{H}_{M-1}.  
    \end{cases}
\end{equation*}
We refer the reader to \cite{buffa2017refinable} for details about this construction.
The corresponding \emph{hierarchical mesh} $\mathcal{Q} \equiv \mathcal{Q}(\bm{\omega}_M)$ is given by the set $$ \mathcal{Q} \vcentcolon = \bigcup_{\ell = 0}^{M-1} \{ Q \in \mathcal{Q}_\ell^c \ | \ Q \subset \omega_\ell \ \text{and} \ Q \not\subset \omega_{\ell+1}\}.
$$
\noindent Henceforth, we call \emph{active basis functions of level $\ell$} the functions $\beta \in \mathcal{B}_\ell$ such that $\beta \in \mathcal{H} \cap \mathcal{B}_\ell$, i.e., whose support covers only cells of level $\bar{\ell} \geq \ell$ and at least one cell of level $\ell$; we also name \emph{deactivated basis functions of level $\ell$} the functions $\beta \in \mathcal{B}_\ell$ such that $\beta \in \mathcal{H}_{\ell} \setminus \mathcal{H}_{\ell + 1}$, i.e., whose support is entirely contained in $\omega_{\ell + 1}$.\\
 
\noindent It is well established (see \cite{thb_giannelli_juttler}) that the space spanned by $\mathcal{H}$ can be alternatively described using a THB-spline basis associated to the same hierarchical mesh $\mathcal{Q}$, carrying better properties from a numerical standpoint, such as the partition of unity and a reduced support of the basis functions. To define this new basis, we must introduce a suitable truncation operator trunc$_{\ell+1}$, which is then iteratively applied to eliminate from coarser B-splines the contribution of refined basis functions on the overlapping parts of their supports. Thus, let us call trunc$_{\ell+1}$ the following operator, defined on a tensor product space of B-splines:
$$ \text{trunc}_{\ell+1}(\beta_{i,\ell}) \vcentcolon = \sum_{j = 1}^{N_{\ell + 1}} c^{\tau}_{j, \ell + 1}( \beta_{i, \ell}) \, \beta_{j, \ell + 1}, \qquad \forall \ \beta_{i, \ell} \in \mathcal{B}_\ell, $$
with coefficients:
\begin{align}
     \label{eq:two_scale_trunc}
     c^{\tau}_{j, \ell + 1}( \beta_{i, \ell}) = 
    \begin{cases}
    0 \qquad & \text{if} \ \beta_{j, \ell + 1} \in \mathcal{H}_{\ell + 1} \cap \mathcal{B}_{\ell + 1}, \\
    c_{j, \ell + 1}( \beta_{i, \ell}) \qquad & \text{otherwise}.
    \end{cases}
\end{align}
 
\noindent Now, the truncated hierarchical basis $\mathcal{T} \equiv \mathcal{T}(\bm{\omega}_M ) $ associated to the set $\bm{\omega}_M$ is defined by replacing the simple B-spline functions with their truncation in the same recursive algorithm employed for $\mathcal{H}$, that is:
 \begin{equation*}
\label{eq:def_thbsplines}
    \begin{cases}
    \mathcal{T}_0 &\vcentcolon = \mathcal{B}_0, \\
    \mathcal{T}_{\ell + 1} &\vcentcolon = \{\text{trunc}_{\ell+1}(\beta) \ | \ \beta \in \mathcal{T}_{\ell} \ \text{and} \ \supp{\beta} \not\subset \omega_{\ell+1}\} \cup \{\beta \in \mathcal{B}_{\ell+1} \ | \ \supp{\beta} \subset \omega_{\ell+1}\}, \quad \ell = 0, \dots, M-2, \\
    \mathcal{T} &\vcentcolon = \mathcal{T}_{M-1}.  
    \end{cases}
\end{equation*}
Note that the truncation from~\eqref{eq:two_scale_trunc} eliminates the terms associated to the active and deactivated functions of level $\ell + 1$ in the two-scale relation~\eqref{eq:two_scale}. Therefore, trunc$_{\ell+1}(\beta)$, with $\beta \in \mathcal{T}_\ell$, no longer belongs to $\mathcal{B}_\ell$ and must be expressed in terms of the basis of the finest level that truncates it.\\

\noindent In this work, we deal with geometries $B\in\mathcal U$ parameterized by THB splines, that is, given a THB basis $\mathcal{T} =\{\beta_i\}_{i=1}^{N} $, and a sequence $\{ \mathbf{c}_i \}_{i=1}^{N} \subset \R^2$ of control points, $B$ is explicitly given by the following mapping $\mathbf{x}$ of the parametric domain $\widehat{\Omega}$:
$$ \mathbf{x} \vcentcolon \widehat{\Omega} \rightarrow B \qquad \text{such that} \qquad \mathbf{x}(\bm{\xi}) = \sum_{i=1}^{N_\mathcal{T}} \mathbf{c}_i\beta_i(\bm{\xi}) \qquad \forall  \ \bm{\xi} \in \widehat{\Omega}. $$
% ----------------------------------------------------------------------------------------------------------------------------------------
\subsection{Parametrization method and mesh generation}
\label{sect:EGG}
The purpose of this paper is to devise an adaptive strategy for finding an approximation $\Omega_S$ of $\Omega$ such that $\vert \mathcal{J}(\Omega) - \mathcal{J}(\Omega_S) \vert$ is below some user-defined threshold $\varepsilon > 0$. To this end, the methodology employs the concept of shape gradients, as discussed in Section \ref{subsec:shapeder}. On the one hand, this paper assumes that no explicit parameterisation $\mathbf{x}: \widehat{\Omega} \rightarrow \Omega$ of the interior of $\Omega$ is given. On the other hand, both the state (\ref{eq:weak_state_pb}) and the adjoint problem (\ref{eq:weak_adj_pb}) rely on a valid parameterisation $\mathbf{x}_S: \widehat{\Omega} \rightarrow \Omega_S$ of the interior of $\Omega_S$ at each defeaturing iteration. The algorithm therefore requires solving an additional surface-to-volume (StV) problem $\partial \Omega_S \rightarrow \Omega_S$. More precisely, given a parametric description $\mathbf{F}_S: \partial \widehat{\Omega} \rightarrow \partial \Omega_S$ of the defeatured domain's boundary, we need to construct a valid parameterisaton of $\Omega_S$, i.e., a bijective map $\mathbf{x}_S: \hat{\Omega} \rightarrow \Omega_S$. \\

The StV-problem is challenging and notoriously difficult to automate for general geometries. As we seek for an autonomously operating framework, the parametrisation strategy adopted by this paper should:
\begin{itemize}
    \item Reliably compute parameterisations for general geometries $\Omega_S$ that are topologically equivalent to $\widehat{\Omega} = (0, 1)^2$ without manual intervention.
    \item Be compatible with and exploit local refinement made possible by the use of hierarchical splines.
\end{itemize}
With these requirements in mind, the basic parameterisation strategy adopted by this paper approximates the inverse of a map $\mathbf{h}: \Omega_S \rightarrow \widehat{\Omega}$ whose components $\mathbf{h}_i$, $i \in \{1, 2\}$, are harmonic functions in $\Omega_S$, i.e.,
\begin{align}
\label{eq:harmonic_map}
    i \in \{1, 2\}: \enskip \Delta \mathbf{h}_i = 0 \text{ in } \Omega_S, \quad \text{ s.t.} \quad \mathbf{h} \vert_{\partial \Omega_S} = \mathbf{F}_S^{-1}.
\end{align}
In the following, we drop the $S$ subscript for convenience. Equation (\ref{eq:harmonic_map}) can be inverted for the map $\mathbf{x} = \mathbf{h}^{-1}$ that we are interested in. This leads to the following system of coupled quasi-linear second-order elliptic PDEs in nondivergence form (see \cite[Chapter~7]{knupp2020fundamentals} and \cite{hinz2020goal}):
\begin{align}
\label{eq:EGG}
    i \in \{1, 2 \}: \enskip A \left( \mathbf{x} \right) \colon H(\mathbf{x}_i) = 0 \text{ in } \widehat{\Omega}, \quad \text{s.t. } \mathbf{x} \vert_{\partial \widehat{\Omega}} = \mathbf{F},
\end{align}
with
\begin{align}
\label{eq:definition_hessian_diffusivity_EGG}
    H(\phi)_{ij} := \frac{\partial^2 \phi}{\partial \xi_i \partial \xi_j} \quad \text{and} \quad A \left( \mathbf{x} \right) = \begin{pmatrix} g_{22} & -g_{12} \\ -g_{12} & g_{11} \end{pmatrix}.
\end{align}
Here, $\boldsymbol{\xi} = (\xi_1, \xi_2)^T$ contains the free coordinate functions in $\widehat{\Omega}$ while the $g_{ij} = \displaystyle\frac{\partial \mathbf{x}}{\partial \xi_i} \cdot \frac{\partial \mathbf{x}}{\partial \xi_j}$ denote the entries of the mapping's metric tensor and $A \, \colon  B$ the Frobenius inner product. A justification for seeking the geometry parameterisation as the inverse of a map that is harmonic in $\Omega_S$ is provided by the Rad\'o-Kneser-Choquet theorem \cite{gravesen2012planar, choquet1945type} which states that any $\mathbf{h}: \Omega \rightarrow \widehat{\Omega}$ whose entries satisfy a maximum principle is a diffeomorphism in the interior of $\Omega \subset \mathbb{R}^2$, provided its target domain $\widehat{\Omega} \subset \mathbb{R}^2$ is convex. \\

\noindent Seeking a mapping $\mathbf{x}_h: \widehat{\Omega} \rightarrow \Omega_S$ that approximately solves (\ref{eq:EGG}) over a given hierarchical spline space that is compatible with the boundary correspondence $\mathbf{F}: \widehat{\Omega} \rightarrow \Omega$ hence constitutes a natural choice. The finite element treatment of nondivergence form equations, despite still being in its infancy, has become a rich research topic with the first publication due to Lakkis and Pryer \cite{lakkis2011finite} in $2011$. To start with, consider the scalar, linear and homogeneous nondivergence form problem that reads
\begin{align}
\label{eq:nondivergence_strong}
    B \, \colon \, H(v) + \text{ lower order terms}= r, \quad \text{such that} \quad v \vert_{\partial \widehat{\Omega}} = 0,
\end{align}
with $r \in L^2(\widehat{\Omega})$ and $B: \widehat{\Omega} \rightarrow \mathbb{R}^{2 \times 2}$ uniformly elliptic, i.e., $B \in L^{\infty}(\widehat{\Omega}; \mathbb{R}^{2 \times 2})$ while there are constants $0 < c_1 \leq c_2 < \infty$ such that
\begin{align}
    c_1 \leq \inf \limits_{\boldsymbol{\xi} \in \mathbb{R}^2, \left \| \boldsymbol{\xi} \right \| = 1} \boldsymbol{\xi}^T B \boldsymbol{\xi} \leq \sup \limits_{\boldsymbol{\xi} \in \mathbb{R}^2, \left \| \boldsymbol{\xi} \right \| = 1} \boldsymbol{\xi}^T B \boldsymbol{\xi} \leq c_2 \quad \text{almost everywhere in } \widehat{\Omega}.
\end{align}
In the following, let us disregard the lower order terms in (\ref{eq:nondivergence_strong}).
Finite element methods that follow from a variational formulation of (\ref{eq:nondivergence_strong}) seek an approximation $v_h$ of $v$ based on the following general Petrov-Galerkin approach \cite{gallistl2017variational}:
\begin{align}
\label{eq:nondivergence_petrov_galerkin}
    \text{find } v_h \in \mathcal{V}_h \cap H^1_0(\widehat{\Omega}) \quad \text{s.t.} \quad \int \limits_{\widehat{\Omega}} \tau(w_h) B \, \colon H(v_h) \mathrm{d} \boldsymbol{\xi} = \int \limits_{\widehat{\Omega}} \tau(w_h) r \mathrm{d} \boldsymbol{\xi} \quad \forall w_h \in \mathcal{V}_h \cap H^1_0(\widehat{\Omega}),
\end{align}
for some $\tau: H^2(\widehat{\Omega}) \rightarrow L^2(\widehat{\Omega})$ and a finite-dimensional trial and test space $\mathcal{V}_h \subset H^2(\widehat{\Omega})$. Here, the most common choice is $\tau(v) = \gamma \Delta v$, with $\gamma = \operatorname{tr}(B) / \left \| B \right \|_F^2$, for which coercivity of the associated bilinear form is well-established for convex domains $\widehat{\Omega}$ \cite{blechschmidt2019error}. \\

% \noindent Conforming discretisations are difficult to achieve as Lagrangian FEM bases are generally only in $H^1(\widehat{\Omega})$. This can not be remedied by performing partial integration on the weak form as this is prohibited by the essential boundedness of $A$. Hence, approaches based on (\ref{eq:nondivergence_petrov_galerkin}) furthermore aim to reduce the regularity requirements by (a) introducing independent variables for the first- or second-order derivatives of $u_h$ (mixed formulation) \cite{lakkis2011finite, gallistl2017variational} or restoring the coercivity of the bilinear form by the use of interior penalty terms acting on the gradient of $u_h$ ($C^0$ discontinuous Galerkin formulation) \cite{blechschmidt2019error, dedner2013discontinuous, feng2017finite, smears2013discontinuous}. \\

\noindent The THB-spline setting enables a conforming discretisation of~\eqref{eq:nondivergence_petrov_galerkin} by choosing $\mathcal{V}_h \subset H^2(\widehat{\Omega})$. Assuming that the boundary correspondence $\mathbf{F}: \partial \widehat{\Omega} \rightarrow \partial \Omega_S$ is compatible with a given THB-spline basis $\mathcal{T} \subset H^2(\widehat{\Omega})$ whose linear span generates $\mathcal{V}_h$, a computational approach for approximating the solution of (\ref{eq:EGG}) follows from applying the Petrov-Galerkin scheme (\ref{eq:nondivergence_petrov_galerkin}) to a Newton linearisation and iterating until the associated residual is below a user-defined threshold. To be precise, we define $\mathcal{F}: \mathcal{V}_h^2 \times \mathcal{V}_{h, 0}^2 \rightarrow \mathbb{R}$, with $\mathcal{V}_h^2 \vcentcolon= \mathcal{V}_h \times \mathcal{V}_h$ and $\mathcal{V}_{h, 0} \vcentcolon= \mathcal{V}_h \cap H^1_0(\hat{\Omega})$. The semi-linear form $\mathcal{G}(\, \cdot \, , \, \cdot \,)$ reads:
\begin{align}
    \mathcal{G}(\mathbf{x}, \boldsymbol{\sigma}) := \sum_{i=1}^2 \int \limits_{\widehat{\Omega}} \boldsymbol{\tau}_i(\mathbf{x}, \boldsymbol{\sigma}) A(\mathbf{x}) \, \colon H(\mathbf{x}_i) \operatorname{d} \! \boldsymbol{\xi},
\end{align}
with $A(\mathbf{x})$ as in (\ref{eq:definition_hessian_diffusivity_EGG}) and $\boldsymbol{\tau}_i \left(\mathbf{x}, \boldsymbol{\sigma} \right) = \gamma(\mathbf{x}) \Delta \boldsymbol{\sigma}_i$, where $\gamma(\mathbf{x}) = \operatorname{tr}(A(\mathbf{x})) / \left \| A(\mathbf{x}) \right \|_F^2$. We denote by $\mathcal{G}^\prime \left(\mathbf{x}, \boldsymbol{\sigma} \vert \boldsymbol{\psi} \right)$ the G\^ateaux derivative of $\mathcal{G}(\, \cdot \,, \, \cdot \,)$ with respect to its first argument in the direction of $\boldsymbol{\psi}$. Given the $k$-th iterate $\mathbf{x}^k \in \mathcal{V}_h^2$ that satisfies $\mathbf{x}^k \vert_{\partial \widehat{\Omega}} = \mathbf{F}$, we compute the Newton increment by solving
\begin{align}
\label{eq:EGG_discretized_newton}
    \text{find } \partial \mathbf{x}^k \in \mathcal{V}_{h, 0}^2 \quad \text{s.t.} \quad \mathcal{G}^{\prime} \left( \mathbf{x}^k, \boldsymbol{\sigma}_h \, \vert \, \partial \mathbf{x}^k \right) = - \mathcal{G}(\mathbf{x}^k, \boldsymbol{\sigma}_h), \quad \forall \boldsymbol{\sigma}_h \in \mathcal{V}_{h, 0}^2.
\end{align}
Upon completion, the next iterate is given by $\mathbf{x}^{k+1} = \mathbf{x}^k + \kappa \partial \mathbf{x}^k$, for some $\kappa \in (0, 1]$ whose optimal value is estimated using a line search routine. The iteration is terminated once
\begin{align*}
    \max \limits_{\boldsymbol{\sigma}_h \in \tilde{\mathcal{T}}^2} \vert \mathcal{G}(\mathbf{x}^k, \boldsymbol{\sigma}_h) \vert, \quad \text{with} \quad \tilde{\mathcal{T}} := \{ \beta \in \mathcal{T} \, \vert \, \beta \in H^1_0(\widehat{\Omega}) \},
\end{align*}
is below a user-defined threshold $\mu > 0$. Heuristically, a good stopping criterion is choosing $\mu = 10^{-6}$, regardless of the characteristic length scale of the geometry. Hereby, scaling invariance results from multiplying by $\gamma(\cdot)$ in $\boldsymbol{\tau}(\cdot, \cdot)$. With this choice of $\mu$, the scheme typically converges after $5$ iterations whenever the initial iterate $\mathbf{x}^0$ approximates a pair of harmonic functions in $\widehat{\Omega}$ (instead of $\Omega$). Hereby, the approximation follows from discretising the Laplace equation, with pure Dirichlet boundary conditions resulting from $\mathbf{F}: \widehat{\Omega} \rightarrow \Omega$, in the usual way. \\

Due to the numerical scheme's truncation error, the theoretically predicted validity property may fail to carry over to its numerical approximation. In case of grid folding, we locally refine $\mathcal{T}$ by replacing all elements
\begin{align*}
    Q \in \bigcup \limits_{\mathclap{\beta \in \mathcal{T}_{-}}} \operatorname{supp}(\beta),
\end{align*}
by their finer counterparts from the next level in the element hierarchy, where $\mathcal{T}_{-} \subset \mathcal{T}$ denotes the set of all $\beta \in \mathcal{T}$ that are nonvanishing on an element $Q \in \mathcal{Q}$ on which a negative value of
\begin{align*}
    \det D \mathbf{x}_h := \det \frac{\partial \mathbf{x}_h}{\partial \boldsymbol{\xi}}
\end{align*}
has been detected. The value of $\det D \mathbf{x}_h$ is computed on all points that correspond to the abscissae of the Gaussian quadrature scheme that is used to approximate the state problem $u$, the associated adjoint $p$ and the shape gradient. Upon element refinement, the defective solution $\mathbf{x}_h \in \mathcal{V}_h^2$ is prolonged to the enriched THB spline basis $\mathcal{T}^+$ and the Newton scheme is restarted using the prolongation as an initial iterate. These steps are repeated until $\det D \mathbf{x}_h$ is strictly positive on all quadrature points. \\

The above procedure reliably computes valid maps for a wide range of input geometries while being algorithmically lightweight and hence easy to automate.

%%%%%%%%%%%%%%%%%%%%%%%%%%%%%%%%%%%%%%%%%%%%%%%%%%%%%%%%%%%%%%%%%%%%%

\section{Adaptive gradient-based simplification algorithm}\label{sec:algo}
%\textcolor{red}{\begin{itemize}
%\item add simple picture to explain
%\end{itemize}}

\noindent This section proposes an adaptive algorithm aimed at constructing a simplified geometry $\Omega_S$ that yields a sufficiently accurate value of the functional of interest $\mathcal{J}(\cdot)$. Ideally, we would like to control the discrepancy in the quantity of interest such that
\begin{align}
\label{eq:convergence_criterion}
    \abs{\mathcal{J}(\Omega) - \mathcal{J}(\Omega_S)} < \varepsilon
\end{align}
for a given tolerance $\varepsilon$. However, achieving quantitative error control is beyond the scope of this manuscript. Instead, our termination criterion is based on the norm of the shape gradient with respect to a set of shape directions. The relation between the discrepancy in $\mathcal{J}(\cdot)$ and the shape gradient will be investigated experimentally in Section \ref{sec:numerical}. \\

% \noindent Let us start by commenting on its name, to highlight its main features. 
\noindent As in other \emph{adaptive} frameworks, we devised an iterative strategy to build increasingly detailed geometries based on the following steps:
\begin{center}
    FIT $\rightarrow$ SOLVE $\rightarrow$ ESTIMATE $\rightarrow$ MARK $\rightarrow$ REFINE.
\end{center}
\noindent In our context, the FIT step creates a parametric description of $\partial \Omega_n \approx \partial \Omega$ which is compatible with the THB-spline space at the current iteration $n$.
The SOLVE step corresponds to the computation of the functional $\mathcal{J}(\cdot)$ on the current geometry approximation $\Omega_n$, which furthermore requires solving the StV problem $\partial \Omega_n \rightarrow \Omega_n$ beforehand. Regarding the remaining steps, rather than directly selecting a set of local features to create the refined geometry $\Omega_{n+1}$, we construct a locally refined THB-spline space to represent it. Hence, in the ESTIMATE step, we decide whether the geometry has been represented sufficiently accurately. If not, we choose a suitable set of basis functions in the MARK step, while the REFINE step enriches the given THB-spline space starting from this set to increase the complexity of the geometry. \\

In the following, we discuss the key steps of the algorithm in detail.
% \noindent Our algorithm is \emph{gradient based} for two reasons: on the one side, the marking and refinement procedures rely on the evaluation of the (first order) sensitivity of the functional at hand to a modification of the domain (see Subsections~\ref{subsec:algo_mark} and \ref{subsec:algo_refine}); on the other, our heuristic estimator for the modelling error is also defined through the same values of the first order sensitivities (see Subsection~\ref{subsec:algo_estimate}). \\
\subsection{Fit}\label{subsec:algo_fit}
\noindent The $n$-th iteration of the algorithm is associated with a hierarchical mesh $\mathcal{Q}_n$ along with its canonical THB-spline basis $\mathcal{T}_n$ with cardinality $\dim(\mathcal{T}_n) := N_n$. Given the current $(\mathcal{Q}_n , \mathcal{T}_n)$, as a basic ingredient, the algorithm requires approximating the exact boundary correspondence $\mathbf{F}: \partial \widehat{\Omega} \rightarrow \partial \Omega$ over the subset $\partial \mathcal{T}_n \subset \mathcal{T}_n$ with $\dim(\partial \mathcal{T}_n) := \partial N_n$, comprised of all $\beta_i \in \mathcal{T}_n$ that do not vanish on $\partial \widehat{\Omega}$. The approximation is of the form
$$ \mathbf{F}_n \vcentcolon \partial \widehat{\Omega} \rightarrow \partial \Omega_n \qquad \text{such that} \qquad \mathbf{F}_n(\bm{\xi}) = \sum_{i=1}^{\partial N_n} \mathbf{c}_{i}\beta_{i}(\bm{\xi}) \qquad \forall  \ \bm{\xi} \in \partial \widehat{\Omega}.$$
The routine hence requires properly selecting the boundary control points \smash{$\partial \mathcal{C} = \{\mathbf{c}_i \}_{i=1}^{\partial N_n}$}. We accomplish this by fitting the \smash{$\mathbf{c}_k \in \partial \mathcal{C}$} via a scaled $H^1(\partial \widehat{\Omega})$ projection. To be precise, the fitting routine selects the $\mathbf{c}_k \in \partial \mathcal{C}$ such that
\begin{align}
\label{eq:boundary_fit}
    & \kappa_0 \int \limits_{\partial \widehat{\Omega}} \left \| \mathbf{F}_n - \mathbf{F} \right \|^2 \operatorname{d}\!x + \kappa_1 \int \limits_{\partial \widehat{\Omega}} \left \| \frac{\partial(\mathbf{F}_n - \mathbf{F})}{\partial t} \right \|^2 \operatorname{d}\!x \longrightarrow \min \limits_{\partial \mathcal{C}},
\end{align}
wherein {$\displaystyle\frac{\partial}{\partial t}(\, \cdot \,)$} denotes the directional derivative along the unit tangent of $\partial \widehat{\Omega}$, while $\kappa_0 > 0$ and $\kappa_1 \geq 0$ are constants tuning the magnitude of the $L^2(\partial \widehat{\Omega})$ and $H^1(\partial \widehat{\Omega})$ semi-norms in the fit. The minimisation from (\ref{eq:boundary_fit}) may furthermore be subjected to constraints such as requiring the fit to assume the exact value of $\mathbf{F} \vcentcolon \partial \widehat{\Omega} \rightarrow \partial \Omega$ in the four corners of the parametric domain.

%%%%%%%%%%%%%%%%%%%%%%%%%%%%%%%%%%%%%%%%%%%%%%%%%%%%%%%%%%%%%%%%%%%%%%%%%%
\subsection{Solve}\label{subsec:algo_solve}
Given the boundary fit $\mathbf{F}_n: \partial \widehat{\Omega} \rightarrow \partial \Omega_n$ generated in Subsection \ref{subsec:algo_fit}, this step computes the value of the functional $\mathcal{J}(\Omega_n)$. As a first step, the correspondence $\mathbf{F}_n: \partial \widehat{\Omega} \rightarrow \partial \Omega_n$ is forwarded to the parameterisation pipeline from Section \ref{sect:EGG} to create a valid $\mathbf{x}_n: \widehat{\Omega} \rightarrow \Omega_n$ that parameterises $\Omega_n$. The map $\mathbf{x}_n$ enables us to numerically solve the state problem~\eqref{eq:weak_state_pb} on $\Omega_n$. Moreover, we also solve the adjoint problem~\eqref{eq:weak_adj_pb} on the same geometry, as its solution is needed for the computation of the shape derivatives as explained in Section~\ref{subsec:shapeder}.\\
\noindent Since the main focus of this work concerns the error due to geometric simplification, we assume that an exact solution of the PDE problem is available on the simplified geometry. In practice, we ensure that the THB meshes utilised for computing the state and adjoint problems are sufficiently fine that the impact of the truncation error is negligible.

%%%%%%%%%%%%%%%%%%%%%%%%%%%%%%%%%%%%%%%%%%%%%%%%%%%%%%%%%%%%%%%%%%%%%%%%%%
% \subsection{Estimate}\label{subsec:algo_estimate}
%To achieve our goal of controlling the \textcolor{red}{modelling error} on simplified geometries, we must introduce a suitably computable estimator of the error~\eqref{eq:modelling_err} that we use as a stopping criterion of our iterative algorithm. As mentioned before, our choice is based on the first order shape derivatives of $\mathcal{J}$: indeed, given a finite set $\Theta_n$ of directions $\bm{\theta}_k$,
%a stopping criterion operates on the shape gradient, i.e., on the
%$$j_k \in \left\{ \tfrac{\operatorname{d}\! \mathcal{J}}{\operatorname{d} \!\Omega}(\Omega_n\, |\, \bm{\theta}_k) \enskip \vert \enskip \bm{\theta}_k \in \Theta_n \right\}.$$
%Clearly, several reasonable stopping criteria exist. The choice of $\Theta_n$ as well as reasonable stopping criteria are discussed in the next Subsection~\ref{subsec:algo_mark}.

%The numerical experiments show that this heuristic estimator provides a good upper bound for the modelling error, although a proof of its reliability and efficiency goes beyond the scope of this work.\\

%%%%%%%%%%%%%%%%%%%%%%%%%%%%%%%%%%%%%%%%%%%%%%%%%%%%%%%%%%%%%%%%%%%%%%%%%%
\subsection{Estimate}\label{subsec:algo_est}
\noindent We are given the current iteration's $(\mathcal{Q}_n, \mathcal{T}_n)$ along with the boundary fit $\mathbf{F}_n: \partial \widehat{\Omega} \rightarrow \partial \Omega_n$ expressed in $\operatorname{span} \partial \mathcal{T}_n$ and the mapping $\mathbf{x}_n: \widehat{\Omega} \rightarrow \Omega_n$ that satisfies $\mathbf{x}_n \vert_{\partial \widehat{\Omega}} = \mathbf{F}_n$. Let then $\mathcal{Q}_n^+$ be the hierarchical mesh obtained by dyadic refinement of the boundary elements of $\mathcal{Q}_n$, that is, of those elements $Q \in \mathcal{Q}_n$ with $\partial Q \cap \partial \widehat{\Omega} \neq \emptyset$. By \smash{$\mathcal{T}_n^+ = \{\beta^+_{i}\}_{i=1}^{N_n^+}$} we denote the canonical THB-spline basis associated with $\mathcal{Q}^+_n$. As before, we consider the subset $\partial \mathcal{T}_n^+ \subset \mathcal{T}_n^+$ of the $\partial N_n^+$ basis functions that do not vanish on $\partial \widehat{\Omega}$.
Thanks to the \emph{two-scale relation}~\eqref{eq:two_scale}, we can express $\mathbf{F}_n: \partial \widehat{\Omega} \rightarrow \partial \Omega_n$ in terms of this subset:
$$  \mathbf{F}_n(\bm{\xi}) = \sum_{i=1}^{\partial N^+_n} \mathbf{\overline{c}}_{i} \beta^+_{i}(\bm{\xi}) \qquad \forall  \ \bm{\xi} \in \partial \widehat{\Omega}.$$
Our goal is to create a set of THB splines \smash{$\partial \mathcal{T}_{n+1}$} to express the boundary at the next iteration. To this aim, we choose some of the basis functions in \smash{$\partial \mathcal{T}_n^+$} according to their impact on the functional of interest, measured in terms of the value of the shape derivative of $\mathcal{J}(\cdot)$ in the corresponding direction. In the following, we only need to clarify how we associate a direction of differentiation $\bm{\theta}_k$ to a basis function.\\

\noindent We accomplish this by minimizing~\eqref{eq:boundary_fit} over the set $\partial \mathcal{C}^+$ of boundary control points associated with $\partial \mathcal{T}_n^+$, which yields a more accurate fit $\partial \Omega_n^+$, with
$$ \mathbf{F}^+_n \vcentcolon \partial \widehat{\Omega} \rightarrow \partial \Omega^+_n \qquad \text{such that} \qquad \mathbf{F}^+_n(\bm{\xi}) = \sum_{i=1}^{\partial N^+_{n}} \mathbf{c}^+_{i}\beta^+_{i}(\bm{\xi}) \qquad \forall  \ \bm{\xi} \in \partial \widehat{\Omega}.$$ \\
\noindent Expressing the two fits in the same basis divides their discrepancy into the shape directions
$$\beta_i^+ \underbrace{(\mathbf{c}_i^+ - \overline{\mathbf{c}}_i)}_{\Delta \mathbf{c}_i^+}, \enskip i \in \{1, \ldots, \partial N_n^+ \}.$$
As such, each direction naturally associates with exactly one $\beta_i^+ \in \partial \mathcal{T}^+_n$. We define $\Theta_n$ as the set that follows from splitting aforementioned shape directions into their $x$ and $y$ components, i.e.,
\begin{align}
\label{eq:Theta_n}
    \Theta_n = \left \{ \bm{\theta}_1^x, \bm{\theta}_1^y, \bm{\theta}_2^x, \bm{\theta}_2^y, \ldots, \bm{\theta}_{\partial N^+_n}^x, \bm{\theta}_{\partial N^+_n}^y \right \}, \quad \text{where} \quad \left \{ \begin{array}{ll} \bm{\theta}_i^x(\boldsymbol{x}) = (\Delta \mathbf{c}_{i}^+ \cdot {\mathbf{e}_x}) \beta_i^+(\mathbf{x}^{-1}_n(\boldsymbol{x})) {\mathbf{e}}_x \\ \bm{\theta}_i^y(\boldsymbol{x}) = (\Delta \mathbf{c}_i^+ \cdot {\mathbf{e}}_y) \beta_i^+(\mathbf{x}^{-1}_n(\boldsymbol{x})) {\mathbf{e}}_y \end{array} \right.
\end{align}
for all $\boldsymbol{x}\in \Omega_n$, with $\boldsymbol{e}_x := (1,0)T$ and $\boldsymbol{e}_y=(0,1)^T$. 
Each $\beta_k^+$ is now associated with the directions $\{ \bm{\theta}^{x}_k, \bm{\theta}^{y}_k \} \subset \Theta_n$ along with the basis function wise two-component shape gradient
\begin{align}
\label{eq:shape_gradient_beta_k}
\nabla \mathcal{J}(\Omega_n, \beta_i^+) := \left( \frac{\operatorname{d}\! \mathcal{J}}{\operatorname{d} \!\Omega}(\Omega_n |\, \bm{\theta}_i^x), \frac{\operatorname{d}\! \mathcal{J}}{\operatorname{d} \!\Omega}(\Omega_n |\, \bm{\theta}_i^y) \right)^T.
\end{align}
A reasonable stopping criterion now terminates the algorithm once all gradients are below a threshold, for instance in the $l_2$-norm. Letting
\begin{align}
\label{eq:E_n_max_l2}
E_n^{\text{mod}}(\Omega_n) := \max \limits_{\beta_k^+ \in \partial T^+_n} \left \| \nabla \mathcal{J}(\Omega_n, \beta_k^+) \right \|_2, \quad \text{the algorithm is terminated once} \quad E_n^{\text{mod}}(\Omega_n) < \varepsilon,
\end{align}
which concludes the ESTIMATE step. \\

\subsection{Mark}\label{subsec:algo_mark2}
This section deals with picking the directions of differentiation $\bm{\theta}_k$ and employing the associated shape derivatives in the marking process which is then translated to a proper refinement of the $Q \in \mathcal{Q}_n$ in the REFINE step (see Section \ref{subsec:algo_refine}).\\

\noindent Following the same reasoning as for the ESTIMATE step, we compute the shape derivatives associated to all the boundary degrees of freedom and then select only the most significant ones according to a maximum strategy: given a parameter $\alpha \in (0,1)$, the set \smash{$\mathcal{D}^+_{\text{bdry}}  \subset \partial \mathcal{T}_n^+$} contains all the basis functions $\beta^+_{k}$ such that the associated shape gradient satisfies:
\begin{equation}
\label{eq:selection_criterion}
   \left \| \nabla \mathcal{J}(\Omega_n, \beta_i^+) \right \|_2 \geq \alpha \, E_n^{\text{mod}}(\Omega_n).
\end{equation}
In the context of shape directions, it is useful to define the set
\begin{align}
    \widehat{\Theta}_n = \left \{ \widehat{\bm{\theta}}_1^x, \widehat{\bm{\theta}}_1^y, \ldots, \widehat{\bm{\theta}}_{\partial N_n^+}^x, \widehat{\bm{\theta}}_{\partial N_n^+}^y \right \},
\end{align}
where for all $\boldsymbol{x}=(x,y)\in\Omega_n$,
$$\widehat{\bm{\theta}}_i^x(\boldsymbol{x}) = \beta_i^+(\mathbf{x}_n^{-1}(\boldsymbol{x}))\boldsymbol{e}_x, \quad \widehat{\bm{\theta}}_i^y(\boldsymbol{x}) = \beta_i^+(\mathbf{x}_n^{-1}(\boldsymbol{x}))\boldsymbol{e}_y,$$
and the function-wise unit shape gradient
\begin{align}
    \nabla \widehat{\mathcal{J}}(\Omega_n, \beta_i^+) = \left( \frac{\operatorname{d}\! \mathcal{J}}{\operatorname{d} \!\Omega}(\Omega_n |\, \hat{\bm{\theta}}_i^x), \frac{\operatorname{d}\! \mathcal{J}}{\operatorname{d} \!\Omega}(\Omega_n |\, \hat{\bm{\theta}}_i^y) \right)^T.
\end{align}
Thanks to the linearity of the shape gradient and recalling equations~(\ref{eq:Theta_n})--(\ref{eq:shape_gradient_beta_k}), %we may equivalently write
% \begin{align}
%    \nabla \mathcal{J}(\Omega_n, \beta_i^+) = \nabla \widehat{\mathcal{J}}(\Omega_n, \beta_i^+) \odot \Delta \mathbf{c}_i^+,
% \end{align}
% wherein $A \odot B$ denotes the Hadamard product \cite{horn1990hadamard} between vectors or matrices $A$ and $B$ (of equal dimension). Again, by linearization, 
we may write
\begin{align}
\label{eq:discrepancy_J_Omega_plus_Omega_approx}
    \mathcal{J}(\Omega_n^+) - \mathcal{J}(\Omega_n) \approx \sum_{i = 1}^{\partial N_n^+} \nabla \widehat{\mathcal{J}}(\Omega_n, \beta_i^+) \cdot \Delta \mathbf{c}_i^+.
\end{align}

%%%%%%%%%%%%%%%%%%%%%%%%%%%%%%%%%%%%%%%%%%%%%%%%%%%%%%%%%%%%%%%%%%%%%%%%%%
\subsection{Refine}\label{subsec:algo_refine}
\noindent The purpose of the refinement step is to translate the marked $\beta^+ \in \mathcal{D}^+_{\text{bdry}}$, cf.~\eqref{eq:selection_criterion}, to an appropriate dyadic refinement of boundary elements $Q \in \mathcal{Q}_n$. Upon performing the refinement $\mathcal{Q}_n \rightarrow \mathcal{Q}_{n+1}$, the associated canonical THB-spline basis should:
\begin{itemize}
    \item contain additional boundary degrees of freedom (DOFs) on segments of $\partial \widehat{\Omega}$ where the marked basis functions in $\mathcal{D}_{\text{bdry}}^+$ suggest that a large impact on the cost function evaluation can be expected;
    \item Avoid over- and under-refinement.
\end{itemize}

\begin{figure}[h!]
\centering
\begin{subfigure}[t]{0.45\textwidth}
    \centering
    \includegraphics[width=0.9\linewidth]{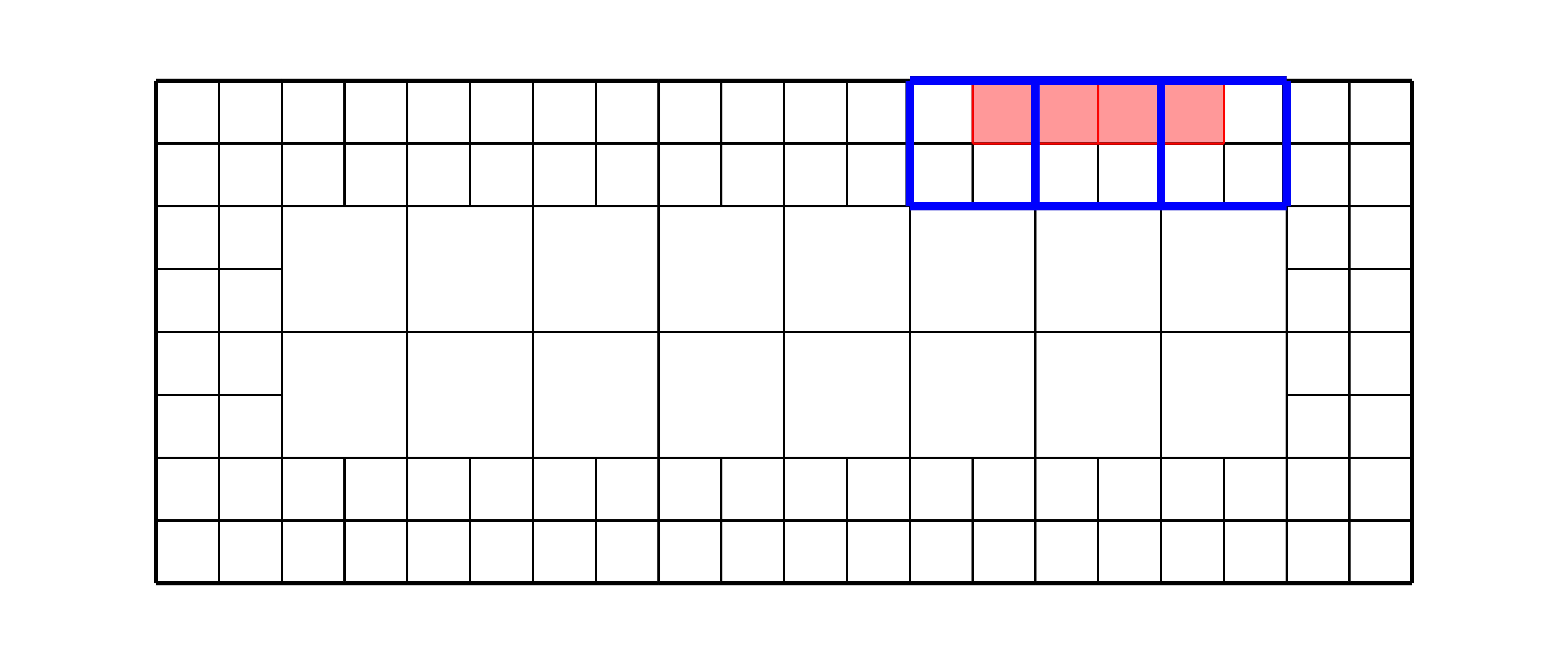}
    \caption{The $Q \in \mathcal{Q}^+$ highlighted in $\mathcal{Q}_n^+$ (red) and the overlapping elements from $\mathcal{Q}_n$ (blue).}
\end{subfigure} $\quad$
\begin{subfigure}[t]{0.45 \textwidth}
    \centering
    \includegraphics[width=0.9\linewidth]{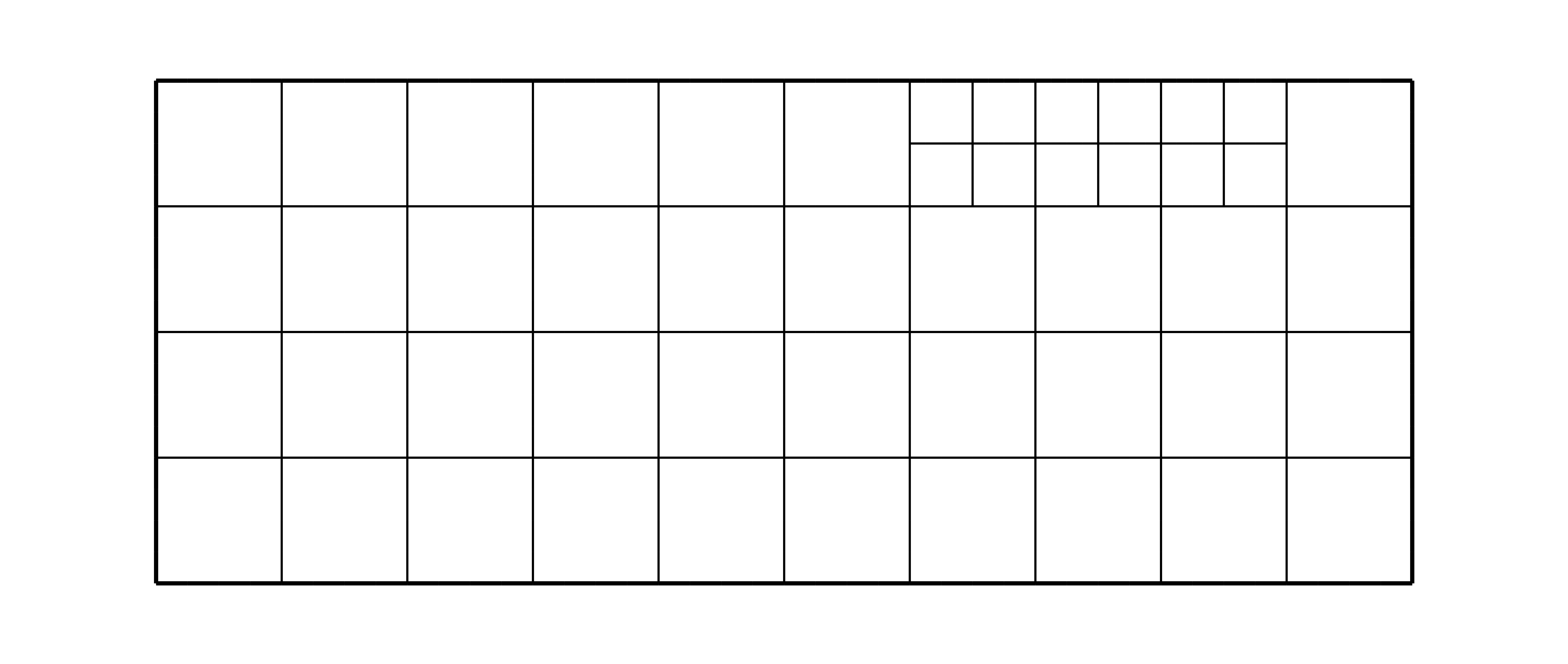}
    \caption{Illustration of the slight over-refinement in $\mathcal{Q}_{n+1}$.}
\end{subfigure}
\caption{Figure that illustrates the slight over-refinement resulting from refining all $Q \in \mathcal{Q}_n$ that overlap with $\mathcal{Q}^+$.}
\label{fig:slight_over_refinement}
\end{figure}

\noindent Here, over-/under-refinement is to be understood in the following ways: 1) providing too few additional local DOFs to accurately capture required local changes upon refitting the boundary contours; 2) introducing more additional DOFs than strictly necessary to reach the convergence threshold~\eqref{eq:convergence_criterion}; and 3) accidentally performing element refinement that has no effect on the enriched THB-spline basis. Clearly, a reasonable starting point is to define the set 
$$\mathcal{Q}^+ = \bigcup \limits_{\mathclap{\beta^+ \in \mathcal{D}^+_{\text{bdry}}}} \left\{ Q^+ \in \mathcal{Q}_n^+ \, \vert \, Q^+ \subseteq \operatorname{supp}\left(\beta^+ \right) \right \},$$
with $\mathcal{Q}_n^+$ as defined in Subsection \ref{subsec:algo_est},
and refining all $Q \in \mathcal{Q}_n$ that partially overlap with the $Q^+ \in \mathcal{Q}^+$. This selection ensures that upon refinement, for all $Q^+ \in \left(\mathcal{Q}_n^+ \cap \mathcal{Q}^+ \right)$, we have $Q^+ \in \mathcal{Q}_{n+1}$. Hereby, over- and under-refinement is avoided by suitably choosing the parameter $\alpha \in (0, 1)$ in~\eqref{eq:selection_criterion}. It should be noted that this method may produce elements $Q \in \left(\mathcal{Q}_{n+1} \cap \mathcal{Q}_n^+ \right)$ with $Q \notin \mathcal{Q}^+$ (see Figure \ref{fig:slight_over_refinement}), however, this effect is mild and thus deemed acceptable.

\begin{figure}[h!]
\centering
\begin{subfigure}[t]{0.45\textwidth}
    \centering
    \includegraphics[width=0.9\linewidth]{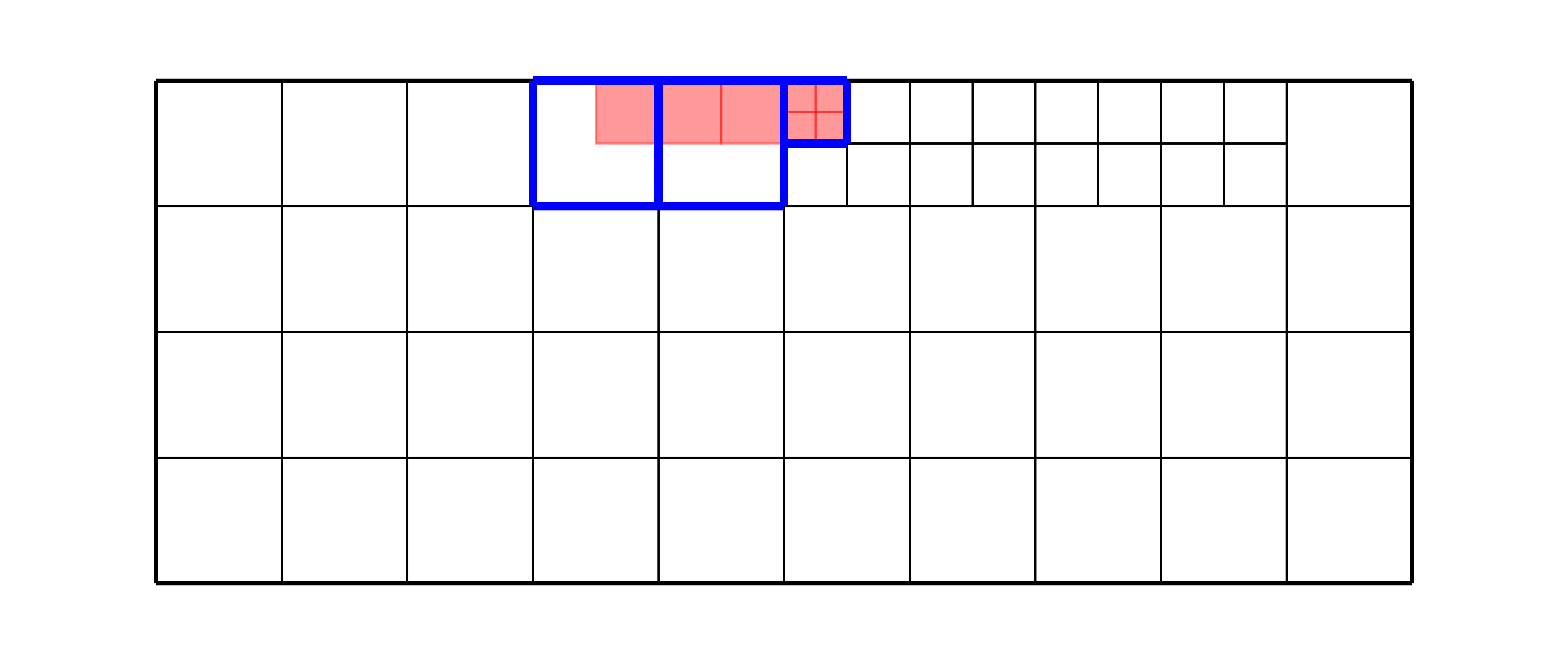}
    \caption{The $Q \in \mathcal{Q}^+$ highlighted on top of $\mathcal{Q}_n^+$ (red) and the overlapping elements from $\mathcal{Q}_n$ (blue).}
\end{subfigure} $\quad$
\begin{subfigure}[t]{0.45 \textwidth}
    \centering
    \includegraphics[width=0.9\linewidth]{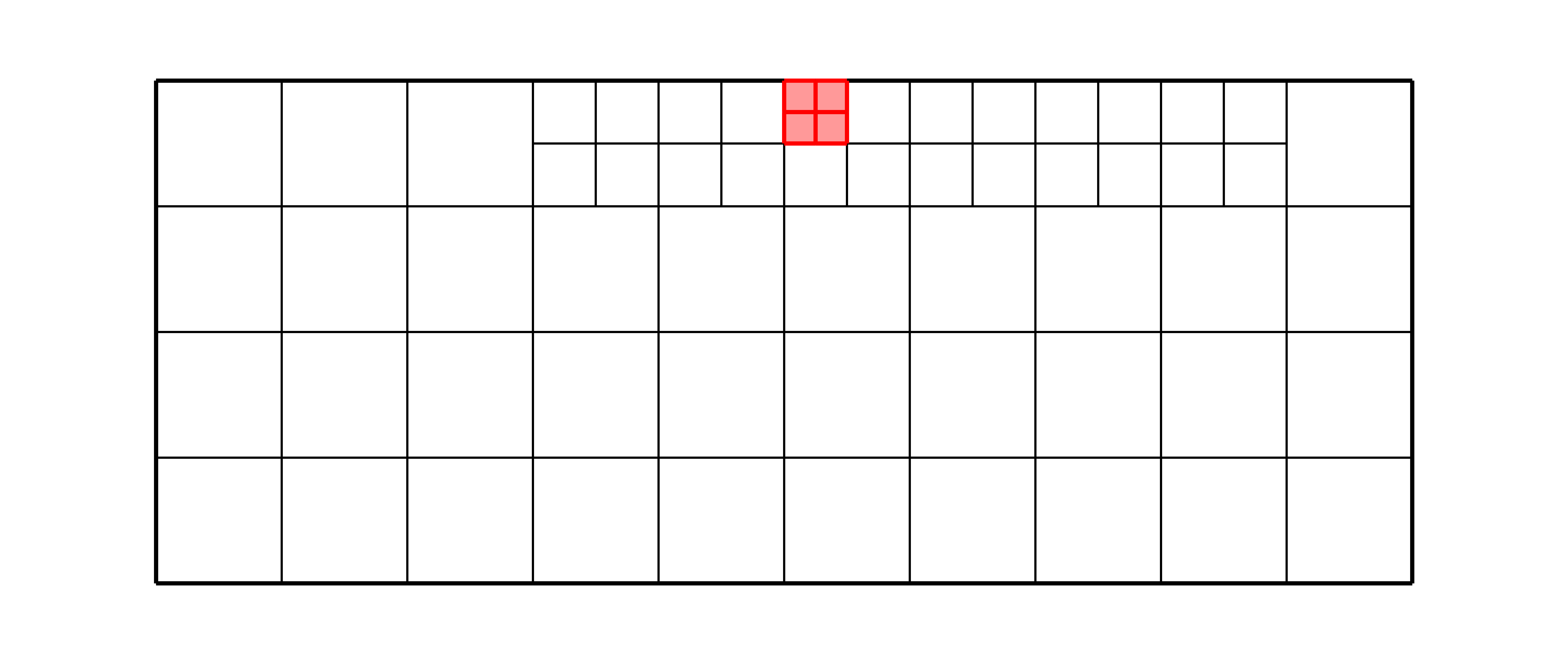}
    \caption{Illustration of over-refinement in $\mathcal{Q}_{n+1}$. The highlighted elements do not impact the basis cardinality for polynomial degrees $p \geq 2$.}
\end{subfigure}
\caption{Figure that illustrates how refining all $Q \in \mathcal{Q}_n$ that intersect with $\mathcal{Q}^+$ can produce redundant elements.}
\label{fig:over_refinement}
\end{figure}

\noindent A more critical unwanted effect is introducing elements $Q \in \mathcal{Q}^+$ that have no effect on the underlying basis cardinality and hence serve no purpose. For an example, see Figure \ref{fig:over_refinement}.

\begin{figure}[h!]
\centering
\begin{subfigure}[t]{0.45\textwidth}
    \centering
    \includegraphics[width=0.9\linewidth]{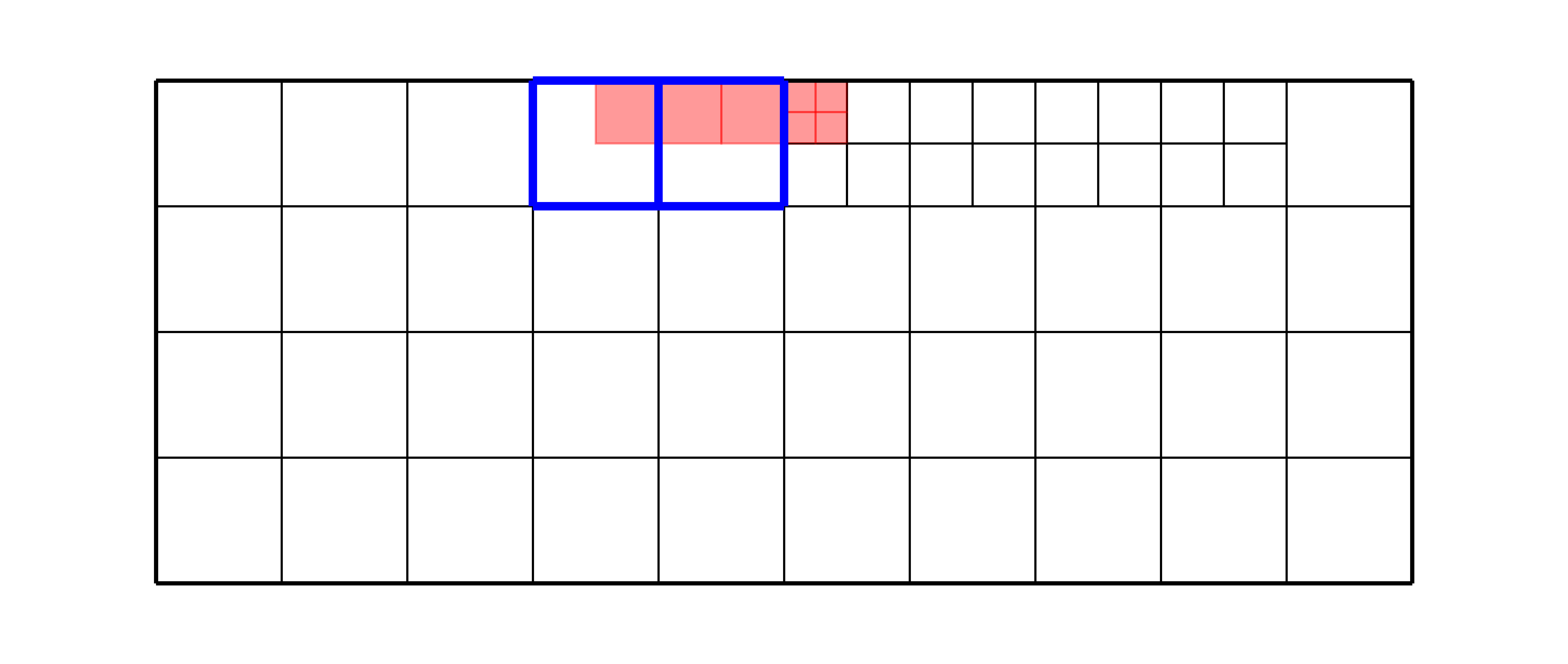}
    \caption{Supporting elements of one $\beta \in \mathcal{D}_{\text{bdry}}^+$ (red) along with the $Q \in \mathcal{Q}_{\text{act}}^{*}\left(\beta \right)$ (blue).}
\end{subfigure} $\quad$
\begin{subfigure}[t]{0.45 \textwidth}
    \centering
    \includegraphics[width=0.9\linewidth]{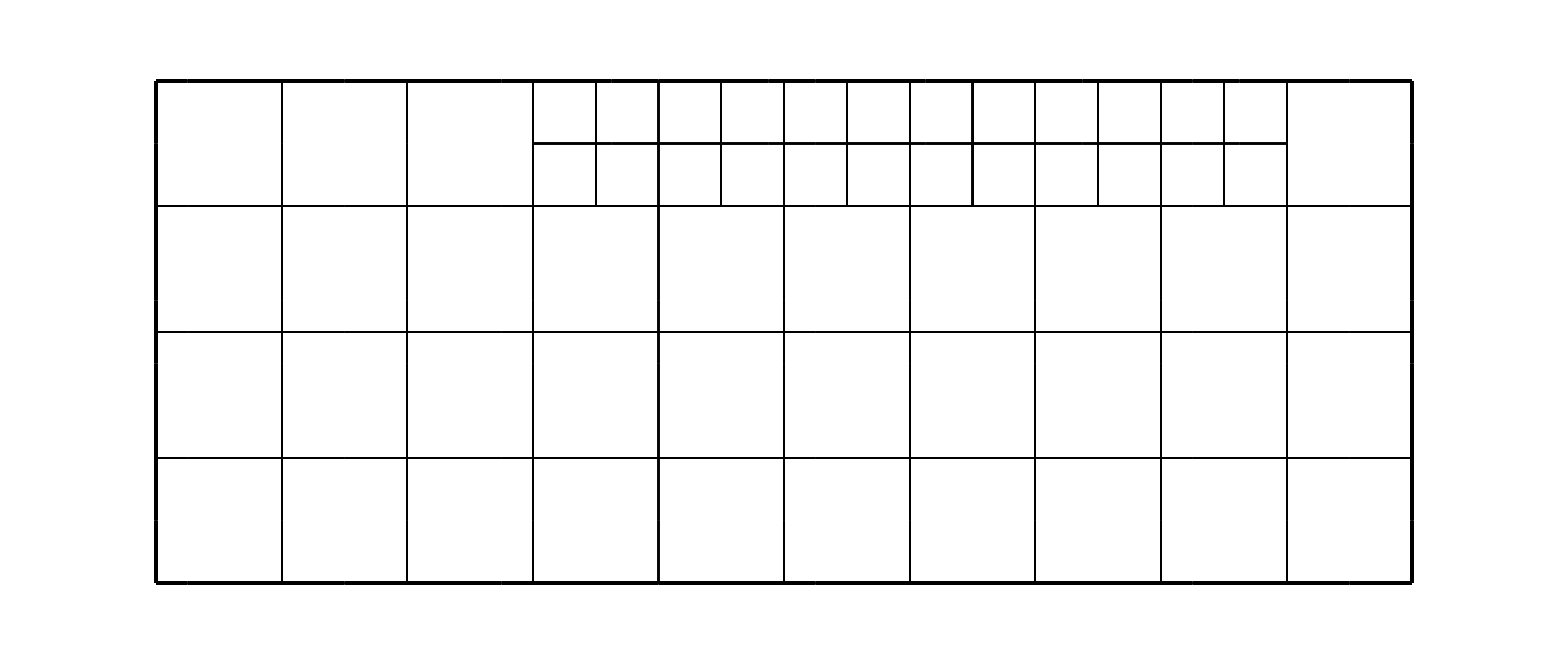}
    \caption{Minimal refinement does not introduce redundant elements.}
\end{subfigure}
\caption{Figure that illustrates how minimal refinement avoids redundant refinements.}
\label{fig:min_refinement}
\end{figure}

\noindent To avoid redundant refinements, in the following, we introduce the concept of \textit{minimal refinement}. Let
\begin{align}
    \mathcal{Q}_{\text{act}}\left(\beta^+ \right) = \left \{ Q \in \mathcal{Q}_n \enskip \vert \enskip Q \cap \operatorname{supp}\left(\beta^+ \right) \neq \emptyset \right \},
\end{align}
where the subscript \textit{`act'} stands for \textit{`active'}. Let $\ell_0 \geq 0$ be the index of the lowest level of the cells in $\mathcal{Q}_{\text{act}}(\beta^+)$ and let $\mathcal{Q}^*_{\text{act}}(\beta^+)$ be the subset of cells $Q \in \mathcal{Q}_{\text{act}}(\beta^+)$ composed of cells from the level $\ell_0$ in the mesh hierarchy. Minimal refinement then refines the set of elements
\begin{align}
\label{eq:min_refinement}
    \mathcal{Q}_{\text{min}}(\mathcal{D}_{\text{bdry}^+}) := \bigcup \limits_{\mathclap{\beta^+ \in \mathcal{D}_{\text{bdry}}^+}} \mathcal{Q}_{\text{act}}^{*}(\beta^+).
\end{align}
For an illustration of how minimal refinement avoids redundant refinements, see Figure \ref{fig:min_refinement}.

\noindent As refinement based on $\mathcal{Q}_{\text{min}}$ refines only a subset of the boundary elements of $\mathcal{Q}_n$, we can express both $\mathbf{F}_n: \partial \hat{\Omega} \rightarrow \partial \Omega_n$ and $\mathbf{F}_{n+1}: \partial \hat{\Omega} \rightarrow \partial \Omega_{n+1}$ as well as their discrepancy in $\partial \mathcal{T}_n^+$:
\begin{align}
\label{eq:discrepancy_J_Omega_next_Omega_approx}
    \mathbf{F}_{n+1}(\bm{\xi}) - \mathbf{F}_n(\bm{\xi}) = \sum_{i=1}^{\partial N_n^+} \Delta \mathbf{c}_i \beta_i^+(\bm{\xi}).
\end{align}
Then, we can estimate (c.f. equation \eqref{eq:discrepancy_J_Omega_plus_Omega_approx}):
\begin{align}
    \mathcal{J}(\Omega_{n+1}) - \mathcal{J}(\Omega_n) \approx \sum_{i = 1}^{\partial N_n^+} \nabla \widehat{\mathcal{J}}(\Omega_n, \beta_i^+) \cdot \Delta \mathbf{c}_i.
\end{align}
The purpose of the MARK and REFINE steps is then ensuring that
$$\sum_{i = 1}^{\partial N_n^+} \nabla \widehat{\mathcal{J}}(\Omega_n, \beta_i^+) \cdot \left(\Delta \mathbf{c}_i - \Delta \mathbf{c}_i^+ \right)$$
is small, while the cardinality of $\partial \mathcal{T}_{n+1}$ should be lower than that of $\partial \mathcal{T}_n^+$. This justifies our selection criterion from \eqref{eq:selection_criterion}, which allows for tuning the degrees to which both requirements are satisfied by properly selecting the value of $\alpha$ in~(\ref{eq:selection_criterion}).

%%%%%%%%%%%%%%%%%%%%%%%%%%%%%%%%%%%%%%%%%%%%%%%%%%%%%%%%%%%%%%%%%%%%%%%%%%%%%%%%%%%%%%%%%%%%%%%%%%%%%%%%%%%%%%%%%%%%%
\subsection{Implementation}
This section combines the
\begin{center}
   FIT $\rightarrow$ SOLVE $\rightarrow$ ESTIMATE $\rightarrow$ MARK $\rightarrow$ REFINE
\end{center}
strategies introduced in Sections \ref{subsec:algo_fit} -- \ref{subsec:algo_refine} into an autonomously-operating defeaturing algorithm. As input, the algorithm takes the exact boundary correspondence $\mathbf{F}: \partial \widehat{\Omega} \rightarrow \partial \Omega$, the quantity of interest $\mathcal{J}(\, \cdot \,)$, a coarse initial mesh $\mathcal{Q}_0$ of $\widehat{\Omega}$, the defeaturing tolerance $\varepsilon \in \mathbb{R}^+$ and the marking parameter $\alpha \in (0, 1)$. \\

\noindent In each iteration $n$, the first step approximates the exact correspondence $\mathbf{F}$ over the boundary DOFs of the current THB-spline basis $\mathcal{T}_n$ associated with the current mesh $\mathcal{Q}_n$ which is initialised to $\mathcal{Q}_0$ in the first iteration. The minimisation of~\eqref{eq:boundary_fit} yields the approximate correspondence $\mathbf{F}_n$ which is then forwarded to the parameterisation routine from section \ref{sect:EGG} to compute the mapping \smash{$\mathbf{x}_n: \widehat{\Omega} \rightarrow \Omega_n$}, while performing \textit{a posteriori} (aPos) refinement if necessary. Hereby, we initialise the parameterisation routine's mesh to the current $\mathcal{Q}_n$ plus the previous iteration's inner refinements (if applicable) that were performed to acquire a bijective map. This significantly reduces the required number of aPos refinements until a valid map is found. As a next step, we pick suitable mesh and basis pairs $(\mathcal{Q}_n^u, \mathcal{T}_n^u)$ and $(\mathcal{Q}_n^p, \mathcal{T}_n^p)$ for the approximation of $u_h^n$ and $p_h^n$ of the current geometry's state- and adjoint-equations. Instead of making this choice explicitly, we only list a number of requirements:
\begin{itemize}
    \item $(\mathcal{Q}_n^{u}, \mathcal{Q}_n^{p})$ follow from a sequence of hierarchical refinements of $\mathcal{Q}_0$.
    \item $(u_h^n, p_h^n)$ are suitable approximations to accurately compute the shape derivatives.
    \item The numerical error stemming from the aproximation $u(\Omega_n) \approx u_h^n$ does not dominate the geometrical error caused by evaluating $\mathcal{J}(\, \cdot \,)$ in $\Omega_n$ instead of $\Omega$.
\end{itemize}
Concrete basis selection strategies depend on the application. In this work, the state and adjoint THB meshes are always selected in a way that largely isolates the geometrical error.
Upon acquisition of the state and adjoint approximations $(u_h^n, p_h^n)$, we are in the position to approximately evaluate the cost function $\mathcal{J}(\Omega_n) \approx J_n := J(u_h^n, \Omega_n)$, cf.~\eqref{eq:J_integral}. \\

As a next step, the algorithm refines all boundary elements of the current $\mathcal{Q}_n$ to obtain the refined mesh / basis pair $(\mathcal{Q}_n^+, \mathcal{T}_n^+)$. The pair is forwarded to the minimisation of~\eqref{eq:boundary_fit} which yields the refined boundary fit $\mathbf{F}_n^+$. The discrepancy $\Delta \mathbf{F}_n = \mathbf{F}_n^+ - \mathbf{F}_n$ is expressed in the linear span of $\partial \mathcal{T}_n^+ \subset \mathcal{T}_n^+$ and the associated control points are utilised to define the directions of differentiation $\Theta_n$. With $\Theta_n, u_h^n$ and $p_h^n$ at hand, we assemble the approximate function-wise shape gradients
%$$\frac{\mathrm{d} \mathcal{J}}{\mathrm{d} \Omega_n} = \left \{ \frac{\mathrm{d} L}{\mathrm{d} \Omega}\big(u_h, p_h; \Omega_n \big \vert \boldsymbol{\theta}_k \big) \, \big \vert \, \boldsymbol{\theta}_k \in \boldsymbol{\Theta}_n \right\},$$
$$\nabla \mathcal{J}_h(\Omega_n, \beta_i^+) := \left( \frac{\mathrm{d} L}{\mathrm{d} \Omega}\big(u_h^n, p_h^n; \Omega_n \big \vert \boldsymbol{\theta}_i^x \big), \frac{\mathrm{d} L}{\mathrm{d} \Omega}\big(u_h^n, p_h^n; \Omega_n \big \vert \boldsymbol{\theta}_i^y \big) \right)^T,$$
cf.~\eqref{eq:shape_gradient}, and exit the loop in case 
$$E^{\text{mod}}_n = \max \limits_{\beta_i^+ \in \partial \mathcal{T}_n^+} \left \| \nabla \mathcal{J}_h(\Omega_n, \beta_i^+) \right \|_2 < \varepsilon.$$
If $E_n^{\text{mod}} \geq \varepsilon$, we select the basis functions $\beta^+ \in \partial \mathcal{T}_n^+$ according to~\eqref{eq:selection_criterion} and we set $\mathcal{Q}_{n+1}$ to be the mesh obtained after refining all $Q \in \mathcal{Q}_n$ with $Q \in \mathcal{Q}_{\text{min}}(\mathcal{D}_{\text{bdry}}^+)$ before moving to the next iteration of the loop. \\

All steps are summarised in the pseudocode from Algorithm \ref{algo:simplification}. The algorithm has been implemented in the open-source finite element library \textit{Nutils} \cite{nutils7}.

\begin{algorithm}
	\caption{}
	\label{algo:simplification}
	\small{
	\begin{algorithmic}[1]
		
		\Procedure{Construct simplified geometry $\Omega_S$}{exact correspondence $\mathbf{F}: \partial \widehat{\Omega} \rightarrow \Omega$, functional $\mathcal{J}$, coarse initial mesh $\mathcal{Q}_0$, tolerance $\varepsilon$, marking parameter $\alpha$}
		\Statex
		\For{$n = 0, 1, \ldots$}
		\State $\mathcal{T}_n \longleftarrow$ canonical THB-spline basis associated with $\mathcal{Q}_n$.
		\State $\mathbf{F}_n \longleftarrow$ best approximation of $\mathbf{F}$ over $\partial \mathcal{T}_n \subset \mathcal{T}_n$ that follows from minimizing~\eqref{eq:boundary_fit}. \Comment{Fit}
		\State $\mathbf{x}_n \longleftarrow$ harmonic map of $\Omega_n$ acquired by solving~\eqref{eq:EGG_discretized_newton} with Dirichlet data $\mathbf{F}_n$.
             \State $\left( \mathcal{Q}_n^{\mathbf{x}}, \mathcal{T}_n^{\mathbf{x}} \right) \longleftarrow$ mesh and basis pair representing the bijective harmonic map $\mathbf{x}_n$.
		\State $\left( \mathcal{Q}_n^u, \mathcal{T}_n^u \right) \longleftarrow$ mesh and basis tuple, refinement of $\left( \mathcal{Q}_n^{\mathbf{x}}, \mathcal{T}_n^{\mathbf{x}} \right)$, used to solve for $u_h$.
		\State $\left( \mathcal{Q}_n^p, \mathcal{T}_n^p \right) \longleftarrow$ mesh and basis tuple, refinement of $\left( \mathcal{Q}_n^{\mathbf{x}}, \mathcal{T}_n^{\mathbf{x}} \right)$, used to solve for $p_h$.
		\State $(u_h^n, p_h^n) \longleftarrow$ weak state and adjoint equation solutions over $\mathcal{T}_n^u$ and $\mathcal{T}_n^p$, respectively. \Comment{Solve}
		\State $J_n \longleftarrow$ cost function evaluation that results from substituting $(u_h^n, \Omega_n)$ in~\eqref{eq:J_integral}.
		\State $(\mathcal{Q}_n^+, \mathcal{T}_n^+) \longleftarrow$ mesh and basis obtained after refining all $Q \in \mathcal{Q}_n$ with $\overline{Q} \cap \partial \widehat{\Omega} \neq \emptyset$.
		\State $\mathbf{F}_n^+ \longleftarrow$ best approximation of $\mathbf{F}$ obtained after minimizing~\eqref{eq:boundary_fit} over $\partial \mathcal{T}^+_n$.
		\State \smash{$\left\{ \Delta \mathbf{c}_i^+ \right\}_{i = 1}^{\partial N^+} \longleftarrow$} control points that follow from expressing $\Delta \mathbf{F}_n = \mathbf{F}_n^+ - \mathbf{F}_n$ in $\operatorname{span} \partial \mathcal{T}^+_n$.
		\State $\Theta_n \longleftarrow$ directions of differentiation associated with \smash{$\left\{ \Delta \mathbf{c}_i^+ \beta_i^+ \right\}_{i = 1}^{\partial N_n^+}$}, cf. \eqref{eq:Theta_n}.
		\State $\{\nabla \mathcal{J}_h(\Omega_n, \beta_k^+) \, \vert \, \beta_k^+ \in \partial \mathcal{T}_n^+ \} \longleftarrow$ approximate shape gradients resulting from $(u_h^n, p_h^n)$ and $\Theta_n$.
		\Statex
		\If{$E^{\text{mod}}_n = \max \limits_{\beta_k^+ \in \partial \mathcal{T}_n^+} \left \| \nabla \mathcal{J}_h(\Omega_n, \beta_k^+) \right \|_2 < \varepsilon$}
		    \Return{$\Omega_S = \Omega_n$.} \Comment{Estimate}
		\EndIf
		\Statex
		\State $\mathcal{D}_{\text{bdry}}^+ \longleftarrow$ set of active $\beta^+ \in \partial \mathcal{T}_n^+$ verifying \eqref{eq:selection_criterion}. \Comment{Mark}
		\State $\mathcal{Q}_{n+1} \longleftarrow$ new mesh acquired upon refining all $Q \in \mathcal{Q}_{\text{min}}(\mathcal{D}_{\text{bdry}}^+)$ according to~\eqref{eq:min_refinement}. \Comment{Refine}
		\Statex
		\EndFor
		\EndProcedure{}
		
	\end{algorithmic}
	}
\end{algorithm}

\newpage 
\section{Numerical experiments}\label{sec:numerical}
% say that the estimator is good since it has the same behaviour of the error and the multiplicative constant is not large
This section presents some numerical experiments which have been carefully selected in order to demonstrate the capabilities of the proposed algorithm as well as highlighting the pitfalls that may prevent full automation.

\subsection{A simple validation test case}
\label{sec:numerical_experiments_validation}
\begin{figure}[h!]
\centering
    \includegraphics[width=0.65\linewidth]{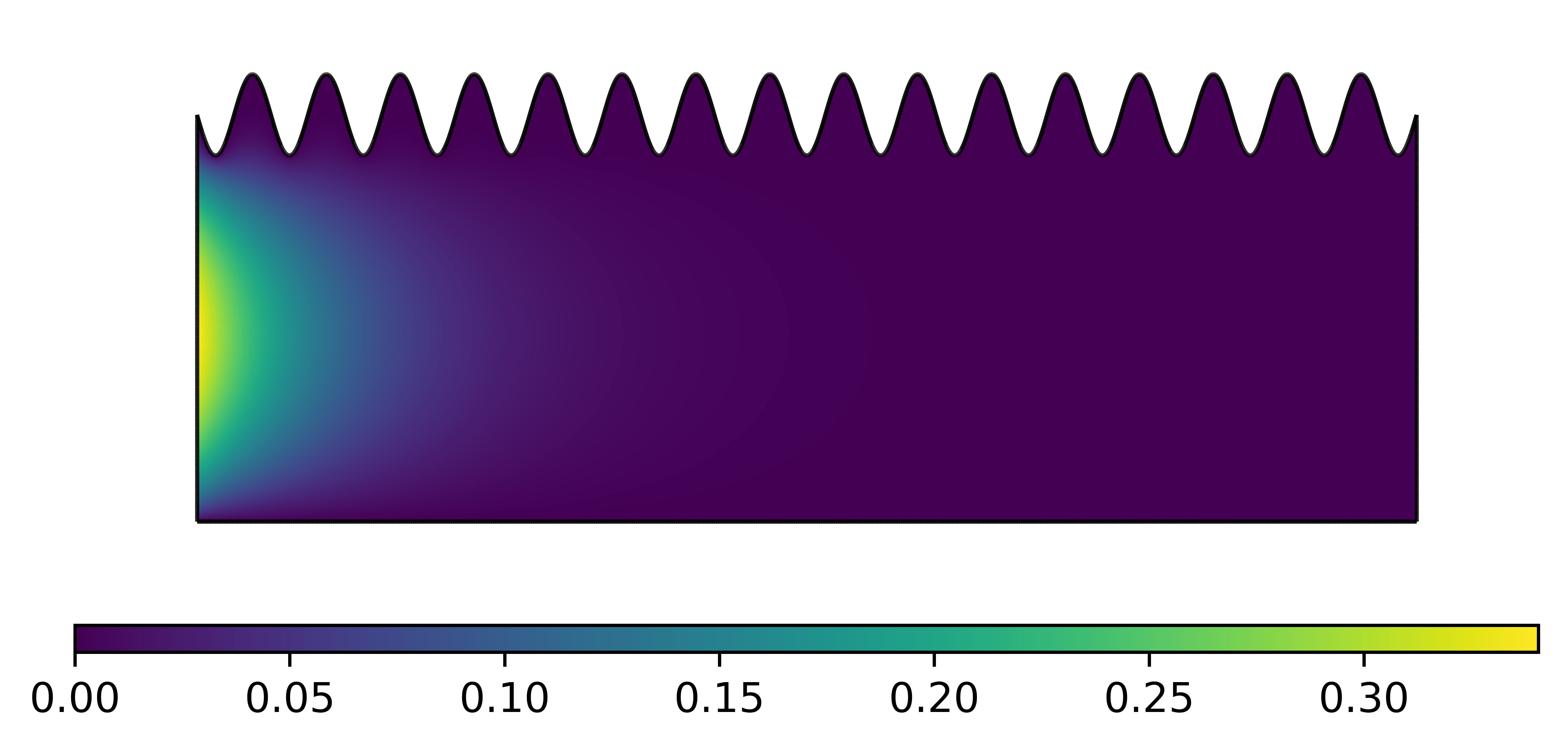}
\caption{Accurate representation of the geometry resulting from (\ref{eq:exact_bc_flag}) showing the state equation solution $u_E$.}
\label{fig:flag_exact}
\end{figure}

To experiment with the proposed methodology, we start with a simple test case that serves the purpose of gauging the plausibility of the results. Furthermore, we investigate the effect that the choice of the marking parameter $\alpha$ in (\ref{eq:selection_criterion}) has on the convergence behaviour. We are considering the flag-shaped domain from Figure \ref{fig:flag_exact}. We divide $\partial \widehat{\Omega}$ into its four constituents, i.e., $\partial \widehat{\Omega} = \gamma_S \cup \gamma_E \cup \gamma_N \cup \gamma_W$, where the $\gamma_\beta, \beta \in \{S, E, N, W\}$ denote its southern, eastern, northern and western parts, respectively. The exact boundary correspondence $\mathbf{F}: \partial \widehat{\Omega} \rightarrow \Omega$ reads
\begin{align}
\label{eq:exact_bc_flag}
    \mathbf{F}(\boldsymbol{\xi}) = \left \{ \begin{array}{ll} (3 \xi_1, 0)^T & \boldsymbol{\xi} \in \gamma_S \\
                                                              (3, \xi_2)^T & \boldsymbol{\xi} \in \gamma_E \\
                                                              (3 \xi_1, 1 - 0.1 \sin{(31 \pi \xi_1)})^T & \boldsymbol{\xi} \in \gamma_N \\
                                                              (0, \xi_2)^T & \boldsymbol{\xi} \in \gamma_W. \end{array} \right.
\end{align}
The state equation results from the model problem (\ref{eq:weak_state_pb}) with $f=0$, $g=1$ and $\Gamma_N = \{0\} \times (0, 1)=\gamma_W$ (i.e., the western part of $\partial \Omega$) while $\Gamma_D = \partial \Omega \setminus \Gamma_N$. The cost function $\mathcal{J}(\Omega)$ integrates the square of the state equation's solution over $\gamma_W$, i.e.,
\begin{align*}
    J(u(\Omega); \Omega) = \int \limits_{\gamma_W} u^2 \mathrm{d} x.
\end{align*}
We start with a uniform open knot vector comprised of $10$ elements in either coordinate direction, which amounts to $9$ internal knots and results in the initial structured mesh $\mathcal{Q}_0$ of $\widehat{\Omega}$ comprised of $100$ elements. There are no internal knot repetitions and throughout this example, we employ bicubic THB-splines with maximum regularity. The exact boundary is expressed in a THB-spline basis by dyadically refining all elements $Q \in \mathcal{Q}$ with $\overline{Q} \cap \overline{\gamma}_N \neq \emptyset$. This refinement is repeated a total of seven times, resulting in a total of $\partial N_{\Omega} = 1332$ degrees of freedom (DOFs) that are nonvanishing on $\partial \widehat{\Omega}$. The exact boundary correspondence is then projected onto the corresponding spline space's restriction to the boundary $\mathcal{T} \vert_{\partial \widehat{\Omega}}$ by minimising~\eqref{eq:boundary_fit} with $\kappa_0 = \kappa_1 = 1$ while constraining the fit to the boundary correspondence's exact value by the four corners of $\partial \widehat{\Omega}$. The corresponding projection error is negligible and in the following, the spline fit will be regarded as the new exact boundary correspondence. We then approximate a harmonic map (see Section \ref{sect:EGG}) $\mathbf{x}_E: \widehat{\Omega} \rightarrow \Omega$ using the boundary correspondence as the Dirichlet data while performing \textit{a posteriori} refinement if necessary. 

\noindent The computed exact mapping is used to numerically solve the weak state equation, yielding the approximation $u_E \in H^1_{0, \Gamma_D}(\Omega)$. The substitution into the cost function yields $J_E$ which is used as a reference value to compare the effect of defeaturing the geometry. We therefore compute $u_E$ very accurately, in order to largely isolate the geometric component on the cost function discrepancy. Here, $u_E$ is computed over a THB-spline basis that results from uniform dyadic refinement of $\mathcal{Q}_0$ a total of four times in both coordinate directions, leading to a basis comprised of $26569$ DOFs. The reference solution is plotted in Figure \ref{fig:flag_exact}. \\

\noindent It should be noted that $\mathbf{F}: \partial \widehat{\Omega} \rightarrow \Omega$ restricted to $\gamma_W \cup \gamma_S \cup \gamma_E$ is exactly contained in all bicubic THB-spline bases defined over $\widehat{\Omega}$ (thanks to the exact correspondence's piecewise linear nature on this subset). Hence, all the shape derivatives associated with boundary control points that correspond to this subset are zero. For $\alpha >0$, we therefore expect the algorithm to only mark elements $Q \in \mathcal{Q}$ with $\overline{Q} \cap \gamma_N \neq \emptyset$ for refinement. The reference solution from Figure \ref{fig:flag_exact} furthermore reveals that large values of $u_E$ are mainly assumed close to the Neumann influx boundary $\Gamma_N$ with a steep decline in the direction of the eastern boundary. It is therefore plausible to assume that refinement will be concentrated along the western part of $\gamma_N$. Throughout the remainder of this example, for fitting the $n$-th boundary correspondence $\mathbf{F}_n: \partial \widehat{\Omega} \rightarrow \partial \Omega_n$ to $\mathbf{F}: \partial \widehat{\Omega} \rightarrow \partial \Omega$, we  minimise~\eqref{eq:boundary_fit} with $\kappa_0 = \kappa_1 = 1$, as before.\\

\noindent As a first example, we study the convergence behaviour for a very small value $\alpha = 10^{-7}$. Clearly, with this choice, we expect the algorithm to over-refine. We do not take $\alpha = 0$ in order to avoid refinement on $\partial \widehat{\Omega} \setminus \gamma_N$. As a measure of the approximation quality, we utilise the relative (absolute) discrepancy between the cost function evaluation at the current iteration $J_n$ and the reference value $J_E$ while taking the value $\partial N_{\Omega_n} / \partial N_{\Omega}$, wherein $\partial N_{\Omega_n}$ denotes the total number of boundary control points at the current iteration, as a measure of the geometrical DOF savings with respect to the reference solution. Furthermore, we also monitor the value of $E^{\text{mod}}_n(\Omega_n)$ relative to $J_E \approx 6.50273 \times 10^{-2}$ at each iteration. The algorithm is terminated once $E^{\text{mod}}_n(\Omega_n) < 10^{-5}$ and we utilise the same (dense) THB-spline basis $\mathcal{T}$ we used to compute $u_E$ to solve the state and adjoint equations at each iteration.\\

\begin{figure}[h!]
\centering
\includegraphics[width=0.7 \textwidth]{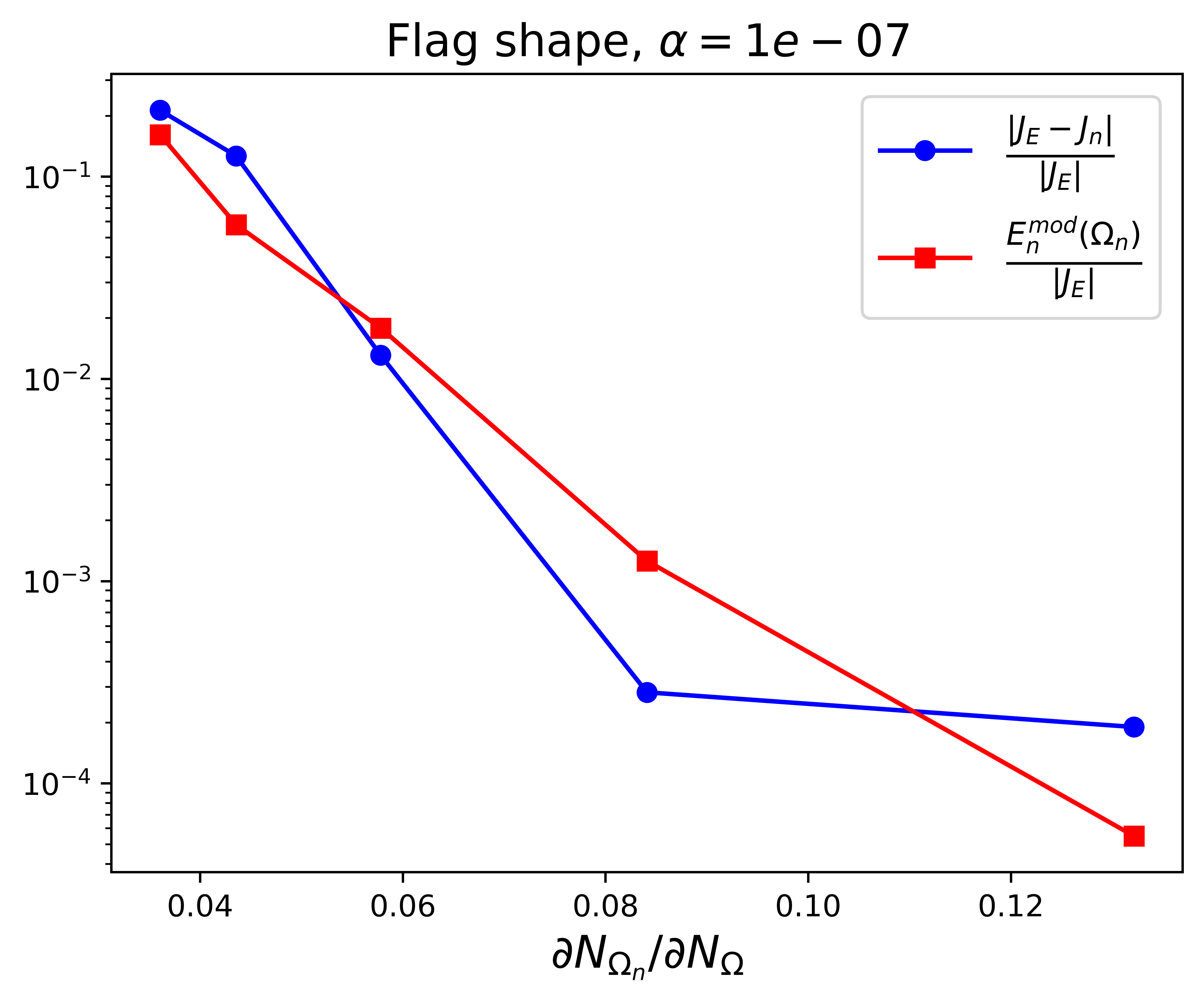}
\caption{The convergence plot of the flag example for $\alpha = 10^{-7}$ and convergence threshold $E^{\text{mod}}_n(\Omega_n) < 10^{-5}$.}
\label{fig:flag_convergence_0}
\end{figure}

\begin{figure}[h!]
\centering
\captionsetup[subfigure]{labelformat=empty}
    \begin{subfigure}[t]{0.9\textwidth}
        \centering
        \includegraphics[align=c, width=0.28\linewidth]{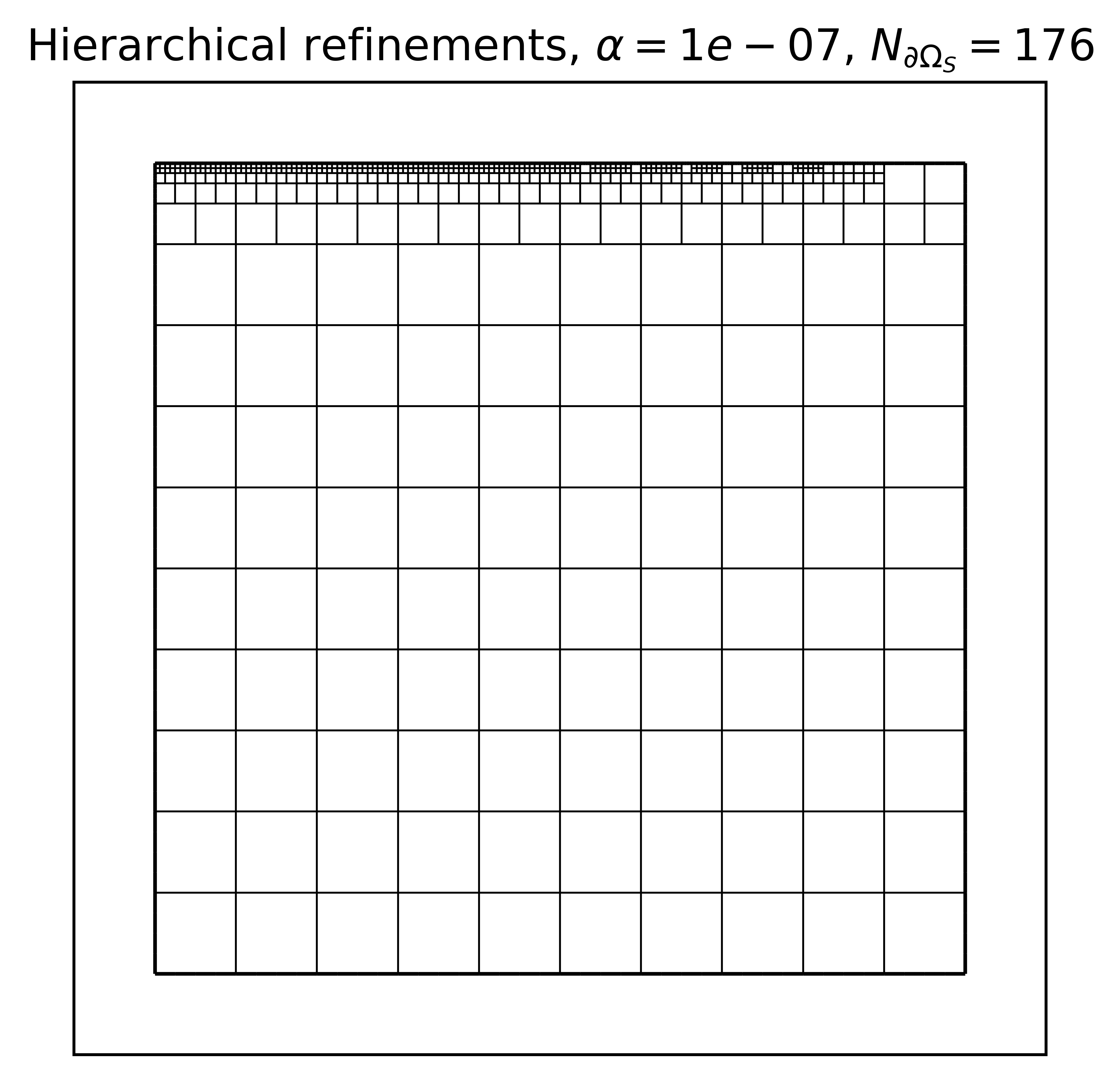} $\quad$ \includegraphics[align=c, width=0.55\linewidth]{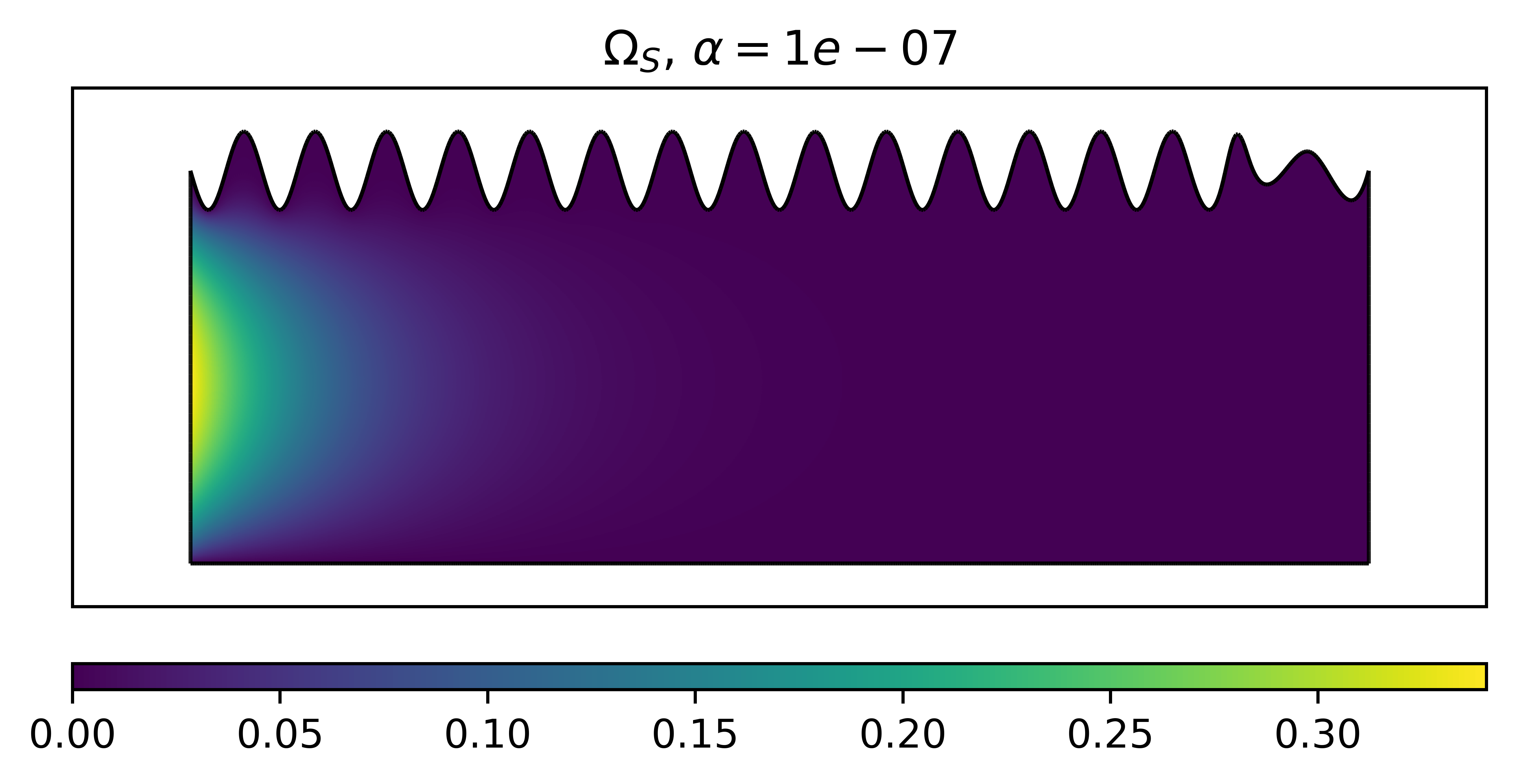}
    \end{subfigure}
\caption{The hierarchical refinements in $\widehat{\Omega}$ and the defeatured domain $\Omega_S$ after the last iteration.}
\label{fig:refinements_geometry_0}
\end{figure}

\noindent The corresponding convergence plot is depicted in Figure \ref{fig:flag_convergence_0}, while Figure \ref{fig:refinements_geometry_0} shows $\Omega_S$ and the hierarchical refinements in $\widehat{\Omega}$. Convergence is reached with $\partial N_{\Omega_S} = 176$, indicating a reduction of the boundary degrees of freedom of almost 90\% while having a relative modelling error around $10^{-4}$. The algorithm has refined a large portion of $\gamma_N$ several times and over-refinement is visible in Figure \ref{fig:refinements_geometry_0} as the geometry becomes virtually indistinguishable from the reference geometry depicted in Figure \ref{fig:flag_exact}, apart from the eastern part of $\gamma_N$ which corresponds to the last period in the sine's cycle. This was expected since the marking parameter $\alpha$ is taken very small, and thus the refinement should almost be uniform. Note that only the boundary $\gamma_N$ has been refined, as expected: indeed, more degrees of freedom on the other parts of boundary will not improve the geometric approximation and thus it will not improve neither the accuracy of the solution. The fact that the eastern part of $\gamma_N$ is the only defeatured part shows that the algorithm constructing $\Omega_S$ is indeed driven by the analysis behind, and not only by geometrical considerations. Indeed, the eastern part of $\Omega_N$ is the region which is the furthest from the domain $\gamma_W$ on which the objective function $J$ is defined, and it is also the region where the solution has a larger gradient. We were therefore expecting this result, which shows a good geometric approximation close to $\gamma_W$, and a defeaturing effect far away from $\gamma_W$. 
\\

\begin{figure}[h!]
\centering
\captionsetup[subfigure]{labelformat=empty}
    \begin{subfigure}[t]{0.45\textwidth}
        \centering
        \includegraphics[align=c, width=0.9\linewidth]{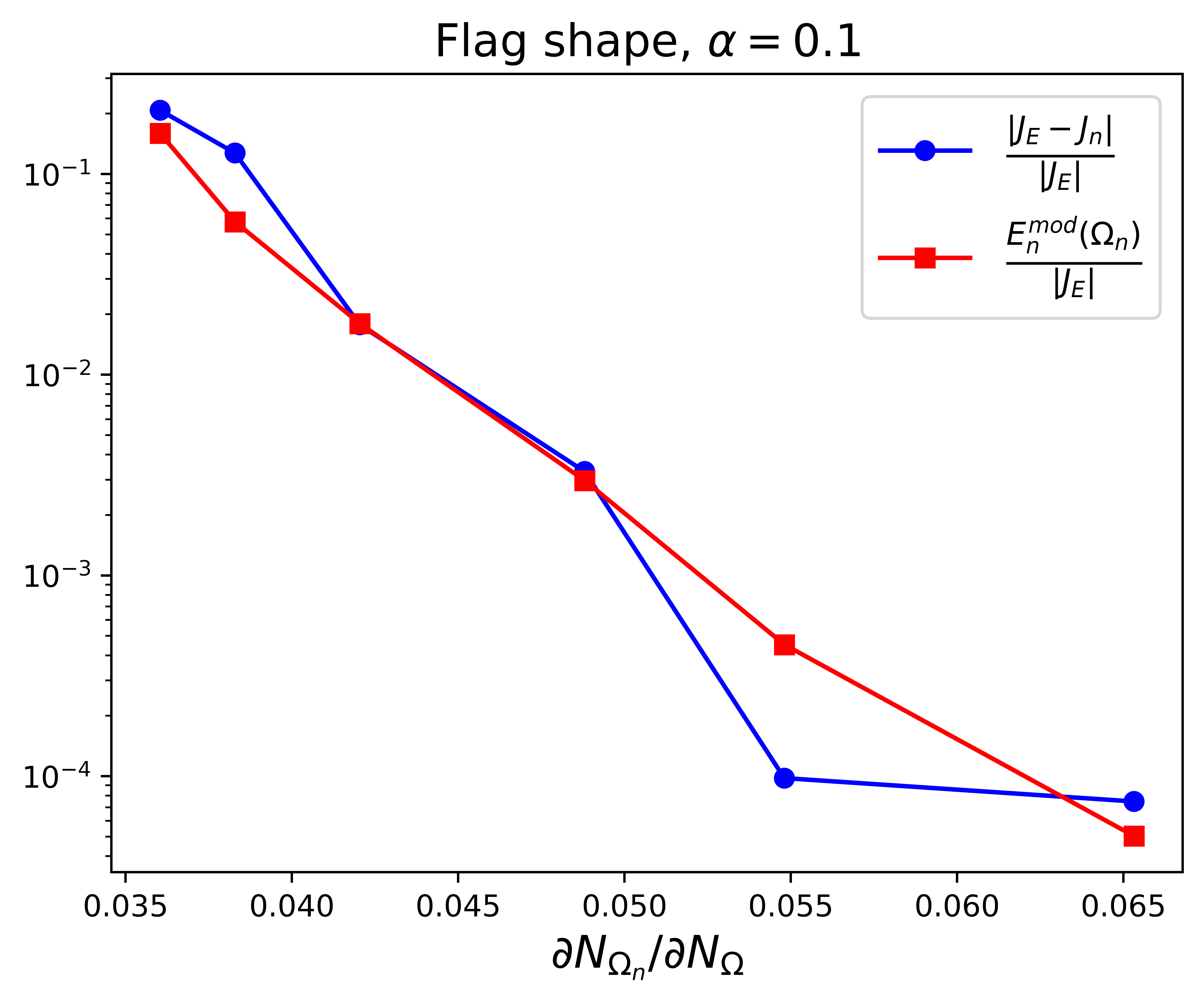}
    \end{subfigure} $\quad$
    \begin{subfigure}[t]{0.45\textwidth}
        \centering
        \includegraphics[align=c, width=0.9\linewidth]{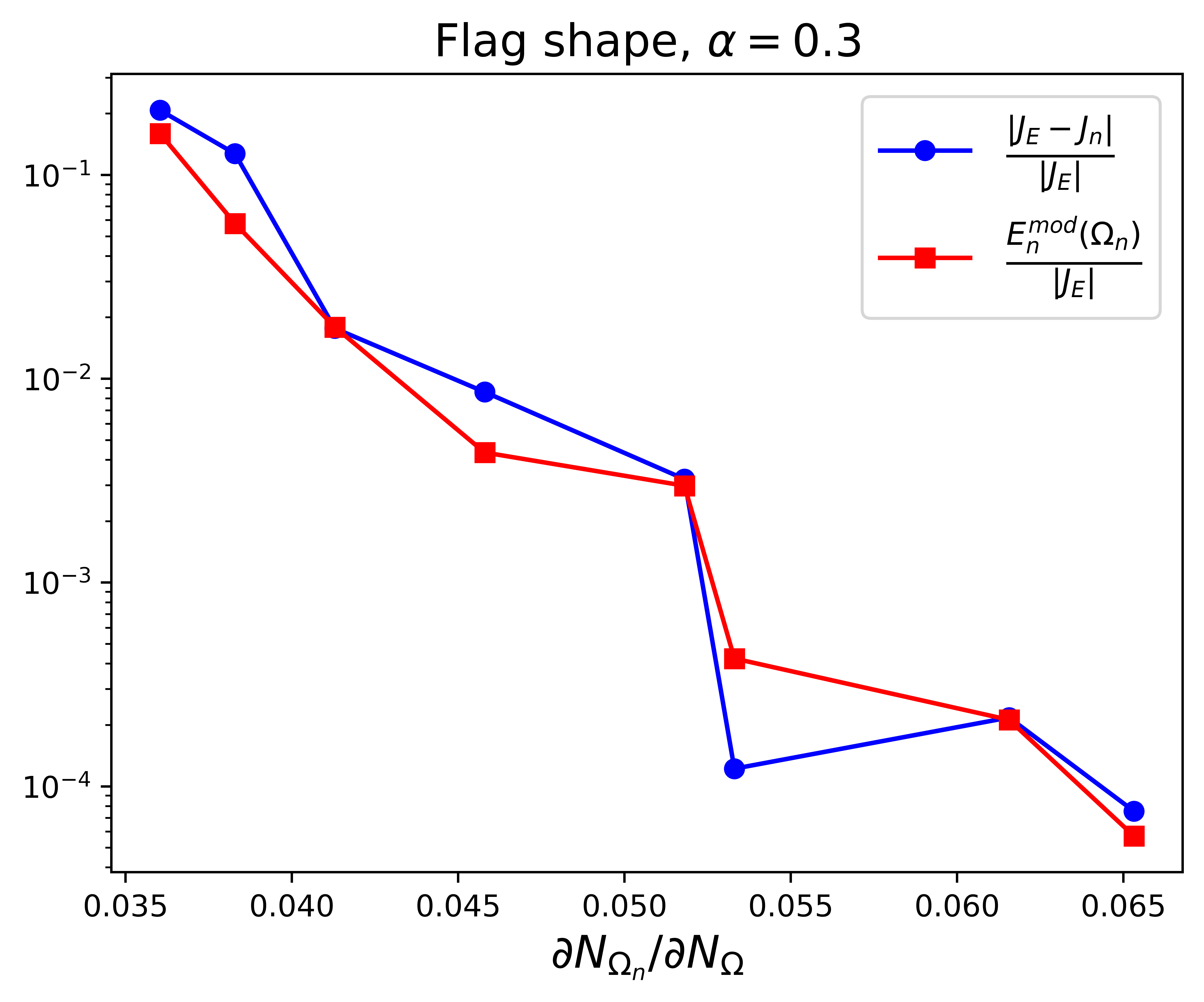}
    \end{subfigure} \\
    \begin{subfigure}[t]{0.45\textwidth}
        \centering
        \includegraphics[align=c, width=0.9\linewidth]{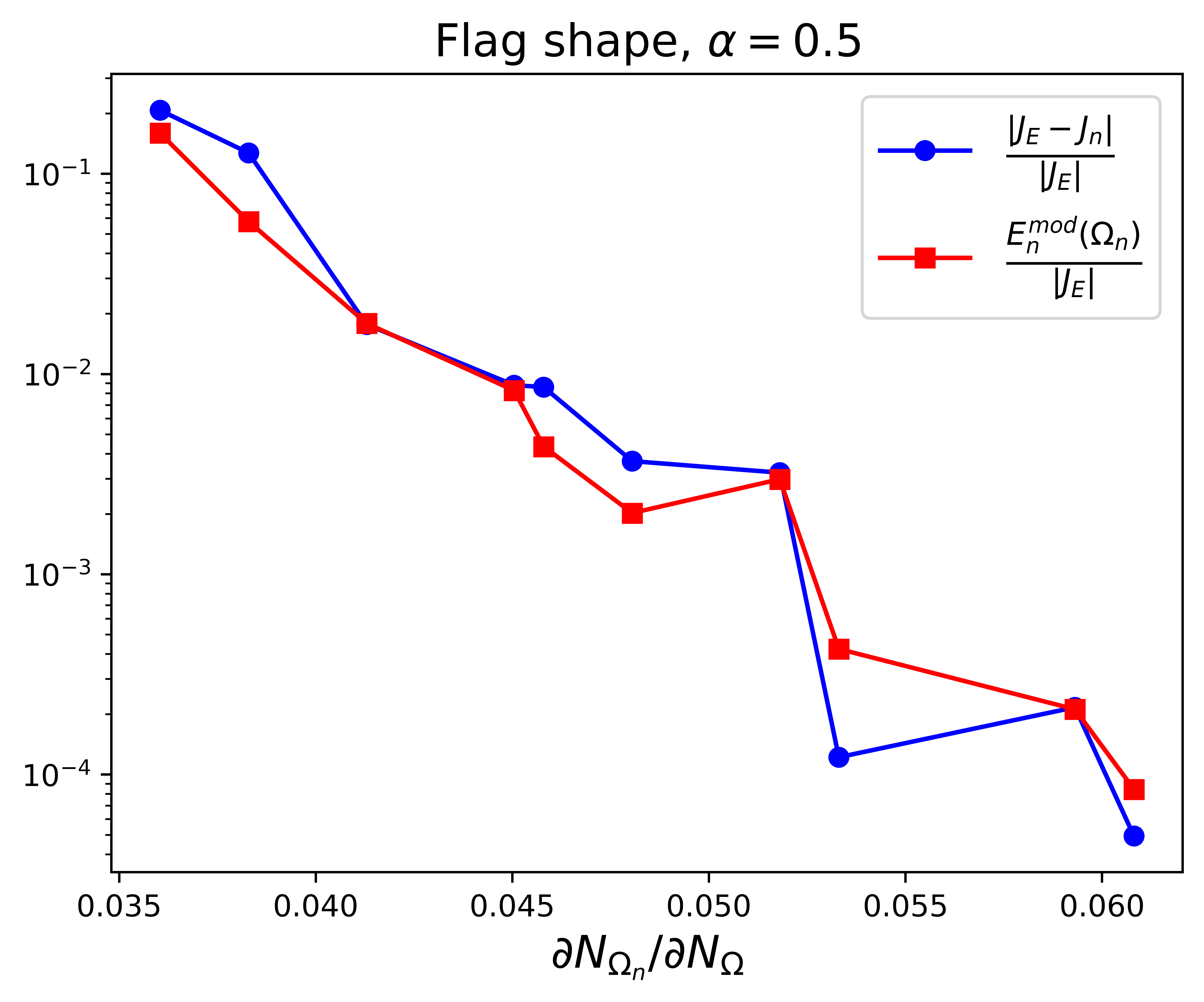}
    \end{subfigure} $\quad$
    \begin{subfigure}[t]{0.45\textwidth}
        \centering
        \includegraphics[align=c, width=0.9\linewidth]{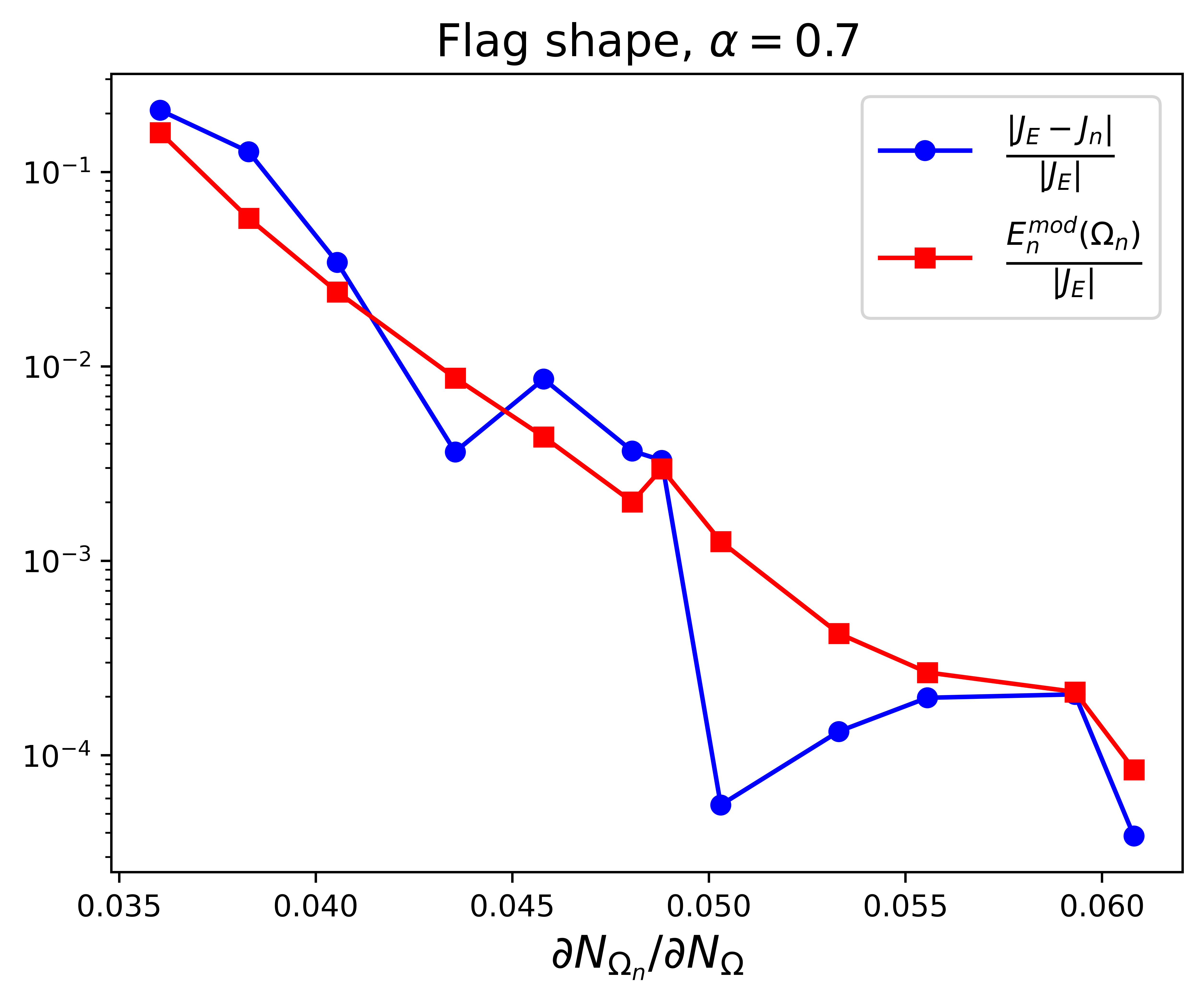}
    \end{subfigure}
\caption{Convergence plots for $\alpha \in \{0.1, 0.3, 0.5, 0.7\}$ and convergence threshold $E^{\text{mod}}_n(\Omega_n) < 10^{-5}$.}
\label{fig:convergence_plots_flag_shape}
\end{figure}

\begin{figure}[h!]
\centering
\captionsetup[subfigure]{labelformat=empty}
    \begin{subfigure}[t]{0.9\textwidth}
        \centering
        \includegraphics[align=c, width=0.28\linewidth]{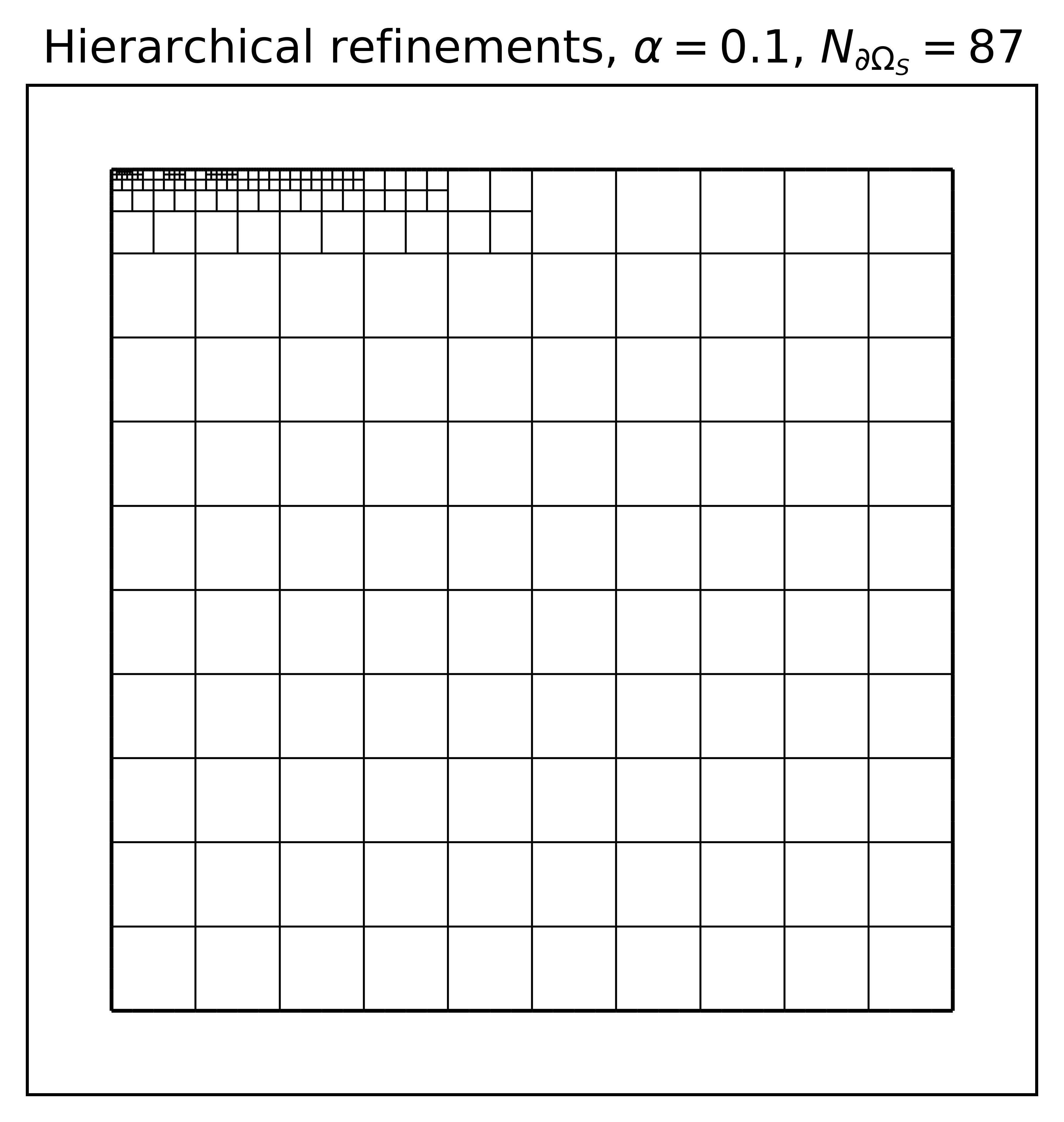} $\quad$ \includegraphics[align=c, width=0.55\linewidth]{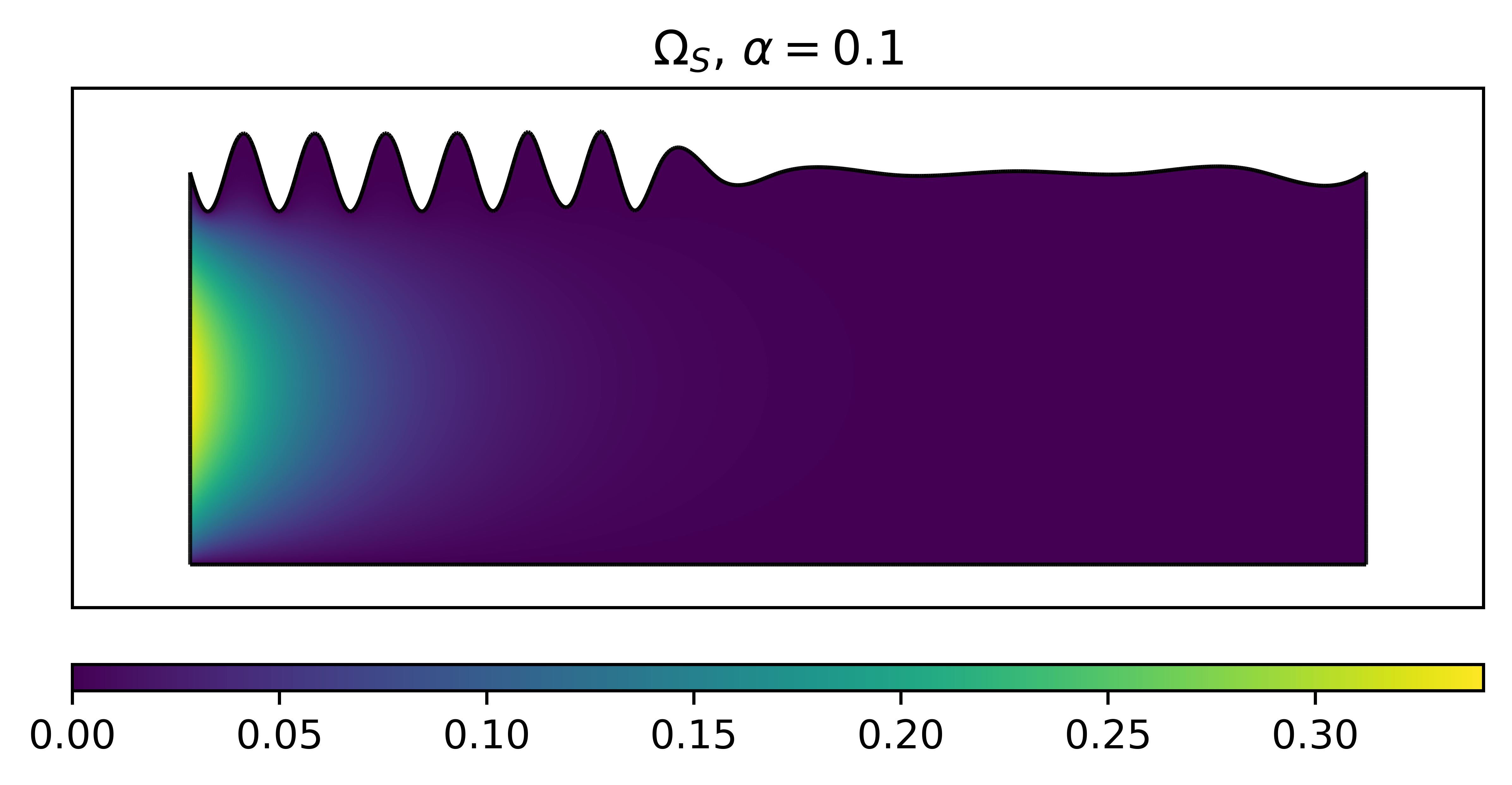}
    \end{subfigure} \\
    % \begin{subfigure}[t]{0.9\textwidth}
    %     \centering
    %     \includegraphics[align=c, width=0.28\linewidth]{defeaturing_shape_der/img/num_test/flag/domain_alpha_0.3.png} $\quad$ \includegraphics[align=c, width=0.55\linewidth]{defeaturing_shape_der/img/num_test/flag/Omega_alpha_0.3.png}
    % \end{subfigure} \\
    % \begin{subfigure}[t]{0.9\textwidth}
    %     \centering
    %     \includegraphics[align=c, width=0.28\linewidth]{defeaturing_shape_der/img/num_test/flag/domain_alpha_0.5.png} $\quad$ \includegraphics[align=c, width=0.55\linewidth]{defeaturing_shape_der/img/num_test/flag/Omega_alpha_0.5.png}
    % \end{subfigure} \\
    \begin{subfigure}[t]{0.9\textwidth}
        \centering
        \includegraphics[align=c, width=0.28\linewidth]{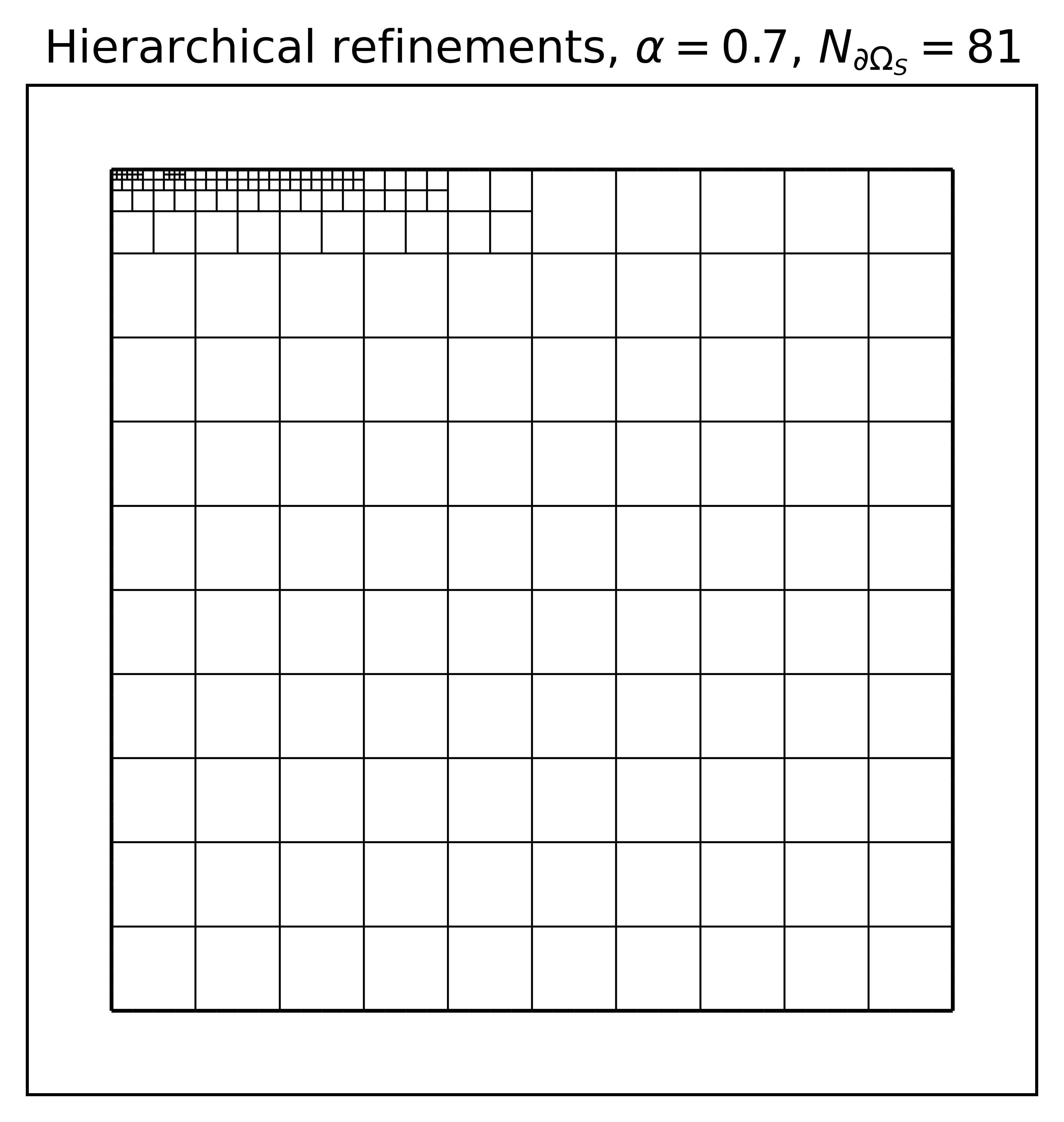} $\quad$ \includegraphics[align=c, width=0.55\linewidth]{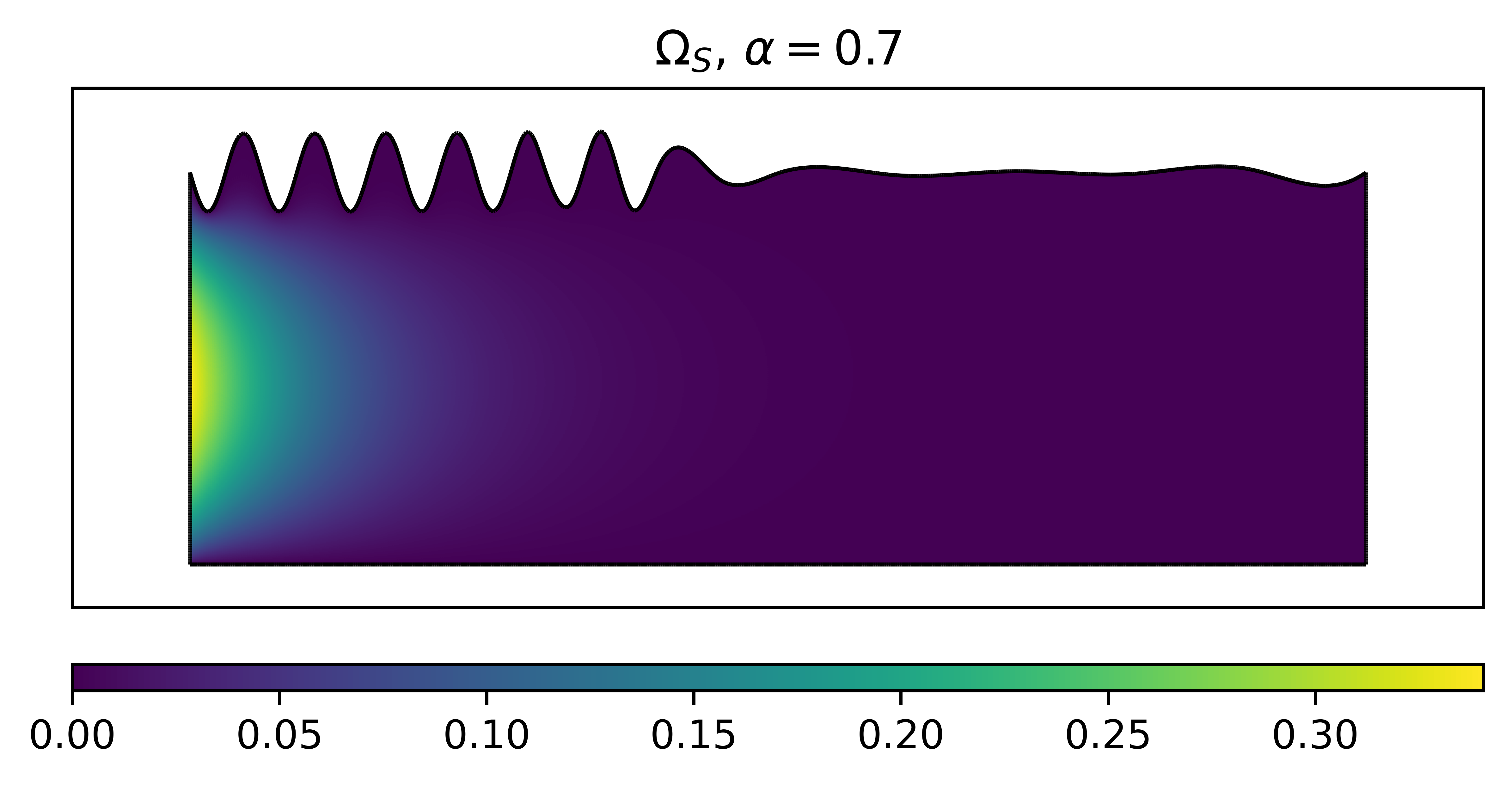}
    \end{subfigure}
\caption{Defeatured domains along with the corresponding boundary refinements in $\widehat{\Omega}$ for $\alpha \in \{0.1, 0.7 \}$.}
\label{fig:flag_solutions}
\end{figure}

\noindent For larger values of $\alpha$, we expect the convergence criterion to be reached with fewer boundary DOFs while potentially requiring more iterations in the case of under-refinement. Figure \ref{fig:convergence_plots_flag_shape} shows the convergence plots for $\alpha \in \{0.1, 0.3, 0.5, 0.7\}$ while Figure \ref{fig:flag_solutions} shows the defeatured domains after the last iteration along with the boundary refinements in $\widehat{\Omega}$ for $\alpha \in \{0.1, 0.7\}$.\\

\noindent The figures confirm our expectations as the required number of iterations range from $n_{\text{iter}} = 6$ for $\alpha = 0.1$ to $n_{\text{iter}}= 12$ for $\alpha = 0.7$ with a mild reduction in the required number of boundary DOFs which are given by $\partial N_{\Omega_S}(\alpha=0.1) = 87$, $\partial N_{\Omega_S}(\alpha=0.3) = 87$, $\partial N_{\Omega_S}(\alpha=0.5) = 81$ and $\partial N_{\Omega_S}(\alpha=0.7) = 81$. Overall, the convergence plot for $\alpha = 0.1$ exhibits the most monotone behaviour and requires the fewest number of iterations at the expense of only a slight increase in the number of boundary DOFs. This suggest that $\alpha = 0.1$ constitutes a good trade off between the required number of boundary DOFs and the total number of iterations.

Furthermore, the figures clearly show that $E^{\text{mod}}_n(\Omega_n) / \vert J_e \vert$ and the relative cost function discrepancy follow the same trend. Note that a monotone convergence cannot be expected since the taken descent direction (in which only some boundary elements are refined) is different from the direction in which the shape gradients are calculated, since the shape gradients are obtained with the direction computed from the geometry in which all boundary elements are refined. Indeed, it is otherwise absolutely computationally unfeasible to compute all possible descent directions obtained with the refinement of every possible subset of the boundary elements. A way to circumvent this would be to develop a more localized strategy in the FIT step presented in Section~\ref{subsec:algo_fit}; this will be discussed in more detail in the next experiment. However, we still obtain a remarkable reduction of degrees of freedom necessary to describe the geometry to obtain a relative error in the objective functional of around $10^{-4}$: indeed, the number of boundary degrees of freedom of $\Omega_S$ are reduced by 94\% compared to $\Omega$. Finally and ever more notably than with a larger $\alpha$, a good geometric approximation is obtained in $\Omega_S$ close to $\gamma_W$, and most of the defeaturing effect is obtained far away from $\gamma_W$, as expected. 

\subsection{A more advanced experiment}
\label{sec:numerical_experiments_advanced}
\begin{figure}[h!]
    \centering
    \includegraphics[width=0.8 \linewidth]{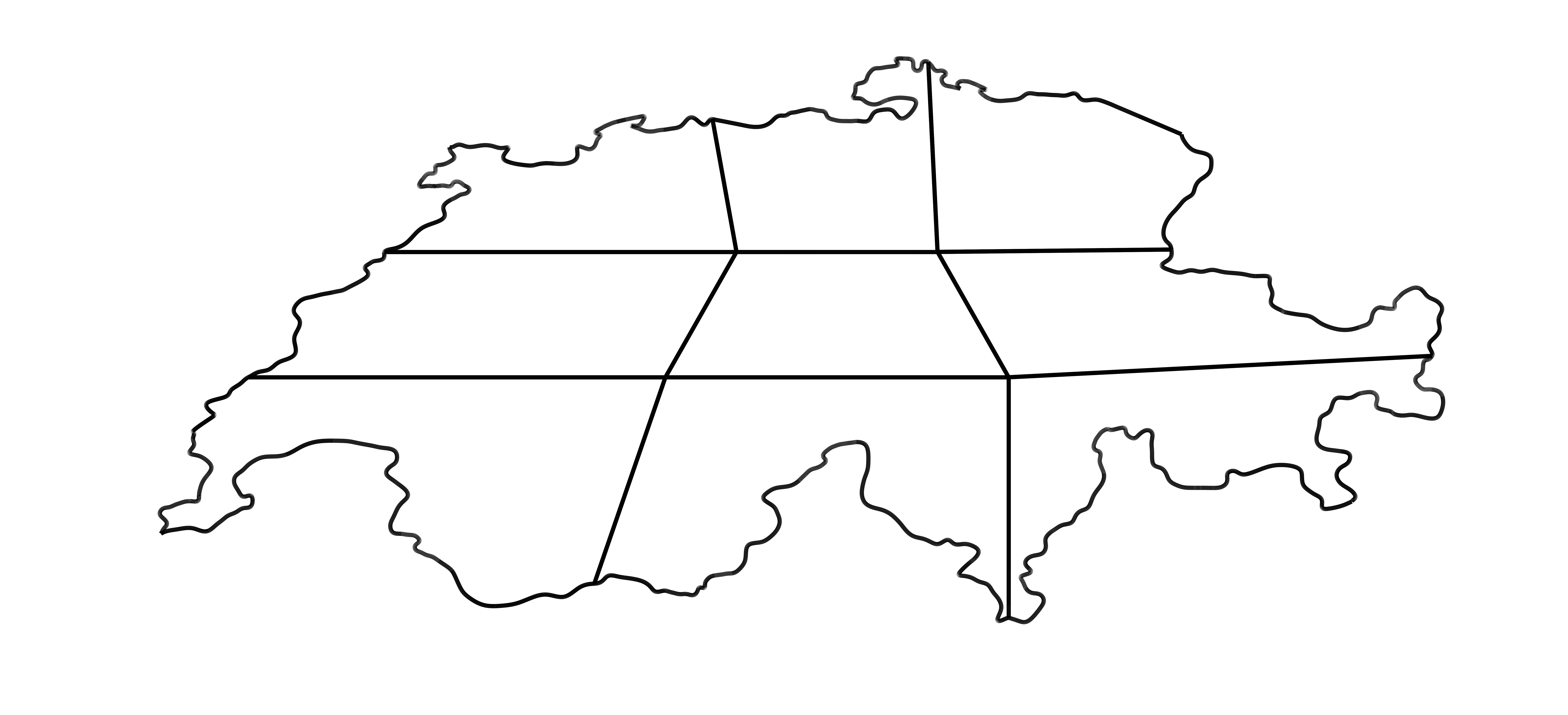}
    \caption{Plot showing the segmentation of the target geometry into nine parts.}
    \label{fig:9_patch_switz_contours.png}
\end{figure}

We are considering the domain $\Omega$ whose contours are given by a very accurate spline representation of the country of Switzerland, see Figure \ref{fig:9_patch_switz_contours.png}. The contours are comprised of $2086$ boundary DOFs. The initial mesh $\mathcal{Q}_0$ of $\widehat{\Omega}$ results from the initial uniform cubic knot vector $\Xi$ with spacing \smash{$\Delta \xi = \tfrac{1}{15}$} and a three-fold knot repetition at \smash{$\xi_1 = \tfrac{1}{3}$} and \smash{$\xi_2 = \tfrac{2}{3}$}. Using $\Xi$ for both coordinate directions naturally divides $\widehat{\Omega}$ into nine subsets $\widehat{\Omega}_i, \enskip i \in \{1, \ldots, 9 \}$ in the usual lexicographical ordering, also called patches in the isogeometric terminology, fenced off by $\partial \widehat{\Omega}$ and the repeated knots at $\xi_1$ and $\xi_2$. On each of the $\widehat{\Omega}_i$, any bicubic THB-spline basis $\mathcal{T}$, resulting from $\mathcal{Q}_0$ or a hierarchical refinement thereof, will satisfy $\mathcal{T} \vert_{\widehat{\Omega}_i} \subset H^2(\widehat{\Omega}_i)$ while the $\beta \in \mathcal{T}$ are interpolatory on the isolines with $i \in \{1, 2\}: \xi_i \in \left\{ 0, \tfrac{1}{3}, \tfrac{2}{3}, 1 \right \}$. \\

\noindent We are considering the model problem (\ref{eq:weak_state_pb}) with $f=0$, $g=100$ and $\Gamma_N = \mathbf{F} \vert_{\gamma_S}$ (i.e., the image of the southern part of $\partial \widehat{\Omega}$ under the boundary correspondence) while, again, $\Gamma_D = \partial \Omega \setminus \Gamma_N$. To each of the $\partial \widehat{\Omega}_i$, we assign its own boundary correspondence $\mathbf{F}_i$ which is comprised of a number of straight line segments in the interior of $\partial \Omega$ and the boundary correspondence $\mathbf{F}: \partial \widehat{\Omega} \rightarrow \partial \Omega$, see Figure \ref{fig:9_patch_switz_contours.png}.
At each iteration, we approximate a total of nine maps $\mathbf{F}_{n, i}: \widehat{\Omega}_i \rightarrow \Omega_{n, i}$ mapping each of the $\widehat{\Omega}_i$ onto the domains $\Omega_{n, i}$ which result from combining the interior line segments with the defeatured boundary correspondence. For this, we pass the defeatured boundary correspondences $\mathbf{F}_{n, i}: \partial \widehat{\Omega}_i \rightarrow \partial \Omega_{n, i}$, along with the restrictions $\mathcal{T} \vert_{\widehat{\Omega}_i} \subset H^2(\widehat{\Omega}_i)$, to the parameterisation routine from Section \ref{sect:EGG} one-by-one while performing \textit{a posteriori} refinement if necessary. Clearly, with this methodology, the centering subset $\mathcal{C} = \Omega_{n, 4} \subset \Omega_n$ is kept fixed and hence suitable for defining a cost function on it. The cost function reads:
\begin{align}
    J\left(u(\Omega); \, \Omega \right) = \int \limits_{\mathclap{\mathcal{C}}} u^2 \mathrm{d}\mathbf{x}.
\end{align}

\begin{figure}[h!]
\centering
\captionsetup[subfigure]{labelformat=empty}
    \begin{subfigure}[t]{0.45\textwidth}
        \centering
        \includegraphics[align=c, width=1\linewidth]{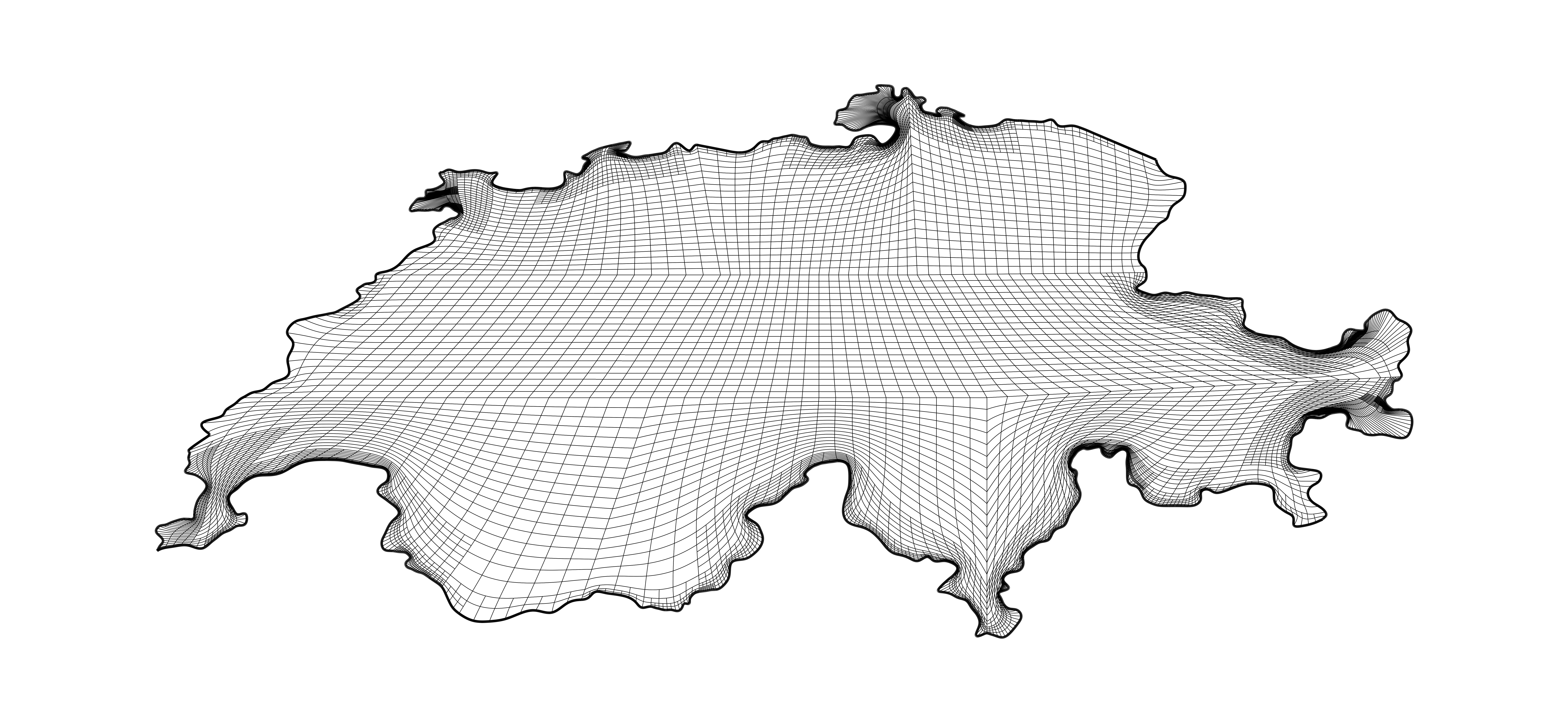}
    \end{subfigure} $\enskip$
    \begin{subfigure}[t]{0.45\textwidth}
        \centering
        \includegraphics[align=c, width=1\linewidth]{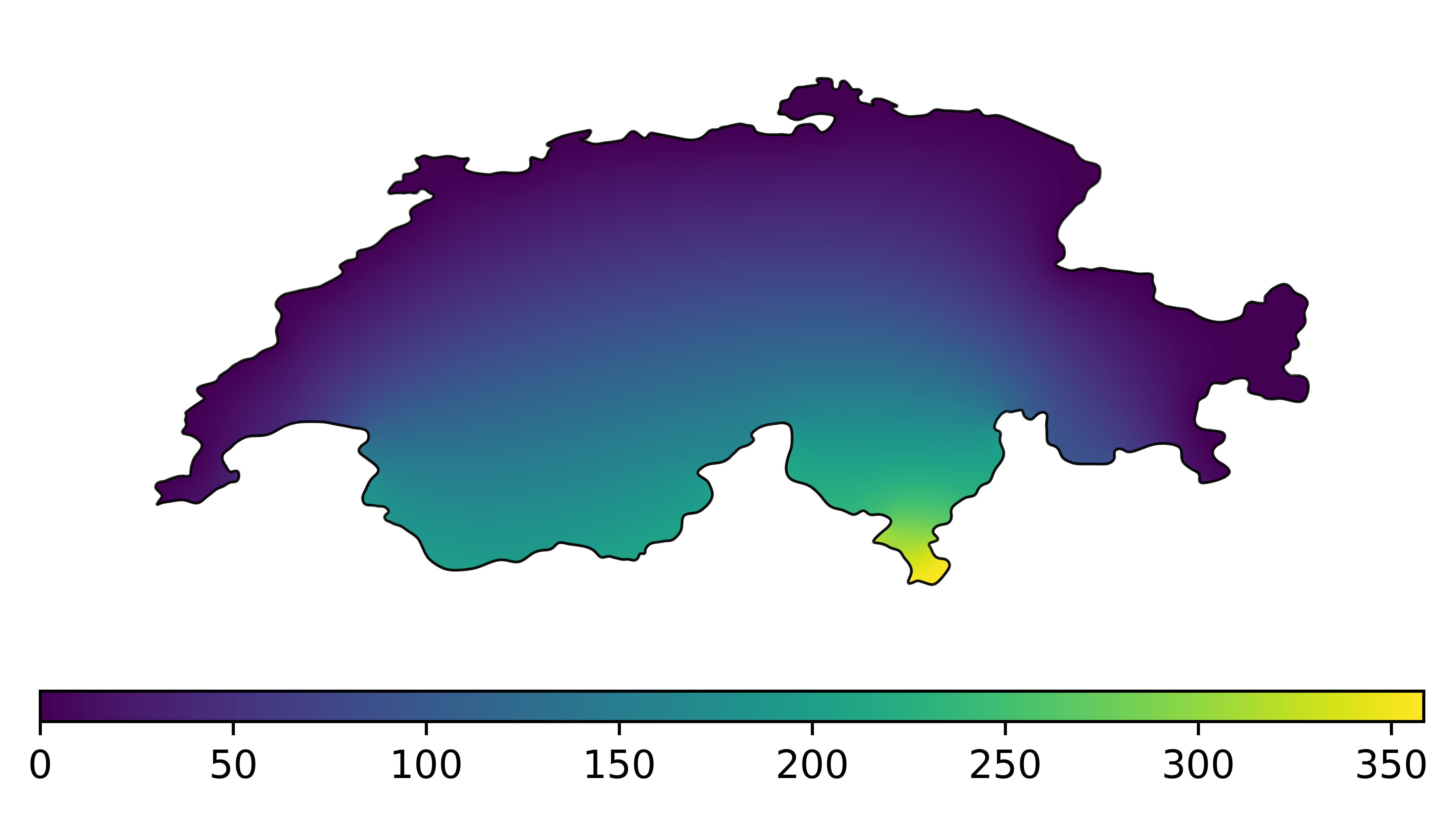}
    \end{subfigure}
\caption{Reference map (left) showing the inner refinements and the reference solution $u_E$ (right).}
\label{fig:swtizerland_exact}
\end{figure}

In order to monitor the convergence behaviour, as in the previous example, we compute an accurate reference solution by passing the exact correspondence $\mathbf{F}: \partial \widehat{\Omega} \rightarrow \partial \Omega$ to aforementioned routine to compute a map $\mathbf{x}: \widehat{\Omega} \rightarrow \Omega$. It is then utilised to compute $u_E$ and the associated $J_E$ over a basis comprised of $239576$ DOFs which results from uniformly refining $\mathcal{Q}_0$ several times. The map $\mathbf{x}: \widehat{\Omega} \rightarrow \Omega$ showing the refinements to acquire a folding-free map, along with the reference solution are plotted in Figure \ref{fig:swtizerland_exact}.

At each iteration the defeatured boundary correspondence $\mathbf{F}_n$ results from minimising~\eqref{eq:boundary_fit} with $\kappa_0 = \kappa_1 = 1$ over the current iteration's THB-spline basis $\mathcal{T}_n$ while constraining the fit as well as its tangent to the exact correspondence's values by the points $\boldsymbol{\xi} \in \partial \widehat{\Omega}$ with $\boldsymbol{\xi} \in \{0, \tfrac{1}{3}, \tfrac{2}{3}, 1\} \times \{0, \tfrac{1}{3}, \tfrac{2}{3}, 1\}$. These additional constraints avoid accidentally creating concave corners in the approximation of the piecewise smooth exact correspondence. In rare cases we have noticed the fit to lead to self-intersections in the boundary correspondence. In case self-intersections are detected, we minimally refine all $\beta \in \partial \mathcal{T}_n$ that are nonvanishing on the intersection point and perform the fit again until no more self-intersections are found. \\

\begin{figure}[h!]
\centering
\captionsetup[subfigure]{labelformat=empty}
    \begin{subfigure}[t]{0.9\textwidth}
        \centering
        \includegraphics[align=c, width=0.35\linewidth]{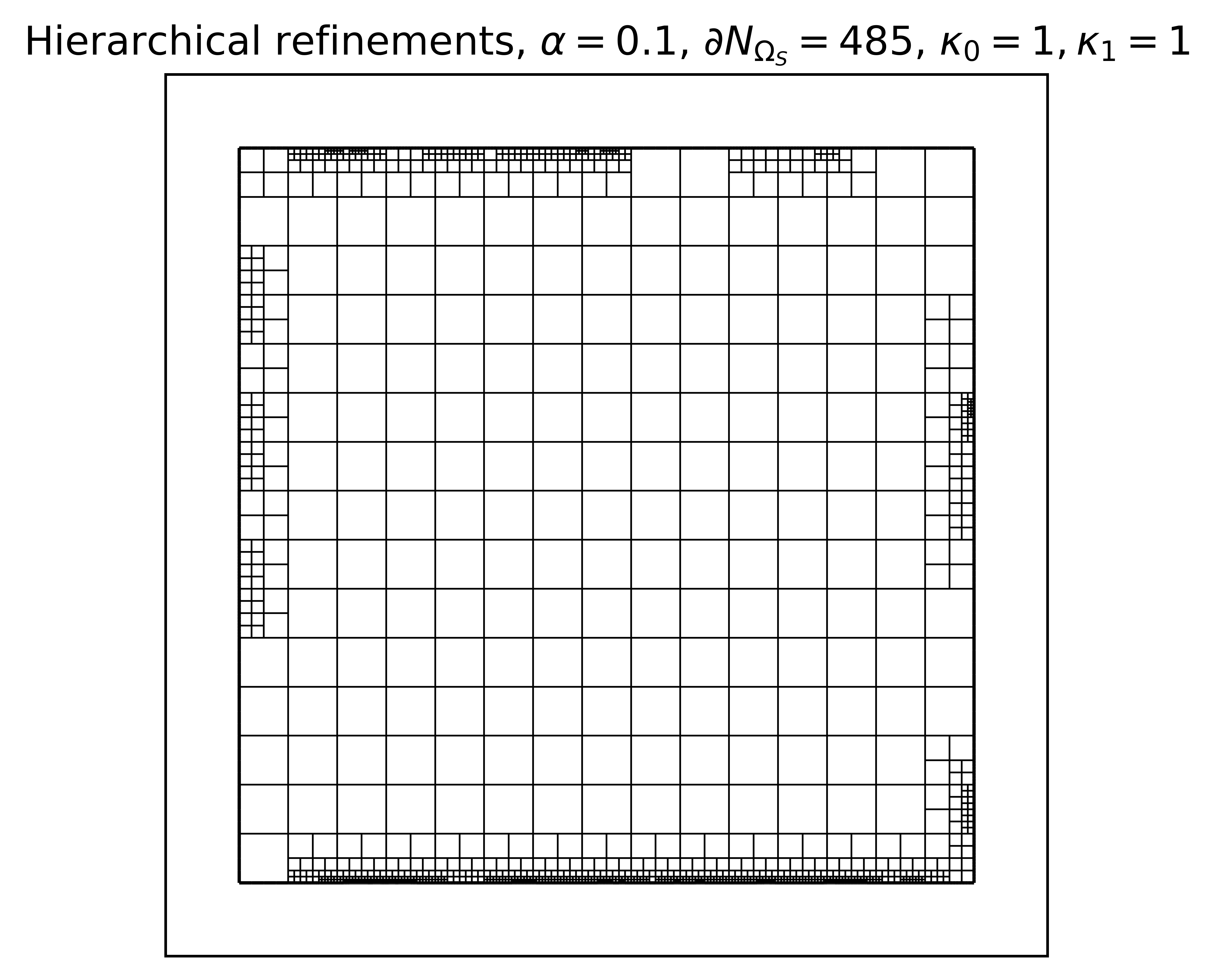} $\quad$ \includegraphics[align=c, width=0.5\linewidth]{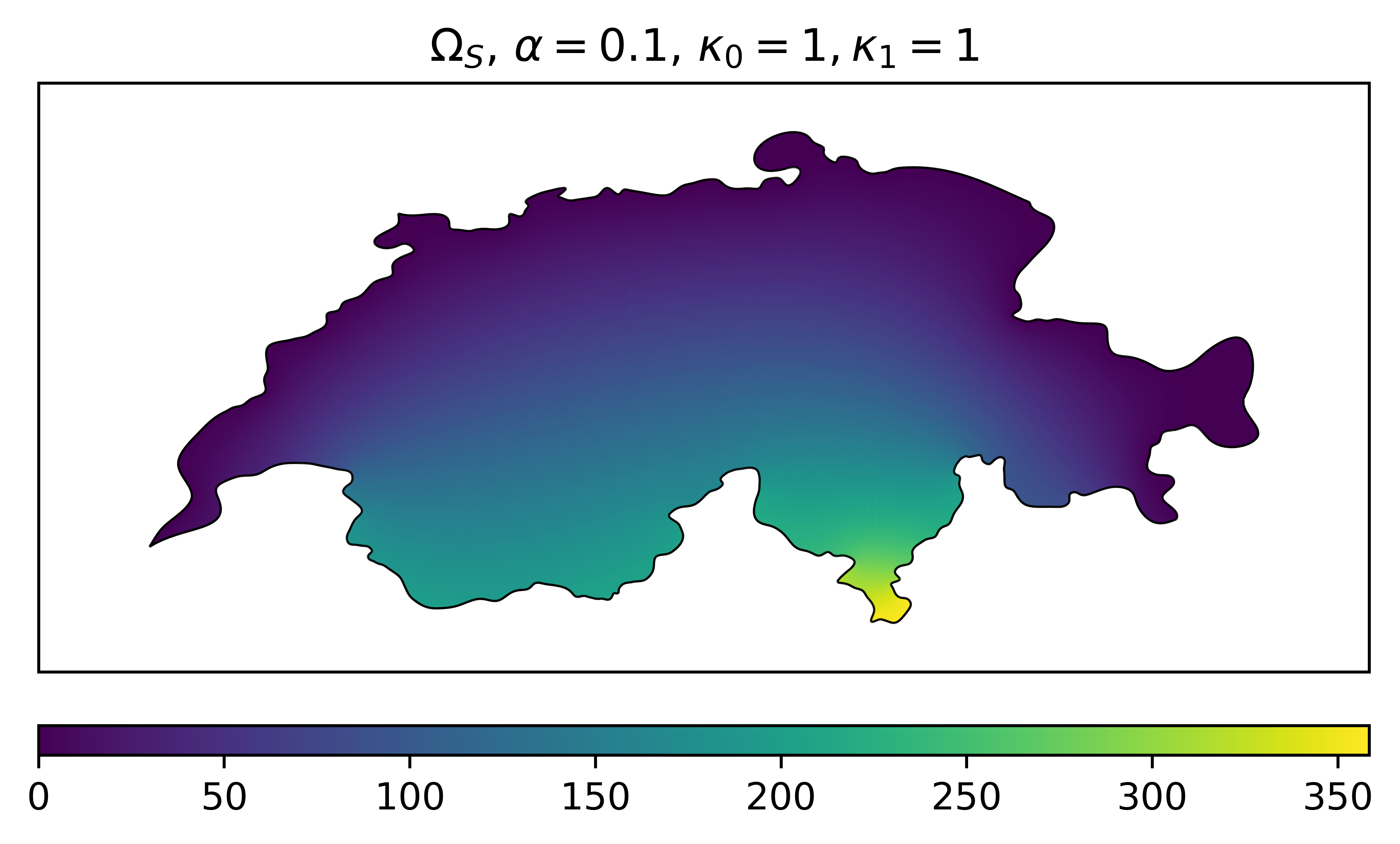}
    \end{subfigure}
\caption{Final iteration's boundary refinements in $\widehat{\Omega}$ (left) and $u_h^n$ plotted on top of the defeatured geometry $\Omega_S$ for $\alpha=0.1$ and $\kappa_0 = \kappa_1 = 1$.}
\label{fig:switz_refinements_geometry_alpha_0.1_H1}
\end{figure}

\begin{figure}[h!]
\centering
    \includegraphics[width=0.8 \linewidth]{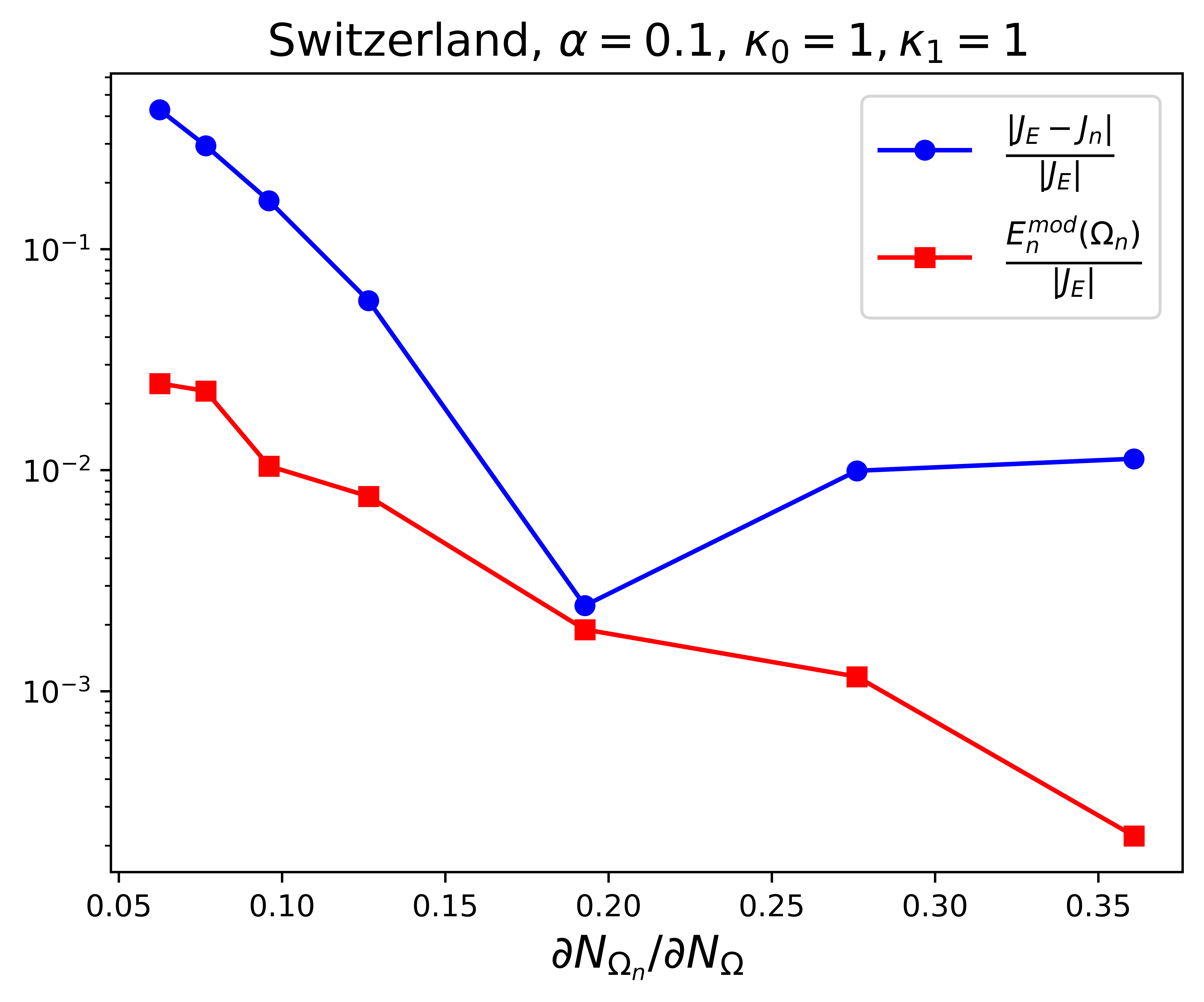}
    \caption{Convergence plot for $\alpha = 0.1$, $\kappa_0 = \kappa_1 = 1$ and convergence threshold $E^{\text{mod}}_n(\Omega_n) < 1$}
\label{fig:switz_convergence_alpha_0.1_H1}
\end{figure}

\noindent To compute each iteration's $(u_h^n, p_h^n)$ sufficiently accurately, we employ an adaptive refinement strategy, see \cite{reviewadaptiveiga}. 
The refinements in $\widehat{\Omega}$ along with $u_h^n$ at the final iteration plotted on $\Omega_S$ are depicted in Figure \ref{fig:switz_refinements_geometry_alpha_0.1_H1}, while Figure \ref{fig:switz_convergence_alpha_0.1_H1} shows the corresponding convergence plot.\\

\noindent First of all, we can see that the final geometry $\Omega_S$ is defeatured from $\Omega$ everywhere at the boundary, but more refinements are performed in the southern part. This is expected as this is the region where the solution gradient is larger, and thus where defeaturing seems to have a larger impact on the solution accuracy, see \cite{buffa2022analysis} for instance. The convergence plot reveals that a relative cost function discrepancy below $10^{-2}$ is achieved for only $20 \%$ of the boundary DOFs after the fifth iteration, which is a remarkable decrease of boundary degrees of freedom while keeping the modelling error small. 

\begin{figure}[h!]
\centering
\captionsetup[subfigure]{labelformat=empty}
    \begin{subfigure}[t]{0.9\textwidth}
        \centering
        \includegraphics[align=c, width=0.35\linewidth]{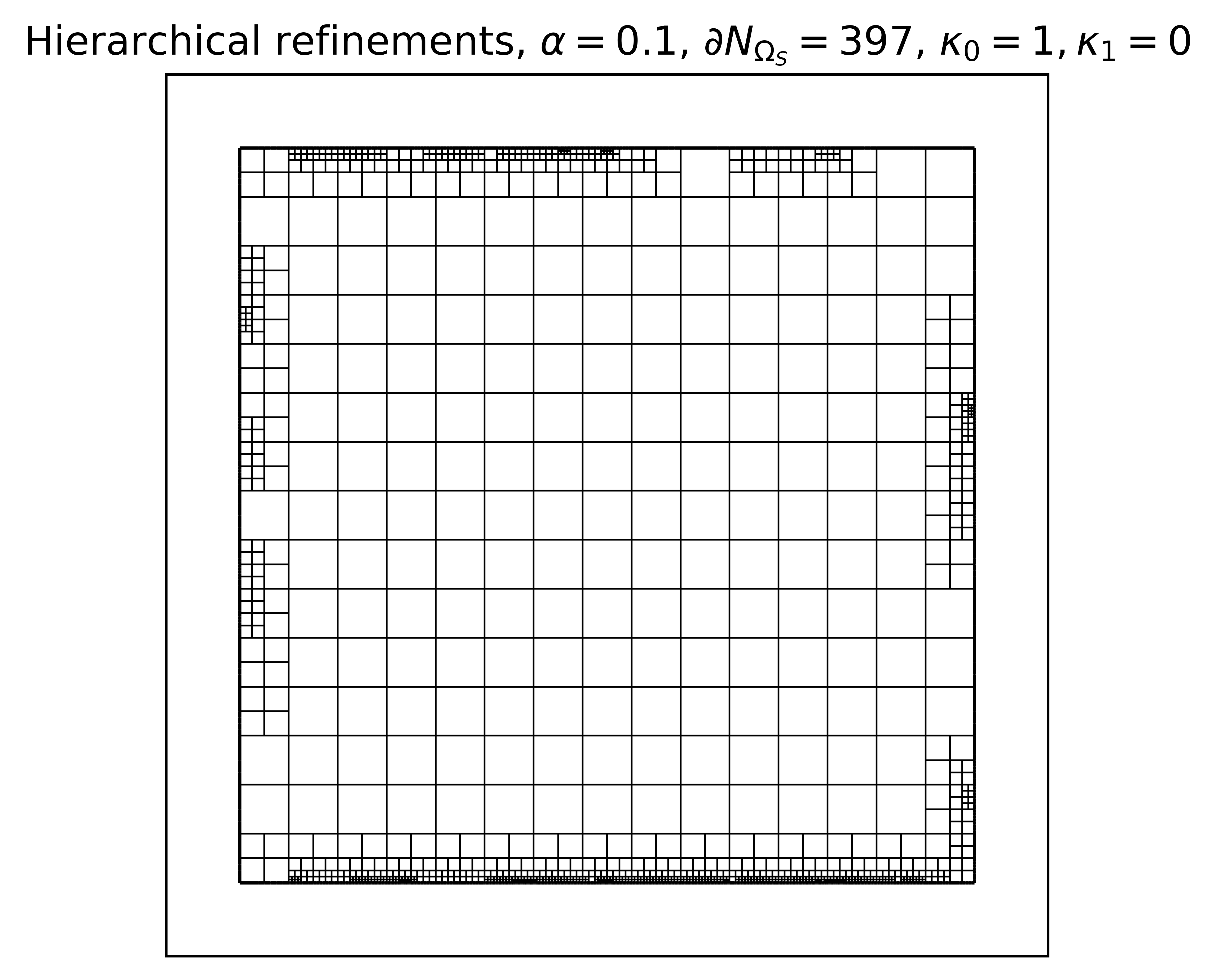} $\quad$ \includegraphics[align=c, width=0.5\linewidth]{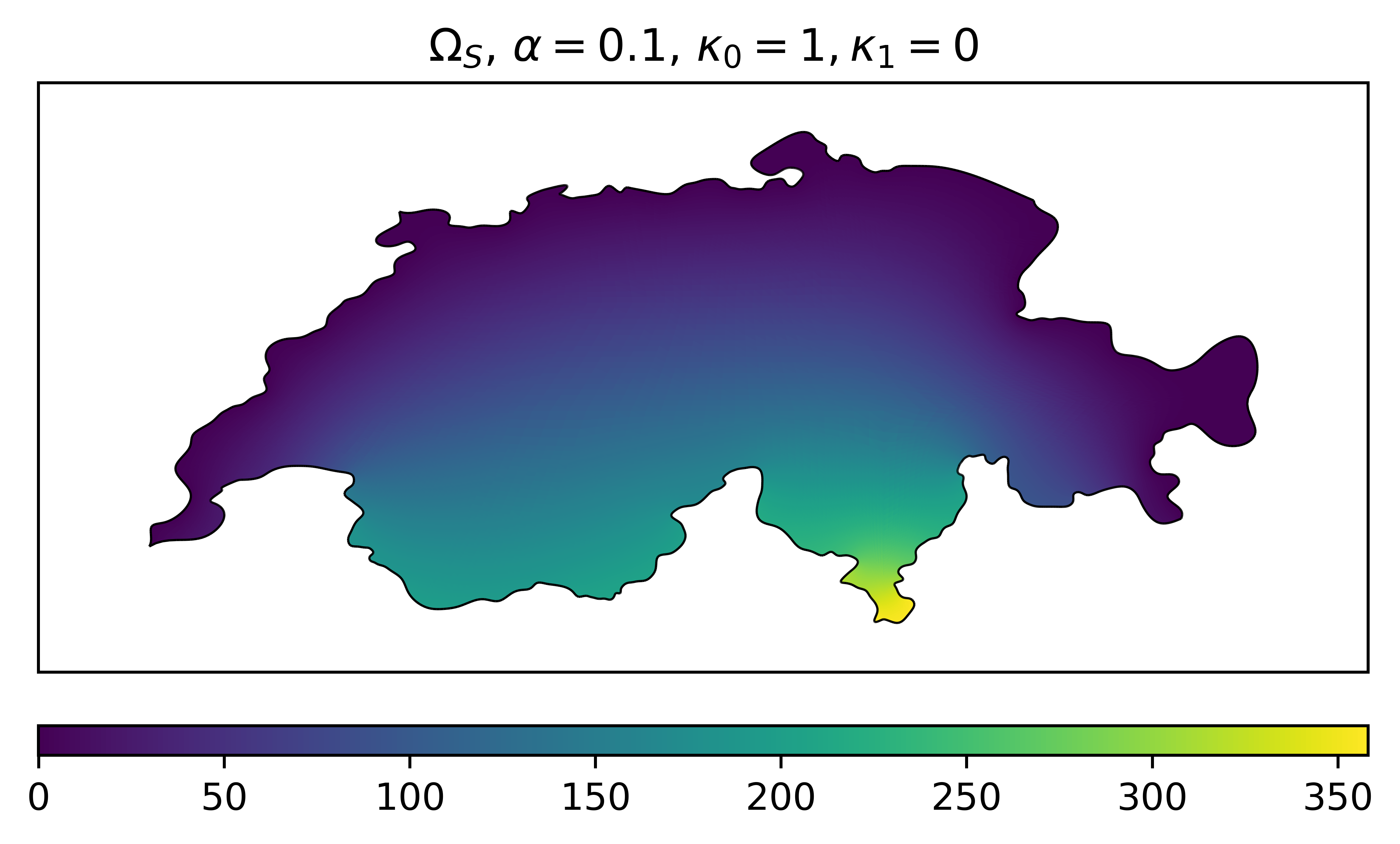}
    \end{subfigure}
\caption{Final iteration's boundary refinements in $\widehat{\Omega}$ (left) and $u_h^n$ plotted on top of the defeatured geometry $\Omega_S$ for $\alpha=0.1$, $\kappa_0 = 1$ and $\kappa_1 = 0$.}
\label{fig:switz_refinements_geometry_alpha_0.1_L2}
\end{figure}

\begin{figure}[h!]
\centering
    \includegraphics[width=0.8 \linewidth]{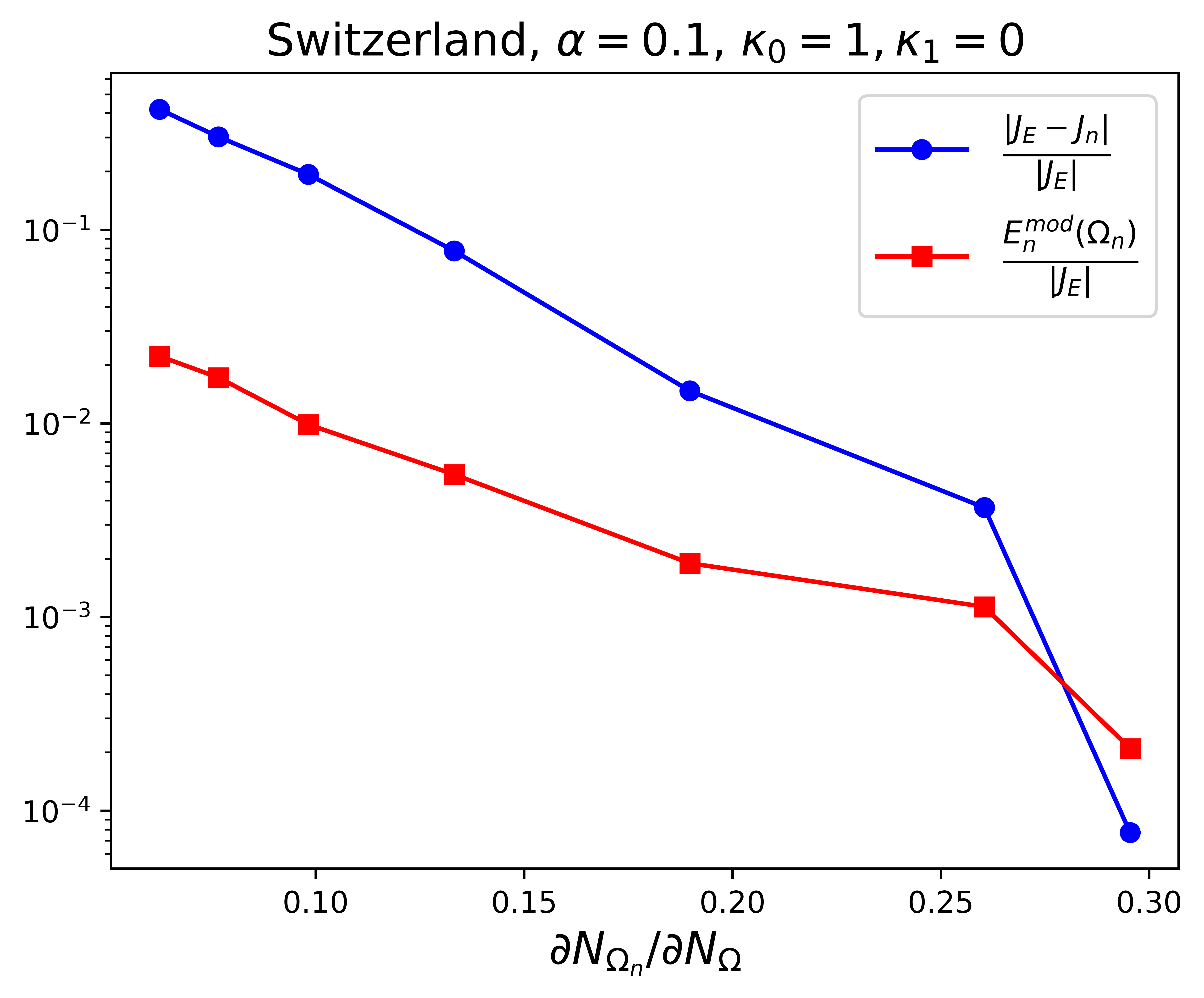}
    \caption{Convergence plot for $\alpha = 0.1$, $\kappa_0 = 1$, $\kappa_1 = 0$ and convergence threshold $E^{\text{mod}}_n(\Omega_n) < 1$}
\label{fig:switz_convergence_alpha_0.1_L2}
\end{figure}

However, the relative error goes back up upon further refinement, then settling for a value slightly above $10^{-2}$. This behavior is not changed if $(u_h^n, p_h^n)$ are computed on finer numerical meshes $\mathcal T_n^u$ and $\mathcal T_n^p$ (without changing the mesh $\mathcal T_n^\mathbf{x}$ that generates the geometry). It therefore suggests that the numerical error is not the culprit. As discussed in the previous numerical experiment, the monotonicity of the modelling error cannot be expected in this algorithm because of the chosen descent directions. Nevertheless, a way to circumvent this would be to try to make the fitting strategy of Section~\ref{subsec:algo_fit} more local. To do so, we repeat the above steps, but with $(\kappa_0, \kappa_1) = (1, 0)$. This yields the convergence plot depicted in Figure \ref{fig:switz_convergence_alpha_0.1_L2} with the corresponding $(\Omega_S, u_h^n)$ at the final iteration depicted in Figure \ref{fig:switz_refinements_geometry_alpha_0.1_L2}. 

The monotone behaviour is restored, suggesting that the fitting parameters have a significant impact on the convergence behaviour. Indeed, increasing the weight of the $L^2(\partial \widehat{\Omega})$ with respect to that of the $H^1(\partial \widehat{\Omega})$ semi-norm in~\eqref{eq:boundary_fit} produces fits that react more drastically to local refinements. This suggests that the choice $\kappa_0 = \kappa_1 = 1$ may have caused the fits to be too non-local in later iterations. % This is further evidenced by the observation that we can estimate (cf. \eqref{eq:discrepancy_J_Omega_plus_Omega_approx}) $|\mathcal{J}(\Omega_5^+) - \mathcal{J}(\Omega_5)| \approx 1.48$, while % (cf. \eqref{eq:discrepancy_J_Omega_next_Omega_approx}) 
% $|\mathcal{J}(\Omega_6) - \mathcal{J}(\Omega_5)| \approx 10.64$. In other words, full refinement of the boundary suggests that a reduction in the cost function can be expected, while a similar estimate suggests that the partially refined boundary increases the cost function in the next iteration. 
% We conclude that the non-local nature of the $(\kappa_0, \kappa_1) = (1, 1)$ fit causes $\Delta \mathbf{F}_5 = \mathbf{F}_6 - \mathbf{F}_5$ to no longer be a descent direction. 
Contrary to this, in the simulation with $(\kappa_0, \kappa_1) = (1, 0)$, the two estimates are always within ten percent of one another. 


%%%%%%%%%%%%%%%%%%%%%%%%%%%%%%%%%%%%%%%%%%%%%%%%%%%%%%%%%%%%%%%%%%%%%

\section{Conclusions} \label{sec:conclusions}
We have successfully conceptualised and implemented an analysis-aware defeaturing algorithm based on the concept of shape-derivatives that exploits the adaptive nature of the THB-spline technology. Here, the algorithm completely avoids the costly parameterisation of the original geometry's interior while only computing parameterisations of intermediate defeatured geometries based on the concept of harmonic maps. The algorithm was able to largely avoid any form of manual intervention while operating on elliptic problems posed over domains $\Omega \subset \mathbb{R}^2$ with varying characteristics.
We have presented two numerical experiments whose results meet our expectations. The numerical experiments furthermore clearly demonstrate that defeaturing can result in substantial (geometrical) DOF savings when one is interested in computing an accurate objective functional.

However in the numerical experiments, we were able to identify one major automation pitfall, namely the possibility of convergence stagnation due to the choice of descent directions and to the nonlocality of the boundary fit~\eqref{eq:boundary_fit}. Some further research should be conducted to avoid the manual fine tuning of the marking parameter $\alpha$ and of the fitting parameters $(\kappa_0, \kappa_1)$. One could for instance introduce a more local fitting strategy making the most of quasi-interpolants for instance, see \cite{leelychemorken}. However, erratic behaviour in the boundary correspondence $\mathbf{F}_n: \partial \widehat{\Omega} \rightarrow \partial \Omega_n$ caused by drastic reactions to local refinements increases the chances of self-intersecting boundary contours $\partial \Omega_n$ as well as generally requiring more aPos refinements in the computation of the harmonic map. A possible remedy may be introducing position dependent parameters, i.e., $(\kappa_0, \kappa_1) \rightarrow (\kappa_0(\bm{\xi}), \kappa_1(\bm{\xi}))$, even though the automation of this choice remains unclear. We therefore expect the algorithm's robustness to greatly benefit from a more thorough investigation of the effect that various contour approximation / fitting procedures have on the convergence behaviour.

% More precisely, the boundary fit may cause $\Delta \mathbf{F}_n = \mathbf{F}_{n+1} - \mathbf{F}_n$ to no longer be a descent direction for certain choices of $(\kappa_0, \kappa_1)$. In our example, we were able to restore monotone convergence by changing the values of the fitting parameters $(\kappa_0, \kappa_1)$, whereby lower values of $\kappa_1$ generally produce fits that react more drastically to local refinements. possibly keeping in mind that, upon re-fitting, $\Delta \mathbf{F}_n$ should be a descent direction.

%\appendix
%\section{Appendix} \label{sec:appendix}
%%%%%%%%%%%%%%%%%%%%%%%%%%%%%%%%%%%%%%%%%%%%%%%%%%%%%%%%%%%%%%%%%%%%%

\section*{Acknowledgments}
The authors gratefully acknowledge the support of the Swiss National Science Foundation through the project ‘‘Design-through-Analysis (of PDEs): the litmus test’’ n. 40B2-0 187094 (BRIDGE Discovery 2019). 
Ondine Chanon also thanks the support of the Swiss National Science Foundation via the project n.P500PT\_210974.

%Référence et bibliographie :
%\newpage{\tiny {\tiny }}
%\nocite{*}
\addcontentsline{toc}{section}{Bibliography}
\bibliography{defeaturing_shape_der_bib}

\begin{thebibliography}{10}

\bibitem{cottrell2009isogeometric}
J.~A. Cottrell, T.~J.~R. Hughes, and Y.~Bazilevs, {\em Isogeometric analysis:
  toward integration of {CAD} and {FEA}}.
\newblock John Wiley \& Sons, 2009.

\bibitem{hughes2005isogeometric}
T.~J.~R. Hughes, J.~A. Cottrell, and Y.~Bazilevs, ``Isogeometric analysis:
  {CAD}, finite elements, {NURBS}, exact geometry and mesh refinement,'' {\em
  cmame}, vol.~194, pp.~4135--4195, 2005.

\bibitem{white2003meshing}
D.~R. White, S.~Saigal, and S.~J. Owen, ``Meshing complexity of single part
  {CAD} models.,'' in {\em IMR}, pp.~121--134, 2003.

\bibitem{trimmingbletzinger}
R.~Schmidt, R.~Wüchner, and K.-U. Bletzinger, ``Isogeometric analysis of
  trimmed {NURBS} geometries,'' {\em Computer Methods in Applied Mechanics and
  Engineering}, vol.~241-244, pp.~93--111, 2012.

\bibitem{trimming}
B.~Marussig and T.~J.~R. Hughes, ``A review of trimming in isogeometric
  analysis: challenges, data exchange and simulation aspects,'' {\em Archives
  of Computational Methods in Engineering}, vol.~25, pp.~1059--1127, Nov 2018.

\bibitem{antolinvreps}
P.~Antol\'in, A.~Buffa, and M.~Martinelli, ``Isogeometric analysis on {V}-reps:
  first results,'' {\em Computer Methods in Applied Mechanics and Engineering},
  vol.~355, pp.~976--1002, 2019.

\bibitem{weimarrusigantolin}
X.~Wei, B.~Marussig, P.~Antol\'in, and A.~Buffa, ``Immersed boundary-conformal
  isogeometric method for linear elliptic problems,'' {\em Computational
  Mechanics}, vol.~68, no.~6, pp.~1385--1405, 2021.

\bibitem{antolin2021quadrature}
P.~Antol\'in and T.~Hirschler, ``Quadrature-free immersed isogeometric
  analysis,'' {\em Engineering with Computers}, 2022.

\bibitem{unionkargaran}
S.~Kargaran, B.~Jüttler, S.~K. Kleiss, A.~Mantzaflaris, and T.~Takacs,
  ``Overlapping multi-patch structures in isogeometric analysis,'' {\em
  Computer Methods in Applied Mechanics and Engineering}, vol.~356,
  pp.~325--353, 2019.

\bibitem{unionzuo}
B.-Q. Zuo, Z.-D. Huang, Y.-W. Wang, and Z.-J. Wu, ``Isogeometric analysis for
  {CSG} models,'' {\em Computer Methods in Applied Mechanics and Engineering},
  vol.~285, pp.~102--124, 2015.

\bibitem{antolinwei}
P.~Antol\'in, A.~Buffa, R.~Puppi, and X.~Wei, ``Overlapping multipatch
  isogeometric method with minimal stabilization,'' {\em SIAM Journal on
  Scientific Computing}, vol.~43, no.~1, pp.~A330--A354, 2021.

\bibitem{multipatch}
A.~Buffa, R.~Vázquez, G.~Sangalli, and L.~Beirão~da Veiga, ``Approximation
  estimates for isogeometric spaces in multipatch geometries,'' {\em Numerical
  Methods for Partial Differential Equations}, vol.~31, no.~2, pp.~422--438,
  2015.

\bibitem{bracco2020isogeometric}
C.~Bracco, C.~Giannelli, M.~Kapl, and R.~V{\'a}zquez, ``Isogeometric analysis
  with ${C}^1$-hierarchical functions on planar two-patch geometries,'' {\em
  Computers \& Mathematics with Applications}, vol.~80, no.~11, pp.~2538--2562,
  2020.

\bibitem{overview_thb}
C.~Bracco, A.~Buffa, C.~Giannelli, and R.~V\'{a}zquez, ``Adaptive isogeometric
  methods with hierarchical splines: an overview,'' {\em Discrete Contin. Dyn.
  Syst.}, vol.~39, no.~1, pp.~241--261, 2019.

\bibitem{thakur2009survey}
A.~Thakur, A.~G. Banerjee, and S.~K. Gupta, ``A survey of {CAD} model
  simplification techniques for physics-based simulation applications,'' {\em
  Computer-Aided Design}, vol.~41, no.~2, pp.~65--80, 2009.

\bibitem{fine2000automated}
L.~Fine, L.~Remondini, and J.-C. Leon, ``Automated generation of {FEA} models
  through idealization operators,'' {\em International Journal for Numerical
  Methods in Engineering}, vol.~49, no.~1-2, pp.~83--108, 2000.

\bibitem{foucault2004mechanical}
G.~Foucault, P.~M. Marin, and J.-C. L{\'e}on, ``Mechanical criteria for the
  preparation of finite element models.,'' in {\em IMR}, pp.~413--426, 2004.

\bibitem{rahimi2018cad}
N.~Rahimi, P.~Kerfriden, F.~C. Langbein, and R.~R. Martin, ``{CAD} model
  simplification error estimation for electrostatics problems,'' {\em SIAM
  Journal on Scientific Computing}, vol.~40, no.~1, pp.~B196--B227, 2018.

\bibitem{ferrandes2009posteriori}
R.~Ferrandes, P.~Marin, J.-C. L{\'e}on, and F.~Giannini, ``A posteriori
  evaluation of simplification details for finite element model preparation,''
  {\em Computers \& Structures}, vol.~87, no.~1-2, pp.~73--80, 2009.

\bibitem{sokolowski1999topological}
J.~Sokolowski and A.~Zochowski, ``On the topological derivative in shape
  optimization,'' {\em SIAM journal on control and optimization}, vol.~37,
  no.~4, pp.~1251--1272, 1999.

\bibitem{choi2004structural}
K.~K. Choi and N.-H. Kim, {\em Structural sensitivity analysis and optimization
  1: linear systems}.
\newblock Springer Science \& Business Media, 2004.

\bibitem{li2011estimating}
M.~Li, S.~Gao, and R.~R. Martin, ``Estimating the effects of removing negative
  features on engineering analysis,'' {\em Computer-Aided Design}, vol.~43,
  no.~11, pp.~1402--1412, 2011.

\bibitem{becker2001optimal}
R.~Becker and R.~Rannacher, ``An optimal control approach to a posteriori error
  estimation in finite element methods,'' {\em Acta numerica}, vol.~10,
  pp.~1--102, 2001.

\bibitem{oden2002estimation}
J.~T. Oden and S.~Prudhomme, ``Estimation of modeling error in computational
  mechanics,'' {\em Journal of Computational Physics}, vol.~182, no.~2,
  pp.~496--515, 2002.

\bibitem{li2013engineering}
M.~Li, S.~Gao, and R.~R. Martin, ``Engineering analysis error estimation when
  removing finite-sized features in nonlinear elliptic problems,'' {\em
  Computer-Aided Design}, vol.~45, no.~2, pp.~361--372, 2013.

\bibitem{li2013goal}
M.~Li, S.~Gao, and K.~Zhang, ``A goal-oriented error estimator for the analysis
  of simplified designs,'' {\em Computer Methods in Applied Mechanics and
  Engineering}, vol.~255, pp.~89--103, 2013.

\bibitem{zhang2016estimation}
K.~Zhang, M.~Li, and J.~Li, ``Estimation of impacts of removing arbitrarily
  constrained domain details to the analysis of incompressible fluid flows,''
  {\em Communications in Computational Physics}, vol.~20, no.~4, pp.~944--968,
  2016.

\bibitem{buffa2022analysis}
A.~Buffa, O.~Chanon, and R.~V{\'a}zquez, ``Analysis-aware defeaturing: problem
  setting and a posteriori estimation,'' {\em Mathematical Models and Methods
  in Applied Sciences}, vol.~32, no.~02, pp.~359--402, 2022.

\bibitem{antolin2022multifeature}
P.~Antol{\'i}n and O.~Chanon, ``Analysis-aware defeaturing of complex
  geometries,'' {\em arXiv preprint arXiv:2212.03141}, 2022.

\bibitem{buffa2022adaptive}
A.~Buffa, O.~Chanon, and R.~V{\'a}zquez, ``Adaptive analysis-aware
  defeaturing,'' {\em arXiv preprint arXiv:2212.05183}, 2022.

\bibitem{buffa2021adaptiveapproxshapes}
A.~Buffa, R.~Hiptmair, and P.~Panchal, ``Adaptive approximation of shapes,''
  {\em Numerical Functional Analysis and Optimization}, vol.~42, no.~2,
  pp.~132--154, 2021.

\bibitem{heydarov2022unrefinement}
T.~Heydarov, A.~Buffa, and B.~J{\"u}ttler, ``An unrefinement algorithm for
  planar {THB}-spline parameterizations,'' {\em Computer Aided Geometric
  Design}, vol.~99, p.~102157, 2022.

\bibitem{introshapeopti_zol}
J.~Soko\l~\!\!owski and J.-P. Zol\'{e}sio, {\em Introduction to shape
  optimization}, vol.~16 of {\em Springer Series in Computational Mathematics}.
\newblock Springer-Verlag, Berlin, 1992.
\newblock Shape sensitivity analysis.

\bibitem{shape&geom_zol_delf}
M.~C. Delfour and J.-P. Zol\'{e}sio, {\em Shapes and geometries}, vol.~22 of
  {\em Advances in Design and Control}.
\newblock Society for Industrial and Applied Mathematics (SIAM), Philadelphia,
  PA, second~ed., 2011.
\newblock Metrics, analysis, differential calculus, and optimization.

\bibitem{article_paganini}
R.~Hiptmair, A.~Paganini, and S.~Sargheini, ``Comparison of approximate shape
  gradients,'' {\em BIT}, vol.~55, no.~2, pp.~459--485, 2015.

\bibitem{Verani_adaptiveFEM}
P.~Morin, R.~H. Nochetto, M.~S. Pauletti, and M.~Verani, ``Adaptive finite
  element method for shape optimization,'' {\em ESAIM Control Optim. Calc.
  Var.}, vol.~18, no.~4, pp.~1122--1149, 2012.

\bibitem{Verani_discreteGradFlows}
G.~Do\v{g}an, P.~Morin, R.~H. Nochetto, and M.~Verani, ``Discrete gradient
  flows for shape optimization and applications,'' {\em Comput. Methods Appl.
  Mech. Engrg.}, vol.~196, no.~37-40, pp.~3898--3914, 2007.

\bibitem{Henrot_bookshapeder_clear}
A.~Henrot and M.~Pierre, {\em Shape variation and optimization}, vol.~28 of
  {\em EMS Tracts in Mathematics}.
\newblock European Mathematical Society (EMS), Z\"{u}rich, 2018.

\bibitem{geopde_hier}
E.~M. Garau and R.~V\'{a}zquez, ``Algorithms for the implementation of adaptive
  isogeometric methods using hierarchical {B}-splines,'' {\em Appl. Numer.
  Math.}, vol.~123, pp.~58--87, 2018.

\bibitem{buffa2017refinable}
A.~Buffa and E.~M. Garau, ``Refinable spaces and local approximation estimates
  for hierarchical splines,'' {\em IMA Journal of Numerical Analysis}, vol.~37,
  no.~3, pp.~1125--1149, 2017.

\bibitem{thb_giannelli_juttler}
C.~Giannelli, B.~J\"{u}ttler, and H.~Speleers, ``T{HB}-splines: the truncated
  basis for hierarchical splines,'' {\em Comput. Aided Geom. Design}, vol.~29,
  no.~7, pp.~485--498, 2012.

\bibitem{knupp2020fundamentals}
P.~Knupp and S.~Steinberg, {\em Fundamentals of grid generation}.
\newblock CRC press, 2020.

\bibitem{hinz2020goal}
J.~Hinz, M.~Abdelmalik, and M.~M{\"o}ller, ``Goal-oriented adaptive
  {THB}-spline schemes for {PDE}-based planar parameterization,'' {\em arXiv
  preprint arXiv:2001.08874}, 2020.

\bibitem{gravesen2012planar}
J.~Gravesen, A.~Evgrafov, D.-M. Nguyen, and P.~N{\o}rtoft, ``Planar
  parametrization in isogeometric analysis,'' in {\em International conference
  on mathematical methods for curves and surfaces}, pp.~189--212, Springer,
  2012.

\bibitem{choquet1945type}
G.~Choquet, ``Sur un type de transformation analytique g{\'e}n{\'e}ralisant la
  repr{\'e}sentation conforme et d{\'e}finie au moyen de fonctions
  harmoniques,'' {\em Bull. Sci. Math.}, vol.~69, no.~2, pp.~156--165, 1945.

\bibitem{lakkis2011finite}
O.~Lakkis and T.~Pryer, ``A finite element method for second order
  nonvariational elliptic problems,'' {\em SIAM Journal on Scientific
  Computing}, vol.~33, no.~2, pp.~786--801, 2011.

\bibitem{gallistl2017variational}
D.~Gallistl, ``Variational formulation and numerical analysis of linear
  elliptic equations in nondivergence form with {C}ßordes coefficients,'' {\em
  SIAM Journal on Numerical Analysis}, vol.~55, no.~2, pp.~737--757, 2017.

\bibitem{blechschmidt2019error}
J.~Blechschmidt, R.~Herzog, and M.~Winkler, ``Error estimation for second-order
  {PDE}s in non-variational form,'' {\em arXiv preprint arXiv:1909.12676},
  2019.

\bibitem{nutils7}
J.~van Zwieten, G.~van Zwieten, and W.~Hoitinga, ``Nutils 7.0,'' 2022.

\bibitem{reviewadaptiveiga}
A.~Buffa, G.~Gantner, C.~Giannelli, D.~Praetorius, and R.~V{\'a}zquez,
  ``Mathematical foundations of adaptive isogeometric analysis,'' {\em Archives
  of Computational Methods in Engineering}, vol.~29, no.~7, pp.~4479--4555,
  2022.

\bibitem{leelychemorken}
B.-G. Lee, T.~Lyche, and K.~M{\o}rken, ``Some examples of quasi-interpolants
  constructed from local spline projectors,'' {\em Mathematical methods for
  curves and surfaces: Oslo}, pp.~243--252, 2000.

\end{thebibliography}
\bibliographystyle{ieeetr}
%\bibliographystyle{unsrt}

%\bibliographystyle{model1-num-names}
%\bibliography{sample.bib}

%% Authors are advised to submit their bibtex database files. They are
%% requested to list a bibtex style file in the manuscript if they do
%% not want to use model1-num-names.bst.

%% References without bibTeX database:

% \begin{thebibliography}{00}

%% \bibitem must have the following form:
%%   \bibitem{key}...
%%

% \bibitem{}

% \end{thebibliography}

\end{document}